%&amstex          
\input amstex\documentstyle{amsppt}  
\pagewidth{12.5cm}\pageheight{19cm}\magnification\magstep1
\topmatter
\title Distinguished conjugacy classes and elliptic Weyl group elements\endtitle
\author G. Lusztig\endauthor
\rightheadtext{Distinguished classes and elliptic Weyl group elements}
\address{Department of Mathematics, M.I.T., Cambridge, MA 02139}\endaddress
\thanks{Supported in part by National Science Foundation grant DMS-0758262.}\endthanks
\endtopmatter   
\document

\define\sneq{\subsetneqq}

\define\ul{\un l}

\define\uW{\un W}

\define\mpb{\medpagebreak}

\define\pe{\perp}
\define\si{\sim}

\define\sqc{\sqcup}
\define\ovs{\overset}
\define\qua{\quad}

\define\bin{\binom}
\define\op{\oplus}
   
\define\part{\partial}
\define\emp{\emptyset}
\define\imp{\implies}

\define\n{\notin}

\define\m{\mapsto}
\define\do{\dots}

\define\sub{\subset}    

\define\T{\times}
\define\ti{\tilde}
\define\nl{\newline}
\redefine\i{^{-1}}

\define\un{\underline}

\define\Aut{\text{\rm Aut}}

\define\clu{\clubsuit}

\define\a{\alpha}
\redefine\b{\beta}

\define\g{\gamma}
\redefine\d{\delta}
\define\e{\epsilon}

\redefine\o{\omega}
\define\p{\pi}
\define\ph{\phi}
\define\ps{\psi}

\define\s{\sigma}
\redefine\t{\tau}
\define\th{\theta}
\define\k{\kappa}

\define\z{\zeta}
\define\x{\xi}

\redefine\G{\Gamma}
\redefine\D{\Delta}
\define\Om{\Omega}

\define\Ph{\Phi}

\define\boc{\bold c}

\define\kk{\bold k}

\define\nn{\bold n}

\define\NN{\bold N}

\define\ZZ{\bold Z}

\define\cb{\Cal B}
\define\cc{\Cal C}

\define\ci{\Cal I}

\define\cl{\Cal L}

\define\co{\Cal O}

\define\cz{\Cal Z}

\define\fB{\frak B}

\define\tg{\ti g}

\define\tv{\ti v}
\define\tw{\ti w}
\define\tz{\ti z}

\define\tcc{\ti\cc}

\define\bc{\bar c}

\define\che{\check}

\define\chg{\che g}

\define\BC{BC}
\define\GP{GP}
\define\WEU{L1}
\define\WES{L2}
\define\WEH{L3}
\define\WEII{L4}
\define\WEIII{L5}
\define\WEC{L6}

\head Introduction \endhead
\subhead 0.1\endsubhead
Let $\kk$ be an algebraically closed field of characteristic $p\ge0$
and let $G$ be a (possibly disconnected) reductive algebraic
group over $\kk$. Let $W$ be the Weyl group of $G^0$. (For an algebraic
group $H$, $H^0$ denotes the identity component of $H$.) We view $W$ as an indexing set for the orbits of 
$G^0$ acting diagonally on $\cb\T\cb$ where $\cb$ is the variety of Borel subgroups of $G^0$; we denote by 
$\co_w$ the orbit corresponding to $w\in W$. Note that $W$ is naturally a Coxeter group; its length function
is denoted by $\ul:W@>>>\NN$. 
Let $I$ be the set of simple reflections of $W$; for any $J\sub I$ let $W_J$ be the subgroup of $W$ generated
by $J$. 

Now any $\D\in G/G^0$ defines a group automorphism $\e_\D:W@>>>W$ preserving length, by the requirement that 
$$(B,B')\in\co_w,g\in\D\imp(gBg\i,gB'g\i)\in\co_{\e_\D(w)}.$$
The orbits of the $W$-action $w_1:w\m w_1\i w\e_\D(w_1)$ on $W$ are said to be the $\e_\D$-conjugacy classes
in $W$. Let $\uW_\D$ be the set of $\e_D$-conjugacy classes in $W$. 
We say that $C\in\uW_\D$ is elliptic if for any $J\sneq I$ such that $\e_D(J)=J$ we have $C\cap W_J=\emp$.
Let $\uW^{el}_\D$ be the set of elliptic $\e_D$-conjugacy classes in $W$. For any $C\in\uW_\D$ let 
$C_{min}$ be the set of elements of $C$ where the length function $\ul:C@>>>\NN$ reaches its 
minimum value. Let $\boc$ be a $G^0$-conjugacy class of $G$. Let $\D$ be the connected component of $G$
that contains $\boc$ and let $C\in\uW^{el}_\D$. For any $w\in C_{min}$ we set
$$\fB^\boc_w=\{(g,B)\in\boc\T\cb;(B,gBg\i)\in\co_w\}.$$ 
Note that $G^0$ acts on $\fB^\boc_w$ by $x:(g,B)\m(xgx\i,xBx\i)$. 
We write $C\clu\boc$ if the following condition is satisfied: for some/any $w\in C_{min}$, $\fB^\boc_w$ is a
single $G^0$-orbit for the action above (in particular it is nonempty). The equivalence of ``some'' and 
``any'' follows from \cite{\WEIII, 1.15(a)}. 

\subhead 0.2\endsubhead
For an algebraic group $H$ we denote by $\cz_H$ the centre of $H$; for $h\in H$ we denote by $Z_H(h)$ the 
centralizer of $h$ in $H$. An element $g\in G$ or its $G^0$-conjugacy class is said to be 
{\it distinguished} if $Z_G(g)^0/(\cz_{G^0}\cap Z_G(g))^0$ is a unipotent group. 
The notion of distinguished element appeared in \cite{\BC} in the case where $g$ is unipotent and $G=G^0$.

The following is the main result of this paper.
\proclaim{Theorem 0.3} 
Assume that $G^0$ is almost simple and that $|G/G^0|\le2$. If $G^0$ is of exceptional type assume further
that $G=G^0$ and that $p$ is either $0$ or a good prime for $G$. Then for any distinguished $G^0$-conjugacy 
class $\boc$ in $G$ contained in a connected component $\D$ of $G$, there exists $C\in\uW^{el}_\D$ such that
$C\clu\boc$.
\endproclaim
In the case where $\boc$ is unipotent the theorem is known from \cite{\WEH}. In particular the theorem holds
when $p=2$. Thus we may assume that $p\ne2$. We 
may also assume that $G/G^0@>>>\Aut(W)$, $\D\m\e_\D$ is injective. It is enough to verify the theorem 
assuming that $G^0$ is simply connected (the theorem then automatically holds without that assumption). If 
$G^0$ is of type $A$ and $G=G^0$ then $\boc$ must be a regular unipotent class times a central element and 
we can take $C$ to be the Coxeter class. The case where $G=G^0$ is of type $B$ or $C$ is treated in \S1. The
case where $G^0$ is of type $D$ is treated also in \S1. (In this case we may assume that $|G/G^0|=2$.) The 
case where $G^0$ is of type $A$ and $|G/G^0|=2$ is treated in \S2. (In this case we may assume that 
$\boc\n G^0$.) The case where $G$ is of exceptional type is treated in \S3.

We will show elsewhere that $C$ in the theorem is unique (in the case where $\boc$ is unipotent this is 
known from \cite{\WEU}). 

\subhead 0.4\endsubhead
The results of this paper have applications to the study of character sheaves. We will show elsewhere how 
they can be used to prove that an irreducible cuspidal local system on $\boc$ (a distinguished 
$G^0$-conjugacy class in a connected component $\D$ of $G$), extended by $0$ on $\D-\boc$, is (up to shift) 
a character sheaf on $\D$. In the case where $\D=G^0$ this gives a new, constructive proof of a known 
result, but in the case where $\D\ne G^0$, it is a new result.

\head 1. Isometries\endhead
\subhead 1.0\endsubhead
In this section we assume that $p\ne2$. Let $\e\in\{1,-1\}$. Let $V$ be a 
$\kk$-vector space of finite dimension $\nn$ with a given nondegenerate bilinear form $(,):V\T V@>>>\kk$ 
such that $(x,y)=\e(y,x)$ for all $x,y$; we then say that $(,)$ is $\e$-symmetric. Let $Is(V)$ be the group 
of isometries of $(,)$. 

Assume that we are given $g\in Is(V)$. For any $z\in V$ and $i\in\ZZ$ we set $z_i=g^iz\in V$. Similarly, for 
any line $L$ in $V$ and $i\in\ZZ$ we set $L_i=g^iL\sub V$. For any $z,z'$ in $V$ and any $i,j,k\in\ZZ$ we 
have
$$(z_{i+k},z'_{j+k})=(z_i,z'_j).\tag a$$
Let $a_1\ge a_2\ge\do$, $b_1\ge b_2\ge\do$ be two sequences in $\NN$ such that 

if $i\ge1$, $a_i=a_{i+1}$ then $a_{i+1}=0$,

if $i\ge1$, $b_i=b_{i+1}$ then $b_{i+1}=0$,

if $a_i>0$, then $(-1)^{a_i}=-\e$,

if $b_i>0$, then $(-1)^{b_i}=-\e$.
\nl
It follows that $a_i=0$ for large $i$ and $b_i=0$ for large $i$. We assume that 
$$\nn=(a_1+a_2+\do)+(b_1+b_2+\do).$$
Define $\k\in\{0,1\}$ by $\nn-\k\in2\NN$. Note that if $\e=-1$ we have $\k=0$. Define $k\ge0$ by 
$\{i\ge1;a_ib_i>0\}=[1,k]$. For $i\ge1$ we set $c_i=a_i+b_i$. We have $c_1\ge c_2\ge\do$. We define 
$p_i\in\NN$ for $i\ge1$ as follows.

If $\e=-1$ we have $c_i\in2\NN$ and we set $p_i=c_i/2$ for $i\ge1$.

If $\e=1$ and $i\in[1,k]$ we again have $c_i\in2\NN$ and we set $p_i=c_i/2$.

If $\e=1$ and $i>k$ we have $c_i\in2\NN+1$ or $c_i=0$ and we define $p_i$ by requiring that for 
$s=1,3,5,\do$ we have:

$(p_{k+s},p_{k+s+1})=((c_{k+s}-1)/2,(c_{k+s+1}+1)/2)$ if $c_{k+s}\ge1,c_{k+s+1}\ge1$;

$(p_{k+s},p_{k+s+1})=((c_{k+s}-1)/2,0)$ if $c_{k+s}\ge1,c_{k+s+1}=0$;

$(p_{k+s},p_{k+s+1})=(0,0)$ if $c_{k+s}=0,c_{k+s+1}=0$.
\nl
We define $\s$ as follows. We have $p_1\ge p_2\ge\do\ge p_\s$ where $p_i\in\NN_{>0}$ for $i\in[1,\s]$, 
$p_i=0$ if $i>\s$. This defines $\s$. If $\nn=0$ or $\nn=1$ we have $\s=0$. We set $p'_t=p_t$ if 
$t\in[1,\s]$, $p'_t=1/2$ if $\k=1,t=\s+1$. We have 
$$2(p_1+p_2+\do+p_\s)+\k=2(p'_1+p'_2+\do+p'_{\s+\k})=\nn.$$ 

Let $\cc^V_{a_*,b_*}$ be the set of all $g\in Is(V)$ such that $g^2:V@>>>V$ is unipotent and such that on 
the generalized $1$-eigenspace of $g$, $g$ has Jordan blocks of sizes given by the nonzero numbers in 
$a_1,a_2,\do$ and on the generalized $(-1)$-eigenspace of $g$, $-g$ has Jordan blocks of sizes given by the 
nonzero numbers in $b_1,b_2,\do$. 

For $g\in\cc^V_{a_*,b_*}$ let $\tcc^V_{g;a_*,b_*}$ be the set consisting of all $L^1,L^2,\do,L^{\s+\k}$ 
where $L^t (t\in[1,\s+\k])$ are lines in $V$ (the upper scripts are not powers) such that for 
$i,j\in\ZZ$ we have:

$(L^t_i,L^t_j)=0$ if $|i-j|<p_t$, $(L^t_i,L^t_j)\ne0$ if $j-i=p_t$ ($t\in[1,\s+\k]$);

$(L^t_i,L^r_j)=0$ if $i-j\in[-p_r,2p_t-p_r-1]$ and $1\le t<r\le\s+\k$.
\nl
Here $L^t_i=g^iL^t$. We then have:

(b) $V=\op_{t\in[1,\s+\k],i\in[0,2p'_t-1]}L^t_i$.
\nl
(See \cite{\WEH, 1.3}.) Let $\tcc^V_{a_*,b_*}$ be the set of all $(g,L^1,L^2,\do,L^{\s+\k})$ such that
$g\in\cc^V_{a_*,b_*}$ and $(L^1,L^2,\do,L^{\s+\k})\in\tcc^V_{g;a_*,b_*}$.

Now $Is(V)$ acts on $\cc^V_{a_*,b_*}$ by $\g:g\m(\g g\g\i)$ and on $\tcc^V_{a_*,b_*}$ by
$$\g:(g,L^1,L^2,\do,L^{\s+\k})\m(\g g\g\i,\g(L^1),\g(L^2),\do,\g(L^{\s+\k})).\tag c$$
Let $\ci'=\prod_{t\in[1,\s+\k]}\{1,-1\}$. If $\e=-1$ let $\ci=\ci'$. If $\e=1$ let $\ci$ be the subgroup of 
$\ci'$ consisting of all $(\o_t)_{t\in[1,\s+\k]}$ such that $\o_t=\o_{t+1}$ for any $t$ such that 
$\{t,t+1\}\sub[k+1,\s+\k],t=k+1\mod2$. Thus $\ci$ is a finite elementary abelian $2$-group. The folowing is 
the main result of this section.

\proclaim{Theorem 1.1} (a) $\tcc^V_{a_*,b_*}$ is nonempty;

(b) the action 1.0(c) of $Is(V)$ on $\tcc^V_{a_*,b_*}$ is transitive;

(c) the isotropy group in $Is(V)$ at any point of $\tcc^V_{a_*,b_*}$ is canonically isomorphic to $\ci$.
\endproclaim

\subhead 1.2\endsubhead
Let $a\in\NN,b\in\NN,p\in\NN_{>0}$ be such that $a+b=2p$. We set $-\e=(-1)^a=(-1)^b$. For $e\in\NN$ we 
define $n_e\in\ZZ$ by $(1-T)^a(1+T)^b=\sum_{e\in\NN}n_eT^e$. We have $n_0=0$, $n_{2p-i}=-\e n_i$ for 
$i\in[0,2p]$, $n_e=0$ if $e>2p$. We define $x_e\in\ZZ$ for $e\in\NN$ by $x_0=1$ and 
$n_0x_e+n_1x_{e-1}+\do+n_ex_0=0$ for $e\ge1$. 

\subhead 1.3\endsubhead
In the setup of 1.2, let $V$ be a $\kk$-vector space with basis $\{w_i;i\in[0,2p-1]\}$. Define $g\in GL(V)$ 
by 
$$gw_i=w_{i+1}\text{ if }i\in[0,2p-2],\qua gw_{2p-1}=\e\sum_{i\in[0,2p-1]}n_iw_i.$$
We have the identity $(1-g)^a(1+g)^b=0:V@>>>V$ that is (setting $\t=\sum_{i\in[0,2p]}n_ig^i:V@>>>V$) we have 
$\t=0$. Indeed, $\t w_0=\sum_{i\in[0,2p-1]}n_iw_i+n_{2p}gw_{2p-1}=0$. Then for $i\in[0,2p-1]$ we have 
$\t w_i=\t g^iw_0=g^i\t w_0=0$. Thus $\t$ maps each element of a basis of $V$ to $0$; hence $\t=0$. Define 
a bilinear form $(,)$ on $V$ by 

$(w_i,w_j)=0$ if $i,j\in[0,2p-1]$, $|i-j|<p$,

$(w_i,w_j)=x_s$ if $i,j\in[0,2p-1]$, $j-i=p+s$, $s\ge0$,

$(w_i,w_j)=\e x_s$ if $i,j\in[0,2p-1]$, $i-j=p+s$, $s\ge0$.
\nl
Clearly $(x,y)=\e(y,x)$ for all $x,y$ and $(,)$ is nondegenerate: the determinant of the matrix 
$((w_i,w_j))$ is $\pm1$. We show that $g$ is an isometry of $(,)$. It is enough to show that

$(gw_i,gw_j)=0$ if $|i-j|<p$,

$(gw_i,gw_j)=x_s$ if $j-i=p+s$, $s\ge0$,

$(gw_i,gw_j)=\e x_s$ if $i-j=p+s$, $s\ge0$.
\nl
This is obvious except if one or both $i,j$ are $2p-1$. If $i=2p-1,p-1<j<2p-1$, we must check that
$$(\e\sum_{i'\in[0,2p-1]}n_{i'}w_{i'},w_{j+1})=0,$$
that is 
$$\sum_{i'\in[0,j+1-p]}n_{i'}x_{j+1-i'-p}=0$$
which is true since $j+1-p>0$. If $i=2p-1,0\le j<p-1$, we must check that
$$(\e\sum_{i'\in[0,2p-1]}n_{i'}w_{i'},w_{j+1})=\e x_{2p-1-j-p}$$
that is, 
$$\sum_{i'\in[j+1+p,2p-1]}n_{i'}x_{i'-j-1-p}=\e x_{p-1-j}$$
that is, 
$$-\e\sum_{i'\in[j+1+p,2p-1]}n_{2p-i'}x_{i'-j-1-p}=\e x_{p-1-j}$$
that is, 
$$\sum_{i'\in[j+1+p,2p]}n_{2p-i'}x_{i'-j-1-p}=0$$
which is true since $p-j-1>0$.

If $i=2p-1,j=p-1$, we must check that
$$(\e\sum_{i'\in[0,2p-1]}n_{i'}w_{i'},w_p)=\e x_0$$
that is, $n_0x_0=x_0$ which is obvious. The case where $j=2p-1,i<2p-1$ is entirely similar. 
It remains to show (in the case where $i=j=2p-1$) that 
$$(\e\sum_{i'\in[0,2p-1]}n_{i'}w_{i'},\e\sum_{i'\in[0,2p-1]}n_{i'}w_{i'})=0.$$
If $\e=-1$ this is obvious since $(x,x)=0$ for any $x$. Now assume that $\e=1$. We must show:
$$2\sum_{i'\in[0,p-1]}\sum_{u\in[0,p-1-i']}n_un_{u+p+i'}x_{i'}=0$$
that is,
$$\sum_{u\in[0,p-1]}n_u\sum_{i'\in[0,p-1-u]}n_{p-u-i'}x_{i'}=0.$$
We have $\sum_{i'\in[0,p-u]}n_{p-u-i'}x_{i'}=0$ if $p>u$ hence it is enough to show that
$$\sum_{u\in[0,p-1]}n_un_0x_{p-u}=0$$ 
that is,
$$\sum_{u\in[0,p-1]}n_ux_{p-u}=0.$$ 
We have 
$$\sum_{u\in[0,p]}n_ux_{p-u}=0$$
since $p>0$. Hence it is enough to show that $n_p=0$. This follows from $n_p=-\e n_p$. (We use that $\e=1$.) 

Now $g\in GL(V)$ is regular in the sense of Steinberg and satisfies $(g-1)^a(g+1)^b=0$ on V. Hence 
$V=V^+\op V^-$ where $g$ acts on $V^+$ as a single unipotent Jordan block of size $a$ and $-g$ acts on $V^-$
as a single unipotent Jordan block of size $b$. Note that if $\e=1$ we have $\det(g)=(-1)^a=-1$. It follows 
that, if $L$ is the line spanned by $w_0$ and $a_*=(a,0,0,\do),b_*=(b,0,0,\do)$, then 
$(g,L)\in\tcc^V_{a_*,b_*}$. In particular, $\tcc^V_{a_*,b_*}\ne\emp$.

\subhead 1.4\endsubhead
We now consider a variant of the situation in 1.3. In the setup of 1.2, let $V,(,)$ be as in 1.0. (Recall 
that $-\e=(-1)^a=(-1)^b$.) Let $g\in Is(V)$. We assume that $\dim V=2p$ and that on the generalized 
$1$-eigenspace of $g$, $g$ is a 
single unipotent Jordan block of size $a$ or is $1$ (if $a=0$) and on the generalized $(-1)$-eigenspace 
of $g$, $-g$ is a single unipotent Jordan block of size $b$ or is $1$ (if $b=0$). Moreover we assume 
that we are given $w\in V$ such that (with notation of 1.0) we have for $i,j\in\ZZ$:

$(w_i,w_j)=0$ if $|i-j|<p$; $(w_i,w_j)=1$ if $j-i=p$. 
\nl
We show:

(a) {\it The following equalities hold for any $i,j$ in $\ZZ$:

(a1) $(w_i,w_j)=0$ if $|i-j|<p$,

(a2) $(w_i,w_j)=x_s$ if $j-i=p+s$, $s\ge0$,

(a3) $(w_i,w_j)=\e x_s$ if $i-j=p+s$, $s\ge0$.}
\nl
Note that (a3) follows from (a2). In (a1),(a2) we can assume that $i=0$. (We use 1.0(a).)
Since $(w_0,w_j)=\e(w_0,w_{-j})$ for any $j$ we can also assume in (a1) that 
$j\ge0$ so that $j\in[0,p-1]$ and (a1) holds. We prove (a2) with $i=0,j=p+s$ by induction on $s\ge0$. If 
$s=0$ the result is already known. Assume now that $s\ge1$.
Applying $(1-g)^a(g+1)^b=0$ to $w_{s-p}$ we obtain $\sum_{e\in[0,2p]}n_ew_{s-p+e}=0$. Taking $(w_0,)$ we
obtain $\sum_{e\in[0,2p]}n_e(w_0,w_{s-p+e})=0$. For $e$ in the sum we have $s-p+e\ge-p+1$; hence by (a1) we 
can assume that we have $s-p+e\ge p$. Thus $\sum_{e\in[0,2p];s-p+e\ge p}n_e(w_0,w_{s-p+e})=0$. By the 
induction hypothesis this implies 
$$\sum_{e\in[0,2p-1];s-p+e\ge p}n_ex_{s-2p+e}+(w_0,w_{s+p})=0.$$
It is then enough to show that 
$$\sum_{e\in[0,2p-1];s-p+e\ge p}n_ex_{s-2p+e}+x_s=0$$
or that
$$\sum_{e\in[0,2p];s-p+e\ge p}n_{2p-e}x_{s-2p+e}=0$$
or that
$$\sum_{h\ge0,h'\ge0;h+h'=s}n_hx_{h'}=0.$$ 
But this holds by the definition of $x_e$ since $s\ge1$.

\subhead 1.5\endsubhead
Let $p\ge0$. For $e\ge0$ we set 
$$n_e=(-1)^e(2p+1)(2p)(2p-1)\do(2p+1-e+1)e!\i.$$
For $e\ge1$ we set $x_e=2(p+e)(2p+1)(2p+2)\do(2p+e-1)e!\i$. (Note that $x_1=2p+2$). We set $x_0=1$ if $p>0$
and $x_0=2$ if $p=0$. If $p>0$ then for any $u\ge2$ we  have 
$$\sum_{j\in[0,u]}n_jx_{u-j}=0.\tag a$$
(See \cite{\WEH, line 4 of p.134}.) This shows by induction on $e$ that $x_e\in\NN$ for any $e\ge0$.

For $u\in\ZZ$ we set $f_p(u)=0$ if $|u|<p$ and $f_p(u)=x_e$ if $|u|=p+e$ with $e\ge0$. For $u\in\ZZ$ we have 
$$f_p(u)=2(2p)!\i\prod_{k\in[0,p-1]}(u^2-k^2).\tag b$$
For example, $f_0(u)=2$. Also $f_p(p)=1$ if $p\ge1$.

Setting $A_p=\sum_{e\ge0}f_p(p+e)T^e=\sum_{e\ge0}x_eT^e$ (where $T$ is an indeterminate) we have 
$(1-T)^{2p+1}A_p=1+T$ hence
$$A_p=(1-T)^{-2p-1}(1+T).\tag c$$

\subhead 1.6\endsubhead
In the setup of 1.5 let $E$ be a $\kk$-vector with basis $w_0,w_1,\do,w_{2p}$. We define a symmetric 
bilinear form $(,):E\T E@>>>\kk$ by $(w_i,w_j)=(-1)^pf_p(i-j)$ for $i,j\in[0,2p]$. We define $g\in GL(E)$ 
by

$gw_i=w_{i+1}$ if $i\in[0,2p-1]$,

$gw_{2p}=\sum_{j\in[0,2p]}n_jw_j$.
\nl
We have $(g-1)^{2p+1}=0$ hence $g:E@>>>E$ is unipotent (with a single Jordan block). We show that $g$ is an 
isometry of $(,)$. We can assume that $p>0$. It is enough to show that $(w_{i+1},gw_{2p})=(w_i,w_{2p})$ for
$i\in[0,2p-1]$ and $(gw_{2p},gw_{2p})=0$. Thus we must show that
$$\sum_{j\in[0,2p+1],e\ge0,|i+1-j|=e+p}n_jx_e=0\text{ if }i\in[0,2p-1],\tag a$$
$$\sum_{j,j'\in[0,2p],e\ge0,|j-j'|=e+p}n_jn_{j'}x_e=0.\tag b$$
Now (a) for $i$ is equivalent to (a) for $2p-1-i$ (we use the substitution $j\m 2p+1-j$); hence it is enough
to prove (a) for $i\in[p,2p-1]$. Now (a) for $i=p$ reads $x_1-(2p+1)x_0-x_0=0$ that is $x_1=2p+2$, which is 
true. For $i\in[p+1,2p-1]$, (a) reads $\sum_{j\in[0,2p+1],i+1-j=\ge p}n_jx_{i+1-j-p}=0$ that is 
(setting $u=i+1-p$), $\sum_{j\in[0,u]}n_jx_{u-j}=0$. This follows from 1.5(a) since $u\ge2$. This proves (a).

We prove (b). The left hand side of (b) equals
$$\align&\sum_{j'\in[0,2p]}n_{j'}\sum_{j\in[0,2p],e\ge0,|j-j'|=e+p}n_jx_e\\&
=\sum_{j\in[0,2p],e\ge0,|j|=e+p}n_jx_e+\sum_{j'\in[1,2p]}n_{j'}\sum_{j\in[0,2p],e\ge0,|j-j'|=e+p}n_jx_e\\&
=\sum_{j\in[0,2p],e\ge0,|j|=e+p}n_jx_e+\sum_{j'\in[1,2p]}n_{j'}\sum_{j\in[0,2p+1],e\ge0,|j-j'|=e+p}n_jx_e\\&
-\sum_{j'\in[1,2p]}n_{j'}\sum_{e\ge0,|2p+1-j'|=e+p}n_{2p+1}x_e.\endalign$$
In the last expression the second sum over $j$ is zero by (a) and the second sum over $j'$ becomes
(setting $j=2p+1-j'$) 
$$\sum_{j\in[1,2p]}n_j\sum_{e\ge0,|j|=e+p}x_e.$$
Hence the left hand side of (b) equals
$$\sum_{j\in[0,2p],e\ge0,|j|=e+p}n_jx_e-\sum_{j\in[1,2p]}n_j\sum_{e\ge0,|j|=e+p}x_e=\sum_{e\ge0,|0|=e+p}x_e$$
and this is zero since $e+p>0$. Thus (b) holds.

For any $i\in\ZZ$ we set $w_i=g^iw_0$. This agrees with the earlier notation when $i\in[0,2p]$. We show:
$$(w_i,w_j)=(-1)^pf_p(i-j)\text{ if }i,j\in\ZZ.\tag c$$
If $p=0$ there is nothing to prove since $g=1$; thus we can assume that $p\ge1$. 
We will prove (c) assuming only the identities

(d1) $(w_{p-1},w_j)=0$ if $j\in[0,2p-2]$

(d2) $(w_{p-1},w_{2p-1})=(-1)^p$.
\nl
If $|i-j|<p$ then (c) follows from (d1); if $|i-j|=p$ then (c) follows from (d2). Thus we can assume that 
$|i-j|\ge p+1$. We can also assume that $i=0$ and $j\ge0$ (hence $j\ge p+1$). We must only prove that
$$(w_0,w_j)=(-1)^px_{j-p}\text{ if }j\ge p.$$
We argue by induction on $j$. For $j=p$ the result is known. Assume that $j\ge p+1$. From 
$(g-1)^{2p+1}w_{j-2p-1}=0$ we deduce $\sum_{h\in[0,2p+1]}n_hw_{j-2p-1+h}=0$. Hence 
$\sum_{h'\in[0,2p+1]}n_{h'}w_{j-h'}=0$ and $\sum_{h\in[0,2p+1]}n_h(w_0,w_{j-h})=0$.
If $j=p+1$ we can assume that $h=0,h=1$ or $h=2p+1$ (the other terms are zero); thus,
$$n_0(w_0,w_{p+1})+n_1(w_0,w_p)+n_{2p+1}(w_0,w_{-p})=0.$$
We see that $(w_0,w_{p+1})-(2p+1)(-1)^p-(-1)^p=0$ so that $(w_0,w_{p+1})=(-1)^p(2p+2)$ as required.
Now assume that $j\ge p+2$. We have $$\sum_{h\in[0,2p+1];j-h\ge p}n_h(w_0,w_{j-h})=0.$$ Using the 
induction hypothesis this implies
$$\sum_{h\in[1,2p+1];j-h\ge p}n_h(-1)^px_{j-h-p}+(w_0,w_j)=0$$
hence it is enough to show that $$\sum_{h\in[0,2p+1];j-h\ge p}n_hx_{j-h-p}=0$$ that is
$$\sum_{h\in[0,j-p]}n_hx_{j-h-p}=0.$$ This follows from 1.5(a) with $u=j-p$ since $j-p\ge2$.

\subhead 1.7\endsubhead
We preserve the setup of 1.6. The subspace 
$E'$ of $E$ spanned by $\{w_i;i\in[0,2p-1]\}$ is clearly nondegenerate for $(,)$ hence there 
exists $\tw\in E$ such that $(w_i,\tw)=0$ for $i\in[0,2p-1]$ and $(\tw,\tw)=2$. Moreover, $\tw$ is
unique up to multiplication by $\pm1$. We have $\tw\n E'$. We can write $\tw=\sum_{i\in[0,2p]}c_iw_i$ 
where $c_i\in\kk$ are uniquely defined and $c_*:=c_{2p}\ne0$. Taking $(w_h,)$ and setting 
$\bc_i=c_i/c_*$ we obtain
$$\sum_{i\in[0,2p]}\bc_if_p(i-h)=0\text{ for }h\in[0,2p-1].\tag $*$ $$
We show (setting $l_j=\bin{2p+1}{j}$):
$$\bc_i=(-1)^{i-1}(l_0+l_1+\do+l_i)\text{ if }i\in[0,p-1],$$ 
$$\bc_i=(-1)^i(l_0+l_1+\do+l_{2p-i})\text{ if }i\in[p,2p].$$
We can assume that $p\ge1$. 
Clearly ($*$) has a unique solution $\bc_i (i\in[0,2p-1])$. Note that $\bc_{2p}=1$.
If $h=p$ then ($*$) is $\bc_0+1=0$. If $h\in[p+1,2p-1]$ then ($*$) is $\sum_{i\in[0,h-p]}\bc_if_p(i-h)=0$. If 
$h\in[0,p-1]$ then ($*$) is $\sum_{i\in[h+p,2p]}\bc_if_p(i-h)=0$. It is enough to show:
$$\sum_{i\in[0,h-p]}(-1)^{i-1}(l_0+\do+l_i)x(h-i-p)=0\text{ if }h\in[p+1,2p-1],\tag a$$
$$\sum_{i\in[h+p,2p]}(-1)^i(l_0+\do+l_{2p-i})x(i-h-p)=0\text{ if }h\in[0,p-1].\tag b$$
We rewrite equation (b) (using $i\m 2p-i$ and $h\m 2p-h$) as
$$\sum_{i\in[0,h-p]}(-1)^i(l_0+\do+l_i)x(h-i-p)=0.\tag c$$
Here $h\in[p+1,2p]$. Note that (c) contains (a) as a special case. Thus it is enough to prove (c).
We prove (c) by induction on $h$. If $h=p+1$ then equation (c) is $l_0x_1-(l_0+l_1)x_0=0$ that is 
$2p+2-(2p+2)=0$, which is correct. If $h\ge p+2$ we have $\sum_{i\in[0,h-p]}(-1)^il_ix(h-i-p)=0$. Hence 
in this case (c) is equivalent to $\sum_{i\in[1,h-p]}(-1)^i(l_0+\do+l_{i-1})x(h-i-p)=0$ which is the same 
as equation (c) with $h$ replaced by $h-1$ (this holds by the induction hypothesis). This proves 
(c) hence (a),(b).

We show:
$$(w_{2p},\tw)c_*=2.\tag d$$
Indeed, we have 
$$2=(\tw,\tw)=(\sum_{i\in[0,2p]}c_iw_i,\tw)=c_{2p}(w_{2p},\tw),$$
as desired. We show:
$$c_*^2=2^{-2p}.\tag e$$
We have 
$$2=(w_{2p},\tw)c_*=(w_{2p},\sum_{i\in[0,2p]}c_iw_i)c_*=\sum_{i\in[0,2p]}c_i(-1)^pf_p(2p-i)c_*.$$
Thus
$$2c_*^{-2}=\sum_{i\in[0,p]}\bc_i(-1)^pf_p(2p-i).$$
If $p=0$, this reads $2c_*^{-2}=\bc_0f_0(0)=2$ hence (e) follows. 
If $p\ge1$, we have $(w_0,\tw)=0$ hence $0=\sum_{i\in[0,2p]}\bc_i(-1)^pf_p(i)$ hence 
$0=\sum_{i\in[p,2p]}\bc_i(-1)^pf_p(i)$ that is, 
$$0=\sum_{i\in[0,p]}\bc_{2p-i}(-1)^pf_p(2p-i).$$
Adding to 
$$2c_*^{-2}=\sum_{i\in[0,p]}\bc_i(-1)^pf_p(2p-i)$$
we get
$$2c_*^{-2}=\sum_{i\in[0,p]}(\bc_i+\bc_{2p-i})(-1)^pf_p(2p-i).$$ 
Now $\bc_i+\bc_{2p-i}=0$ if $i\in[0,p-1]$ hence 
$$2c_*^{-2}=2(-1)^p\bc_p=2(l_0+l_1+\do+l_p)=2^{2p+1}$$
and (e) follows.

From (e) we see that, by replacing if necessary $\tw$ by $-\tw$ we can assume that
$$c_*=2^{-p}.\tag f$$
This condition determines $\tw$ uniquely.

We show that for $h\in\ZZ$:
$$(w_h,\tw)=2^{p+1}h(h-1)(h-2)\do(h-2p+1)(2p)!\i.\tag g$$
We must show that for $h\in\ZZ$:
$$\sum_{i\in[0,2p]}c_i(-1)^pf_p(i-h)=2^{p+1}h(h-1)(h-2)\do(h-2p+1)(2p)!\i$$
or that
$$\sum_{i\in[0,2p]}\bc_i(-1)^pf_p(i-h)=2^{2p+1}h(h-1)(h-2)\do(h-2p+1)(2p)!\i.$$
It is enough to prove this equality in $\ZZ$. The left hand side is a polynomial in $h$ with rational
coefficients of degree $\le 2p$ which vanishes for $h\in[0,2p-1]$ in which the coefficient of $h^{2p}$ is
$$\align&\sum_{i\in[0,2p]}\bc_i(-1)^p2(2p)!\i=(-1)^p\bc_p2(2p)!\\&
=(l_0+l_1+\do+l_p)2(2p)!\i=(-1)^p2^{2p}2(2p)!\i.\endalign$$
Hence it is equal to the right hand side.

For any $h\in\ZZ$, $\tw_h$ is defined as in 1.0. We show:
$$(\tw_0,\tw_h)=2(-1)^h\text{ if }h\in[0,p];\qua(\tw_0,\tw_{p+1})=2(-1)^{p+1}+(-1)^p2^{2p+2}.\tag h$$
We can assume that $h\ge1$. We have 
$$\align&(\tw_0,\tw_h)=(\sum_{i\in[0,2p]}c_iw_i,\tw_h)\\&=
\sum_{i\in[0,2p]}c_i(w_{i-h},\tw_0)=\sum_{i\in[0,2p]}2\bc_i(i-h)(i-h-1)\do(i-h-2p+1)(2p)!\i\\&
=\sum_{i\in[0,h-1];i\ne p}2(-1)^{i-1}(l_0+\do+l_i)(i-h)(i-h-1)\do(i-h-2p+1)(2p)!\i\\&
+\d_{h,p+1}2(-1)^p(l_0+\do+l_p)(p-h)(p-h-1)\do(p-h-2p+1)(2p)!\i\\&
=\sum_{i\in[0,h-1]}2(-1)^{i-1}(l_0+\do+l_i)(i-h)(i-h-1)\do(i-h-2p+1)(2p)!\i\\&
+2\d_{h,p+1}2(-1)^p(l_0+\do+l_p).\endalign$$
Now $4(-1)^p(l_0+\do+l_p)=(-1)^p2^{2p+2}$. It remains to show that 
$$\sum_{i\in[0,h-1]}(-1)^{i-1}(l_0+\do+l_i)(h-i)(h-i+1)\do(h-i+2p-1)(2p)!\i=(-1)^h$$
for $h\in[1,p+1]$, or setting $h'=h-1,u=h'-i$:
$$\sum_{i\ge0,u\ge0,i+u=h'}(-1)^i(l_0+\do+l_i)(u+1)(u+2)\do(u+2p)(2p)!\i=(-1)^{h'}$$
for $h'\in[0,p]$. We shall actually show that this holds for any $h'\ge0$. It is enough to show that
for an indeterminate $T$ we have
$$\sum_{i\ge0,u\ge0}(-1)^i(l_0+\do+l_i)T^i(u+1)(u+2)\do(u+2p)(2p)!\i T^u=\sum_{h'\ge0}(-1)^{h'}T^{h'}$$
or that
$$\sum_{i\ge0}(-1)^i(l_0+\do+l_i)T^i(1-T)^{-2p-1}=(1+T)\i$$
or that 
$$l_0(1-T+T^2-\do)+l_1(-T+T^2-T^3)+\do)(1-T)^{-2p-1}=(1+T)\i$$
or that
$$(1+T)\i(l_0-l_1T+l_2T^2-\do)(1-T)^{-2p-1}=(1+T)\i.$$
This is obvious.

\subhead 1.8\endsubhead
We preserve the setup of 1.7. For $h\in\ZZ$ we show 

(a) {\it $(\tw_0,\tw_h)=\sum_{r\in[0,p]}(-1)^r2^{2r}f_r(h)$. In particular, $(\tw_0,\tw_h)\in2\ZZ$.}
\nl
We must prove the equality
$$\sum_{i\in[0,2p]}2\bc_i(i-h)(i-h-1)\do(i-h-2p+1)(2p)!\i=\sum_{r\in[0,p]}(-1)^r2^{2r}f_r(h)\tag a{}'$$
in $\kk$. It is enough to prove that (a${}'$) holds in $\ZZ$. Let $F_p(h)$ be the left hand side of 
(a${}'$). It can be viewed as a polynomial with rational coefficients in $h$ of degree $\le2p$ in which the
coefficient of $h^{2p}$ is
$$\sum_{i\in[0,2p]}2\bc_i(2p)!\i=2\bc_p(2p)!\i=2(-1)^p(l_0+\do+l_p)(2p)!\i=2(-1)^p2^{2p}(2p)!\i.$$
(We have used that $\bc_i+\bc_{2p-i}=0$ if $i\ne p$.) Thus
$$F_p(h)=(-1)^p2^{2p+1}(2p)!\i h^{2p}+\text{ lower powers of }h.$$
In the case where $p=0$ this implies that $F_p(h)=2$ so that (a') holds. We now assume that $p\ge1$. Note 
that $F_p(-h)=F_p(h)$ for $h\in\ZZ'$; an equivalent statement is that $(\tw_0,\tw_h)=(\tw_0,\tw_{-h})$ which
follows from the definitions. We see that $F_p(-h)=F_p(h)$ as polynomials in $h$. 
Now $F_p-F_{p-1}$ is a polynomial of degree $2p$ in $h$ whose value at
$h\in[0,p-1]$ is $2(-1)^h-2(-1)^h=0$. Using this and $F_p(-h)=F_p(h)$ we see that
$$F_p(h)-F_{p-1}(h)=(-1)^p2^{2p+1}(2p)!\i h^2(h^2-1)\do(h^2-(p-1)^2).$$
From this we see by induction on $p$ that (a${}'$) holds.

It follows that, if $L$ is the line spanned by $w_0$, $L'$ is the line spanned by $\tw_0$ and 
$a_*=(2p+1,0,0,\do),b_*=(0,0,0,\do)$ then $(g,L,L')\in\tcc^E_{a_*,b_*}$. In particular, 
$\tcc^E_{a_*,b_*}\ne\emp$.

$N^{2p_t}L^t_0\sub V'$ for $t=1,2$. But for $t=1,2$ we have 
$N^{2p_t}L^t_0\sub N^{2p_2-1}V\sub V'$ since $2p_t-2p_2+1\ge0$. This proves (d).

 \mpb

\subhead 1.9\endsubhead
We now consider a variant of the situation in 1.6.
In the setup of 1.5, we consider a $\kk$-vector space $E$ of dimension $2p+1$ with a given nondegenerate 
symmetric bilinear form $(,):E\T E@>>>\kk$ and a unipotent isometry $g:E@>>>E$ of $(,)$
such that $g$ is a single unipotent Jordan block (of size $2p+1$). Moreover we assume that we are given
$\tw\in E$ and (if $p\ge1$) $w\in E$ such that (with notation of 1.0) for $i,j\in\ZZ$ we have:

$(w_i,w_j)=0$ if $|i-j|<p$; $(w_i,w_j)=(-1)^p$ if $|i-j|=p$ (with $p\ge1$),

$(w_i,\tw_j)=0$ if $i-j\in[0,2p-1]$,

$(\tw_i,\tw_j)=2$ if $i=j$.
\nl
We show:

(a) {\it After possibly replacing $\tw$ by $-\tw$, the following equalities hold for any $i,h$ in $\ZZ$:

(a1) $(w_i,w_h)=(-1)^pf_p(i-h)$ if $p\ge1$,

(a2) $(w_h,\tw_0)=2^{p+1}h(h-1)(h-2)\do(h-2p+1)(2p)!\i$ if $p\ge1$,

(a3) $(\tw_0,\tw_h)=\sum_{r\in[0,p]}(-1)^r2^{2r}f_r(h)$.}
\nl
Now the proof of (a1) is exactly as in 1.6. We show:

(b) if $p\ge1$ then $\{w_i;i\in[0,2p]\}$ is linearly independent.
\nl
Assume that this is not true. Then $w_{2p}$ belongs to $E'$, the span of $\{w_i;i\in[0,2p-1]\}$; hence $E'$
is a $g$-stable hyperplane. Note that $g$ acts on $E'$ as a unipotent linear map with a single Jordan block 
(of size $2p$). By (a1), $(,)_{E'}$ is nondegenerate. Hence $g:E@>>>E$ has a Jordan block of size $2p$ and
one of size $1$; this contradicts our assumption that $g$ has a single Jordan block of size $2p+1$. This 
contradiction proves (b).

By (b) we can write uniquely (assuming $p\ge1$)
$\tw_0=\sum_{i\in[0,2p+1]}c_iw_i$ where $c_i\in\kk$. Note that $c_{2p+1}\ne0$. 
(Otherwise, $\tw_0$ would be contained in $E'$; on the other hand $\tw_0$ is perpendicular to $E'$ 
contradicting the nondegeneracy of $(,)|_{E'}$.) We set $c_*=c_{2p+1}$, $\bc_i=c_ic_*\i$ $(i\in[0,2p+1])$. 
By repeating the arguments in 1.7 we see that $c_*=\pm2^{-p}$. Replacing if necessary $\tw$ by $-\tw$ we
can assume that $c_*=2^{-p}$. Now (a2) and (a3) are proved exactly as in 1.7, 1.8.
If $p=0$ then $\tw_h=\tw_0$ for any $h\in\ZZ$ hence $(\tw_0,\tw_h)=(\tw_0,\tw_0)=f_0(0)=2$. Thus (a3) holds
again.

\subhead 1.10\endsubhead
We fix two integers $p_1,p_2$ such that $p_1\ge p_2\ge1$.Let $V',V''$ be two $\kk$-vector spaces of 
dimension $2p_1+1,2p_2-1$ respectively. Let $V=V'\op V''$. Assume that $V'$ has a given basis 
$z_0,z_1,\do,z_{2p_1}$ and that $V''$ has a given basis $v_0,v_1,\do,v_{2p_2-2}$.
We define a symmetric bilinear form $(,)$ on $V$ by
$$(z_i,z_j)=(-1)^{p_1}f_{p_1}(i-j)\text{ for }i,j\in[0,2p_1],$$
$$(v_i,v_j)=(-1)^{p_2-1}f_{p_2-1}(i-j)\text{ for }i,j\in[0,2p_2-2],$$
$$(z_i,v_j)=(v_j,z_i)=0\text{ for }i\in[0,2p_1], j\in[0,2p_2-2].$$
(Notation of 1.5.) We define $g\in GL(V)$ by
$$gz_i=z_{i+1}\text{ if }i\in[0,2p_1-1],$$
$$gz_{2p_1}=\sum_{j\in[0,2p_1]}(-1)^j\bin{2p_1+1}{j}z_j,$$
$$gv_i=v_{i+1}\text{ if }i\in[0,2p_2-3],$$
$$gv_{2p_2-2}=\sum_{j\in[0,2p_2-2]}(-1)^j\bin{2p_2-1}{j}v_j.$$
Note that $g:V@>>>V$ is unipotent and that $V',V''$ are $g$-stable ($g$ has a single Jordan block on $V'$
and a single Jordan block on $V''$). By 1.6, $g:V@>>>V$ is an isometry. For $i\in\ZZ$ we set
$z_i=g^iz_0\in V',v_i=g^zv_0\in V''$. This agrees with our earlier notation. By 1.6 we have for $i,j\in\ZZ$:
$$(z_i,z_j)=(-1)^{p_1}f_{p_1}(i-j),\qua (v_i,v_j)=(-1)^{p_2-1}f_{p_2-1}(i-j).$$
As in 1.7, 1.8, there is a unique vector $\tz_0\in V'$ and a unique vector $\tv_0\in V''$ such that for 
any $h\in\ZZ$ we have
$$(z_h,\tz_0)=2^{p_1+1}h(h-1)(h-2)\do(h-2p_1+1)(2p_1)!\i,$$
$$(\tz_0,\tz_h)=\sum_{r\in[0,p_1]}(-1)^r2^{2r}f_r(h),$$
$$(v_h,\tv_0)=2^{p_2}h(h-1)(h-2)\do(h-2p_2+3)(2p_2-2)!\i,$$
$$(\tv_0,\tv_h)=\sum_{r\in[0,p_2-1]}(-1)^r2^{2r}f_r(h).$$
For $i\in\ZZ$ we set $\tz_i=g^i\tz_0\in V',\tv_i=g^i\tv_0\in V''$. By 1.7 we have

$(\tz_0,\tz_h)=2(-1)^h$ if $h\in[0,p_1]$,

$(\tz_0,\tz_{p_1+1})=2(-1)^{p_1+1}+(-1)^{p_1}2^{2p_1+2}$,

$(\tv_0,\tv_h)=2(-1)^h$ if $h\in[0,p_2-1]$,

$(\tv_0,\tv_{p_2})=2(-1)^{p_2}+(-1)^{p_2-1}2^{2p_2}$.
\nl
We fix $\z\in\kk$ such that $\z^2=-1$. We set   
$$\x=2^{-p_2}(\tz_{-p_2}+\z\tv_0)\in V.$$
Let $h\in\ZZ$. We show:
$$(z_h,\x)=2^{p_1-p_2+1}(h+p_2)(h+p_2-1)(h+p_2-2)\do(h+p_2-2p_1+1)(2p_1)!\i.$$
In particular, $(\z_h,\x)\in2\ZZ$.

We have
$$\align&(z_h,\x)=2^{-p_2}(z_h,\tz_{-p_2})=2^{-p_2}(z_{h+p_2},\tz_0)\\&
=2^{-p_2}2^{p_1+1}(h+p_2)(h+p_2-1)(h+p_2-2)\do(h+p_2-2p_1+1)(2p_1)!\i,\endalign$$
as desired. In particular we have 
$$(\z_h,\x)=0\text{ if }h\in[-p_2,2p_1-p_2-1].$$
Let $h\in\ZZ$. We set $\x_h=g^h\x$. Using the definitions we see that
$$(\x_0,\x_h)=2^{-2p_2}((\tz_0,\tz_h)-(\tv_0,\tv_h)).$$
From this we deduce using the formulas above that
$$(\x_0,\x_h)=0\text{ if }h\in[-p_2+1,p_2-1],$$
$$(\x_0,\x_h)=(-1)^{p_2}\text{ if }h=p_2,$$
$$(\x_0,\x_h)=\sum_{r\in[p_2,p_1]}(-1)^r2^{2r-2p_2}f_r(h)\text{ for }h\in\ZZ.$$
It follows that, if $L$ is the line in $V$ spanned by $z_0$, $L'$ is the line in $V$ spanned by $\x$ and
$a_*=(2p_1+1,2p_2-1,0,0,\do),b_*=(0,0,\do)$, then $(g,L,L')\in\tcc^V_{a_*,b_*}$. In particular
$\tcc^V_{a_*,b_*}\ne\emp$. 

\subhead 1.11\endsubhead
We now consider a variant of the situation in 1.10. Let $p_1,p_2$ be as in 1.10; let $V,\e,(,)$ be as in 
1.0. Let $g\in Is(V)$. We assume that $\e=1,\dim V=2p_1+2p_2$ and that $g$ is unipotent with exactly two 
Jordan blocks: one of size $2p_1+1$ and one of size $2p_2-1$. Moreover we assume that we are given 
$z\in V,\x\in V$ such that (with notation of 1.0) we have for $i,j\in\ZZ$:

$(z_i,z_j)=0$ if $|i-j|<p_1$, $(z_i,z_j)=(-1)^{p_1}$ if $|i-j|=p_1$,

$(\x_i,\x_j)=0$ if $|i-j|<p_2$, $(\x_i,\x_j)=(-1)^{p_2}$ if $|i-j|=p_2$,

$(z_i,\x_j)=0$ if $i-j\in[-p_2,2p_1-p_2-1]$.
\nl
For $u\in\ZZ$ we set $\a_u=(z_i,z_j),\b_u=(\x_i,\x_j),\g_u=(z_i,\x_j)$ where $i,j\in\ZZ$ are such that
$i-j=u$. (These are well defined by 1.0(a).) Note that $\a_u=a_{-u},\b_u=\b_{-u}$. We show:

(a) {\it After possibly replacing $\x$ by $-\x$, the following equalities hold for any $u\in\ZZ$:

(a1) $\a_u=(-1)^{p_1}f_{p_1}(u)$,

(a2) $\g_u=2^{p_1-p_2+1}(u+p_2)(u+p_2-1)(u+p_2-2)\do(u+p_2-2p_1+1)(2p_1)!\i$,

(a3) $\b_u=\sum_{r\in[p_2,p_1]}(-1)^r2^{2r-2p_2}f_r(u)$.}
\nl
(Notation of 1.5.) When $z_i,\x_i$ are replaced by the vectors with the same name in 1.10, the quantities
$\a_u,\b_u,\g_u$ become $\a_u^0,\b_u^0,\g_u^0$ (which were computed in 1.10). Then (a1)-(a3) are equivalent 
to the equalities $\a_u=\a_u^0,\b_u=\b_u^0,\g_u=\g_u^0$.

We prove (a1). If $|u|\le p_1$ then (a1) is clear. Thus we can assume that $|u|\ge p_1+1$. We can also 
assume that $u\ge0$ (hence $u\ge p_1+1$). We must only prove that

$(z_0,z_u)=(-1)^{p_1}x_{u-p_1}$ if $u\ge p_1$
\nl
where $x_h$ is as in 1.5 (with $p=p_1$). As in the proof of 1.6(c) we argue by induction on $u$. For 
$u=p_1$ the result is known. Assume that $u\ge p_1+1$. We have $(g-1)^{2p_1+1}=0$ on $V$ hence 
$(g-1)^{2p_1+1}z_{u-2p_1-1}=0$ that is $$\sum_{h\in[0,2p_1+1]}n_hz_{u-2p_1-1+h}=0.$$ Hence 
$$\sum_{h'\in[0,2p_1+1]}n_{h'}z_{u-h'}=0$$ and $$\sum_{h\in[0,2p_1+1]}n_h(z_0,z_{u-h})=0.$$
If $u=p_1+1$ we can assume that $h=0,h=1$ or $h=2p_1+1$ (the other terms are zero); thus,
$$n_0(z_0,z_{p_1+1})+n_1(z_0,z_{p_1})+n_{2p_1+1}(z_0,z_{-p_1})=0.$$
We see that $(z_0,z_{p_1+1})-(-1)^{p_1}(2p_1+1)-(-1)^{p_1}=0$ so that $$(z_0,z_{p_1+1})=(-1)^{p_1}(2p_1+2),$$
as required. Now assume that $u\ge p_1+2$. We have $$\sum_{h\in[0,2p_1+1];j-h\ge p_1}n_h(z_0,z_{u-h})=0.$$ 
Using the induction hypothesis this implies
$$\sum_{h\in[1,2p_1+1];u-h\ge p_1}n_h(-1)^{p_1}x_{u-h-p_1}+(z_0,z_u)=0$$
hence it is enough to show that $$\sum_{h\in[0,2p_1+1];u-h\ge p_1}n_hx_{u-h-p_1}=0$$ that is
$$\sum_{h\in[0,u-p_1]}n_hx_{u-h-p_1}=0.$$ This follows from 1.5(a) with $u$ replaced by $u-p_1$ since 
$u-p_1\ge2$.

The proof of (a2),(a3) will be given in 1.12-1.16 where the setup of this subsection is preserved.

\subhead 1.12\endsubhead
We show: 

(a) the set $\{z_i;i\in[0,2p_1]\}$ is linearly independent.
\nl
Assume that this is not true. Then $z_{2p_1}\in E$, the span of $\{z_i;i\in[0,2p_1-1]\}$. Hence $E$ is 
$g$-stable and its perpendicular $E^\pe$ is $g$-stable. By assumption we have  $\x_{p_2}\in E^\pe$. Since 
$E^\pe$ is $g$-stable we see that $\x_i\in E^\pe$ for all $i\in\ZZ$. Thus $E'$, the span of 
$\{\x_i;i\in[0,2p_2-1]\}$, is contained in $E^\pe$. By assumption, $E'$ has dimension $2p_2$ which is the
same as $\dim E^\pe$. Hence $E'=E^\pe$. Since $V=E\op E'$, we see that $V=E\op E^\pe$ with both summands 
being $g$-stable. Now $g$ acts on $E$ as a single unipotent Jordan block of size $2p_1$. Thus 
$g:V@>>>V$ has a Jordan block of size $2p_1$. This contradicts the assumption that the Jordan blocks of 
$g:V@>>>V$ have sizes $2p_1+1,2p_2-1$. This proves (a).

We set $N=g-1,e=p_1-p_2$. Let $\cl$ be the span of $\{N^iz_0;i\in[2p_2,2p_1]\}$ or equivalently the span of
$\{N^{2p_2}z_i;i\in[0,2e]\}$. We show: 
$$\dim\cl=2e+1.\tag b$$
Let $\cl'$ be the span of $\{N^iz_0;i\in[2p_2,2p_1-1]\}$. We have $\dim\cl'=2e$ since 
$\{N^iz_0;i\in[0,2p_1-1]\}$ is a linearly independent set. If (b) is false we would have 
$N^{2p_2}z_0\in\cl'$. Then the span of 
$\{N^iz_0;i\in[0,2p_1-1]\}$ is $N$-stable. Hence the span of $\{g^iz_0;i\in[0,2p_1-1]\}$ is $g$-stable. 
This contradicts the proof of (a).

We show:
$$N^{2p_2}\x_0\in\cl.\tag c$$
From the structure of Jordan blocks of $N:V@>>>V$ we see that $\dim N^{2p_2}V=2e+1$. Clearly, 
$\cl\sub N^{2p_2}V$. Hence using (b) it follows that $\cl=N^{2p_2}V$ so that (c) holds. 

Using (c) we deduce 
$$N^{2p_2}\x_0=\sum_{i\in[0,2e]}c_iN^{2p_2}z_i\tag d$$
where $c_i\in\kk$ ($i\in[0,2e]$) are uniquely determined. 

\subhead 1.13\endsubhead
For $j\in\NN$ we set $m_j=(-1)^j\bin{2p_2}{j}$ so that $N^{2p_2}=\sum_{j\in[0,2p_2]}n_jg^j$. 
From 1.12(d) we deduce
$$\sum_{j\in[0,2p_2]}m_j\x_j=\sum_{i\in[0,2e],j\in[0,2p_2]}c_im_jz_{i+j}.\tag a$$
Taking $(,z_u)$ with $u\in\ZZ$, we deduce
$$\sum_{j\in[0,2p_2]}m_j\g_{u-j}=\sum_{i\in[0,2e],j\in[0,2p_2]}c_im_j\a_{u-i-j}.\tag b$$
We show:

(c1) If $u\in[p_2,2p_1-p_2-1]$ then the left hand side of (b) is $0$.

(c2) If $u=2p_1-p_2$ then the left hand side of (b) is $\g_{2p_1-p_2}$.
\nl
For (c1) it is enough to show: if $u$ is as in (c1) and $j\in[0,2p_2]$ then $u-j+p_2\in[0,2p_1-1]$. (Indeed, 
we have $u-j+p_2\le 2p_1-p_2-1+p_2=2p_1-1$ and $u-j+p_2\ge p_2-2p_2+p_2=0$.) For (c2) it is enough to show: 
if $j\in[1,2p_2]$ then $2p_1-p_2-j+p_2\in[0,2p_1-1]$. (Indeed we have $2p_1-j\le 2p_1-1$ and 
$2p_1-j\ge 2e\ge0$.)

If $u\in[p_2,p_1-1]$ then in the right hand side of (b) we have $u-i-j<p_1$; we can assume then that 
$u-i-j\le-p_1$ hence $i\ge u-j+p_1\ge p_2-2p_2+p_1=e$. Thus in this case (b) becomes (using (c1) and setting 
$u=p_1-t$):
$$\sum_{i\in[e,2e],j\in[0,2p_2]}c_im_j\a_{p_1-t-i-j}=0\text{ for }t\in[1,e].$$
\nl
Setting $c'_h=c_{2e-h}$ for $h\in[0,e]$ and with the change of variable $j\m 2p_2-j$, $i\m2e-i$ we obtain
$$\sum_{i\in[0,e],j\in[0,2p_2]}c'_im_j\a_{-p_1-t+i+j}=0\text{ for }t\in[1,e].\tag d$$
In the last sum we have $-p_1-t+i+j<p_1$. Indeed, we have 
$$-p_1-t+i+j\le-p_1-1+p_1-p_2+2p_2=p_2-1<p_1.$$ 
Hence we can restrict the sum to indices such that $-p_1-t+i+j\le-p_1$ that is $-t+i+j=-s$ where $s\ge0$. 
Thus we have 
$$\sum_{i\in[0,e],j\ge0,s\ge0;i+s+j=t}c'_im_j\a_{-p_1-s}=0\text{ for }t\in[1,e].$$
Hence 
$$(\sum_{i\in[0,e]}c'_iT^i)(\sum_{j\ge0}m_jT^j)(\sum_{s\ge0}f_{p_1}(p_1+s)T^s)
=c'_0+\text{terms of degree $>e$ in }T.$$
Thus 
$$(\sum_{i\in[0,e]}c'_iT^i)(1-T)^{2p_2}A_{p_1}=c'_0+\text{terms of degree $>e$ in }T$$
where $A_{p_1}$ is as in 1.5. Using 1.5(c) we obtain
$$(\sum_{i\in[0,e]}c'_iT^i)(1-T)^{2p_2}(1+T)(1-T)^{-2p_1-1}=c'_0+\text{terms of degree $>e$ in }T$$
hence
$$\sum_{i\in[0,e]}c'_iT^i=(1+T)\i(1-T)^{2e+1}(c'_0+\text{terms of degree $>e$ in }T).$$
We have $(1-T)^{2e+1}=\sum_{j\in[0,2e+1]}(-1)^jl_jT^j$ where $l_j=\bin{2e+1}{j}$. Hence
$$(1+T)\i(1-T)^{2e+1}=\sum_{j\in[0,e]}(-1)^j(l_0+l_1+\do+l_j)T^j+\text{terms of degree $>e$ in }T.$$
We see that 
$$c'_i=(-1)^ic'_0(l_0+l_1+\do+l_i)\text{ for }i\in[0,e].\tag e$$
In the remainder of this subsection we assume that $e>0$.
If $u=p_1$ then in the right hand side of (b) we have $u-i-j\in[-p_1,p_1]$; we can then assume that 
$u-i-j$ is $-p_1$ or $p_1$. Hence $i+j$ is $2p_1$ or $0$ and $(i,j)$ is $(2e,2p_2)$ or $(0,0)$. Thus in 
this case (b) becomes (using (c1)) $c_0+c_{2e}=0$ that is $c_0=-c'_0$ (to apply (c1) we use that $e>0$).

If $u\in[p_1+1,2p_1-p_2-1]$ then in the right hand side of (b) we have $u-i-j>-p_1$; we can assume then that 
$u-i-j\ge p_1$ hence 
$$i\le u-j-p_1\le 2p_1-p_2-1-p_1=e-1.$$
Using this and (c1) we see that (b) becomes (setting $u=p_1+t$):
$$\sum_{i\in[0,e-1],j\in[0,2p_2]}c_im_j\a_{p_1+t-i-j}=0\text{ for }t\in[1,e-1].$$
Note that in the sum we have $p_1+t-i-j>-p_1$. Indeed, we have 
$$p_1+t-i-j\ge p_1+1-p_1+p_2+1-2p_2=-p_2+2>-p_1.$$
Hence we can restrict the sum to indices such that $p_1+t-i-j\ge p_1$ that is  $p_1+t-i-j=p_1+s$ where 
$s\ge0$. Thus we have 
$$\sum_{i\in[0,e-1],j\ge0,s\ge0;i+s+j=t}c_im_j\a_{p_1+s}=0\text{ for }t\in[1,e-1].$$
For such $t$ we have also
$$\sum_{i\in[0,e-1],j\ge0,s\ge0;i+s+j=t}c'_im_j\a_{-p_1-s}=0$$
as we have seen earlier; the index $i$ cannot take the value $e$ since $i\le t$. Adding the last two 
equations and using $\a_{p_1+s}=\a_{-p_1-s}$ we obtain
$$\sum_{i\in[0,e-1],j\ge0,s\ge0;i+s+j=t}(c_i+c'_i)m_j\a_{-p_1-s}=0\text{ for }t\in[1,e-1].$$
Thus,
$$(\sum_{i\in[0,e-1]}(c_i+c'_i)T^i)(\sum_{j\ge0}m_jT^j)(\sum_{s\ge0}f_{p_1}(p_1+s)T^s)
=c+\text{terms of degree $>e$ in }T$$
where $c\in\kk$. We see that 
$$(\sum_{i\in[0,e-1]}(c_i+c'_i)T^i)(1-T)^{2p_2}A_{p_1}=c+\text{terms of degree $>e$ in }T.$$
Using again 1.5(c), we obtain
$$(\sum_{i\in[0,e-1]}(c_i+c'_i)T^i)(1-T)^{2p_2}(1+T)(1-T)^{-2p_1-1}=c+\text{terms of degree $>e$ in }T$$
hence
$$\sum_{i\in[0,e-1]}(c_i+c'_i)T^i=(1+T)\i(1-T)^{2e+1}(c+\text{terms of degree $>e$ in }T)$$
that is,
$$\sum_{i\in[0,e-1]}(c_i+c'_i)T^i=c+\text{terms of degree $>e$ in }T.$$
We see that $c_i+c'_i=0$ for $i\in[1,e-1]$. Using also (e) we see that
$$c_i=(-1)^{i+1}c'_0(l_0+l_1+\do+l_i)\text{ for }i\in[0,e-1].\tag f$$
(In the case where $i=0$ this is just $c_0=-c'_0$ which is already known.)

\subhead 1.14\endsubhead
If $u=2p_1-p_2$ then, using 1.13(b) and 1.13(c2), we have
$$\g_{2p_1-p_2}=\sum_{i\in[0,2e],j\in[0,2p_2]}c_im_j\a_{2p_1-p_2-i-j}.\tag a$$
Taking $(,\x_{p_2})$ with 1.13(a) we obtain
$$\sum_{j\in[0,2p_2]}m_j\b_{p_2-j}=\sum_{i\in[0,2e],j\in[0,2p_2]}c_im_j\g_{i+j-p_2}.$$
In the left hand side only the contribution of $j=0$ and $j=2p_2$ is $\ne0$; it is $(-1)^{p_2}$; in the 
right hand side we can assume that $i+j-p_2\ge2p_1-p_2$ (since $i+j-p_2\ge-p_2$); hence we have $i+j\ge2p_1$
and $i=2e,j=2p_2$ and the right hand side is $c_{2e}\g_{2p_1-p_2}=c'_0\g_{2p_1-p_2}$. Thus
$$2(-1)^{p_2}=c'_0\g_{2p_1-p_2}.\tag b$$
We see that $c'_0\ne0$ and using (a),(b) we have
$$2(-1)^{p_2}c'_0{}\i=\sum_{i\in[0,2e],j\in[0,2p_2]}c_im_j\a_{2p_1-p_2-i-j}.$$
In the right hand side we have $2p_1-p_2-i-j\ge-p_1$; we can assume then that either $2p_1-p_2-i-j=-p_1$
(hence $i=2e,j=2p_2$) or $2p_1-p_2-i-j\ge p_1$ (hence $i\le e$). The first case can arise only if $e=0$ 
hence it is included in the second case. Thus
$$2(-1)^{p_2}c'_0{}\i=\sum_{i\in[0,e],j\in[0,2p_2]}c_im_j\a_{2p_1-p_2-i-j}.\tag c$$
Assume now that $e>0$. From 1.13(d) with $t=e$, we have
$$0=\sum_{i\in[0,e],j\in[0,2p_2]}c'_im_j\a_{-2p_1+p_2+i+j}.\tag d$$
We now add (c) and (d) and use that $c_i+c'_i=0$ if $i\in[0,e-1]$ and $c_e=c'_e$. We get
$$2(-1)^{p_2}c'_0{}\i=2c'_e\sum_{j\in[0,2p_2]}m_j\a_{p_1-j}.$$
If $j\in[1,2p_2]$ we have $p_1-j\in[-p_1+1,p_1-1]$ hence $\a_{p_1-j}=0$. Thus
$$2(-1)^{p_2}c'_0{}\i=2c'_e\a_{p_1}=2(-1)^{p_1}c'_e.$$
By 1.13(e) we have $c'_e=(-1)^ec'_0(l_0+l_1+\do+l_e)=(-1)^ec'_02^{2e}$ hence
$$2(-1)^{p_2}c'_0{}\i=2(-1)^{p_1}(-1)^ec'_02^{2e}$$
so that $c_0'{}^2=2^{-2e}$ and $c'_0=\pm2^{-e}$. Changing if necessary $\x$ by $-\x$ we can therefore assume
that
$$c'_0=2^{-e}.\tag e$$
Assume now that $e=0$. We have $c'_0=c_0$ and (c) becomes
$$2(-1)^{p_2}c_0{}\i=\sum_{j\in[0,2p_2]}c_0m_j\a_{p_1-j}$$
that is, $2(-1)^{p_2}c_0{}\i=2c_0(-1)^{p_1}$ hence $c_0^2=1$ and $c_0=\pm1$. Changing if necessary $\x$ 
by $-\x$ we can therefore assume that $c_0=1$. Thus (e) holds without the assumption $e>0$.     

Using (e) we rewrite 1.13(e), 1.13(f) as follows:
$$c_{2e-i}=(-1)^i2^{-e}(l_0+l_1+\do+l_i)\text{ for }i\in[0,e],\tag f$$
$$c_i=(-1)^{i+1}2^{-e}(l_0+l_1+\do+l_i)\text{ for }i\in[0,e-1].\tag g$$
When $z_i,\x_i$ are replaced by the vectors with the same name in 1.10, the quantities $c_i$ become the 
quantities $c_i^0$. (Here $i\in[0,2e]$.) We show that
$$c_i=c_i^0\text{ for }i\in[0,2e].\tag h$$
By the analogue of (b) we have $2(-1)^{p_2}=c_{2e}^0\g^0_{2p_1-p_2}$. By results in 1.10 we have 
$\g^0_{2p_1-p_2}=2^{e+1}$. Hence $c_{2e}^0=(-1)^{p_2}2^{-e}$. Using this and the analogues of 1.13(e), 
1.13(f), we see that $c^0_i$ are given by the same formulas as $c_i$ in (e),(f). This proves (h).

\subhead 1.15\endsubhead
Let $C=\sum_{s\ge0}\g_{2p_1-p_2+s}T^s$, $C^0=\sum_{s\ge0}\g^0_{2p_1-p_2+s}T^s$. If $u=2p_1-p_2+t,t\ge0$,
then for any $j$ that contributes to the left hand side of 1.13(b) we have $u-j\ge-p_2$ (indeed,
$u-j\ge2p_1-p_2-2p_2\ge-p_2$) hence we can assume that in the left hand side of 1.13(b) we have 
$u-j\ge 2p_1-p_2$. Muliplying both sides of 1.13(b) by $T^t$ and summing over all $t\ge0$ we thus obtain
$$\sum_{t\ge0}\sum_{j\in[0,2p_2];t-j\ge0}m_j\g_{2p_1-p_2+t-j}T^t=
\sum_{t\ge0}\sum_{i\in[0,2e],j\in[0,2p_2]}c_im_j\a_{2p_1-p_2+t-i-j}T^t$$
The left hand side equals
$$(\sum_{j\in[0,2p_2]}m_jT^j)(\sum_{t'\ge0}\g_{2p_1-p_2+t'}T^{t'})=(1-T)^{2p_2}C.$$
Thus 
$$C=(1-T)^{-2p_2}(\sum_{t\ge0}\sum_{i\in[0,2e],j\in[0,2p_2]}c_im_j\a_{2p_1-p_2+t-i-j}T^t).$$
Similarly we have
$$C^0=(1-T)^{-2p_2}(\sum_{t\ge0}\sum_{i\in[0,2e],j\in[0,2p_2]}c_i^0m_j\a^0_{2p_1-p_2+t-i-j}T^t).$$
By 1.14(h) we have $c_i=c_i^0$. By 1.11(a1) we have
$\a_{2p_1-p_2+t-i-j}=\a^0_{2p_1-p_2+t-i-j}$ for any $i,j,t$. It follows that $C=C^0$. Hence
$$\g_{2p_1-p_2+s}=\g^0_{2p_1-p_2+s}\tag a$$ 
for any $s\ge0$. We set $C'=\sum_{t\ge0}\g_{-p_2-1-t}T^t$, $C'{}^0=\sum_{t\ge0}\g^0_{-p_2-1-t}T^t$.
If $u=p_2-1-t, t\ge0$, then for any $j$ that contributes to the left hand side of 1.13(b) we have 
$u-j\le 2p_1-p_2-1$ (indeed $u-j\le p_2-1-j\le p_2-1\le2p_1-p_2-1$) hence we can assume that in the left 
hand side of 1.13(b) we have $u-j\le-p_2-1$. With the substitution $j\m2p_2-j$ the previous inequality 
becomes $j-t\le0$ and the left hand side of 1.13(b) becomes
$$\sum_{j\in[0,2p_2]}m_j\g_{u-2p_2+j}=\sum_{j\in[0,2p_2]}m_j\g_{-p_2-1+j-t}.$$
Muliplying both sides of 1.13(b) by $T^t$ and summing over all $t\ge0$ we thus obtain
$$\sum_{t\ge0,j\ge0;t-j\ge0}m_j\g_{-p_2-1+j-t}T^t=
\sum_{t\ge0}\sum_{i\in[0,2e],j\in[0,2p_2]}c_im_j\a_{p_2-1-t-i-j}T^t.$$
The left hand side equals
$$(\sum_{j\in[0,2p_2]}m_jT^j)(\sum_{t'\ge0}\g_{-p_2-1-t'}T^{t'})=(1-T)^{2p_2}C'.$$
Thus 
$$C'=(1-T)^{-2p_2}(\sum_{t\ge0}\sum_{i\in[0,2e],j\in[0,2p_2]}c_im_j\a_{p_2-1-t-i-j}T^t).$$
Similarly we have
$$C'{}^0=(1-T)^{-2p_2}(\sum_{t\ge0}\sum_{i\in[0,2e],j\in[0,2p_2]}c_i^0m_j\a^0_{p_2-1-t-i-j}T^t).$$
By 1.14(h) we have $c_i=c_i^0$. By 1.11(a1) we have
$\a_{p_2-1-t-i-j}=\a^0_{p_2-1-t-i-j}$ for any $i,j,t$. It follows that $C'=C'{}^0$. Hence
$$\g_{-p_2-1-t}=\g^0_{-p_2-1-t}\tag b$$
for any $t\ge0$. Clearly (a),(b) imply 1.11(a2).

\subhead 1.16\endsubhead
We set $B=\sum_{s\ge0}\b_{p_2+s}T^s$, $B^0=\sum_{s\ge0}\b^0_{p_2+s}T^s$.
Let $t\ge1$. Taking $(,\x_{p_2+t})$ with 1.13(a) we obtain
$$\sum_{j\in[0,2p_2]}m_j\b_{p_2+t-j}=\sum_{i\in[0,2e],j\in[0,2p_2]}c_im_j\g_{i+j-p_2-t}.\tag a$$
For any $j$ that contributes to the left hand side of (a) we have $p_2+t-j\ge-p_2+1$ (indeed,
$p_2+t-j\ge p_2+1-2p_2=-p_2+1$) hence we can assume that in the left hand side of (a) we have 
$p_2+t-j\ge p_2$ that is $t\ge j$. Multiplying both sides of (a) by $T^t$ and summing over all $t\ge1$ we 
thus obtain
$$\sum_{t\ge1}\sum_{j\in[0,2p_2];t\ge j}m_j\b_{p_2+t-j}T^t
=\sum_{t\ge1}\sum_{i\in[0,2e],j\in[0,2p_2]}c_im_j\g_{i+j-p_2-t}T^t.$$
The left hand side equals
$$-(-1)^{p_2}+(\sum_{j\in[0,2p_2]}m_jT^j)(\sum_{t'\ge0}\b_{p_2+t'}T^{t'})=-(-1)^{p_2}+(1-T)^{2p_2}B.$$
Thus
$$B=(1-T)^{-2p_2}((-1)^{p_2}+\sum_{t\ge1}\sum_{i\in[0,2e],j\in[0,2p_2]}c_im_j\g_{i+j-p_2-t}T^t).$$
Similarly we have
$$B^0=(1-T)^{-2p_2}((-1)^{p_2}+\sum_{t\ge1}\sum_{i\in[0,2e],j\in[0,2p_2]}c_i^0m_j\g^0_{i+j-p_2-t}T^t).$$
By 1.14(h) we have $c_i=c_i^0$. By 1.11(a2) we have $\g_{i+j-p_2-t}=\g^0_{i+j-p_2-t}$ for any $i,j,t$. It 
follows that $B=B^0$. Hence $\b_{p_2+s}=\b^0_{p_2+s}$ for any $s\ge0$. This clearly implies 1.11(a3).

\subhead 1.17\endsubhead
In the setup of 1.1 we show that 1.1(a) holds by induction on $\nn$. If $\nn=0$ we have $V=0$ and 
$a_i=b_i=c_i=p_i=0$ for all $i$. We take $g=0$ and $(L^t)$ to be the empty set of lines. We obtain an 
element of $\tcc^V_{a_*,b_*}$. Now assume that $\nn>0$.

Assume first that either $a_1\ge1,b_1\ge1$ or that $\e=-1$.
We can find a direct sum decomposition $V=V'\op V''$ such that $(V',V'')=0$ and $\dim V'=a_1+b_1=2p_1$. Let 
$a'_*$ be the sequence $a_1,0,0,\do$; let $b'_*$ be the sequence $b_1,0,0,\do$; let $a''_*$ be the sequence 
$a_2,a_3,\do$; let $b''_*$ be the sequence $b_2,b_3,\do$. By the induction hypothesis we have 
$\tcc^{V''}_{a''_*,b''_*}\ne0$. By 1.3 we have $\tcc^{V'}_{a'_*,b'_*}\ne\emp$. Let 
$(g',L^1)\in\tcc^{V'}_{a'_*,b'_*}$ and let $(g'',L^2,L^3,\do)\in\tcc^{V''}_{a''_*,b''_*}$. Clearly,
$(g'\op g'',L^1,L^2,\do)\in\tcc^V_{a_*,b_*}$ hence 1.1(a) holds in this case. Thus we can assume that 
$\e=1$ and either 

(i) $a_1>0$ and $b_1=0$ or 

(ii) $a_1=0$ and $b_1>0$.
\nl
Assume that we are in case (i). We have $b_1=b_2=\do=0$ and $g$ is unipotent. If $a_2=0$ then 1.1(a) holds 
by 1.6 with $p=(a_1-1)/2$. If $a_2>0$ we can find a direct sum decomposition $V=V'\op V''$ such that 
$(V',V'')=0$ and $\dim V'=a_1+a_2$. Let $a'_*$ be the sequence $a_1,a_2,0,\do$; let $a''_*$ be the sequence 
$a_3,a_4,\do$; let $b'_*=b''_*$ be the sequence $0,0,\do$. By the induction hypothesis we have 
$\tcc^{V''}_{a''_*,b''_*}\ne\emp$. By 1.10 we have $\tcc^{V'}_{a'_*,b'_*}\ne\emp$. Let 
$(g',L^1,L^2)\in\tcc^{V'}_{a'_*,b'_*}$ and let $(g'',L^3,L^4,\do)\in\tcc^{V''}_{a''_*,b''_*}$. Clearly
$(g'\op g'',L^1,L^2,\do)\in\tcc^V_{a_*,b_*}$ hence 1.1(a) holds in this case. This completes the proof in 
case (i).

Assume now that we are in case (ii) so that $-g$ is unipotent. It is easy to check that
$\tcc^V_{g;a_*,b_*}=\tcc^V_{-g;b_*,a_*}$ and the last set is nonempty by the earlier part of the argument.
Hence $\tcc^V_{g;a_*,b_*}\ne\emp$. This completes the inductive proof of 1.1(a).

In the following result we preserve the setup of 1.1. 

\proclaim{Proposition 1.18} Let $(g,L^1,L^2,\do,L^{\s+\k})\in\tcc^V_{a_*,b_*}$. Let $f_r$ be as in 1.5. 
There exist vectors $z^t\in L^t-\{0\}$ for $t\in[1,\s+\k]$ such that (i),(ii),(iii) below hold for any
$i,j\in\ZZ$.

(i) Assume that either $t\in[1,\s],\e=-1$ or $t\in[1,k]$. Then $(z^t_i,z^t_j)=0$ if $|i-j|<p_t$,
$(z^t_i,z^t_j)=x_s$ if $j-i=p_t+s,s\ge0$ ($x_s$ as in 1.5 with $p=p_t$); $(z^t_i,z^{t'}_j)=0$ if 
$t'\in[1,\s+\k],t'\ne t$.

(ii) Assume that $\{t,t+1\}\sub[k+1,\s+\k],t=k+1\mod2$ and $\e=1$. We set $\d=1$ if $a_t>0$, $\d=-1$ if 
$b_t>0$. Then 
$$(z^t_i,z^t_j)=(-1)^{p_t}\d^{i-j}f_{p_t}(i-j),$$
$$(z^{t+1}_i,z^{t+1}_j)=\d^{i-j}\sum_{r\in[p_{t+1},p_t]}(-1)^r2^{2r-2p_{t+1}}f_r(i-j),$$
$$\align&(z^t_i,z^{t+1}_j)=\d^{i-j}2^{p_t-p_{t+1}+1}(i-j+p_{t+1})(i-j+p_{t+1}-1)\\&
\T(i-j+p_{t+1}-2)\do(i-j+p_{t+1}-2p_t+1)(2p_t)!\i,\endalign$$
$$(z^t_i,z^{t'}_j)=0\text{ if }t'\in[1,\s+\k],t'\n\{t,t+1\}.$$
(iii) Assume that $\e=1,\k=1,t=\s+1$. We set $\d=1$ if $a_t>0$, $\d=-1$ if $b_t>0$. (We have $p_t=0$.) Then 
$$(z^t_i,z^t_j)=2\d^{i-j},$$
$$(z^t_i,z^{t'}_j)=0\text{ if }t'\in[1,\s].$$
\endproclaim
We argue by induction on $\nn$. When $\nn=0$ the result is obvious. Now assume that $\nn\ge1$.

{\it Case 1.} Assume first that either $a_1\ge1,b_1\ge1$ or that $\e=-1$. We have $a_1+b_1=2p_1$. Let 
$V'=\op_{i\in[0,2p_1-1]}L^1_i\sub V$. We show that 
$$gV'=V'.\tag a$$
It is enough to show that $gL^1_{2p_1-1}\sub V'$ that is $g^{2p_1}L^1_0\sub V'$. Since $g^iL^1_0\sub V'$ for
$i\in[0,2p_1-1]$ and $a_1+b_1=2p_1$, it is enough to show that $(g-1)^{a_1}(g+1)^{b_1}L_0^1=0$. It is also 
enough to show that $(g-1)^{a_1}(g+1)^{b_1}=0$ on $V$. But this follows from the fact that 
$g\in\cc^V_{a_*,b_*}$.

Now let $V''=\op_{t\in[2,\s+\k],i\in[0,2p_t-1]}L^t_i\sub V$. We show that
$$V''=V'{}^\pe\text{ (the perpendicular to $V'$) and }V=V'\op V'{}^\pe.\tag b$$
For $t\in[2,\s]$, $i\in[0,2p_1-1]$ we have $(L^1_i,L_{p_t}^t)=0$; thus $L_{p_t}^t\in V'{}^\pe$. Since 
$V'{}^\pe$ is $g$-stable it follows that $L^t_i\sub V'{}^\pe$ for $t\in[2,\s],i\in\ZZ$. If $\k=1$ we have 
$(L_i^1,L_0^{\s+1})=0$ for $i\in[0,2p_1-1]$; thus $L_0^{\s+1}\sub V'{}^\pe$. Hence $V''\sub V'{}^\pe$. But 
these two vector spaces have the same dimension so that $V''=V'{}^\pe$. Since $V=V'\op V''$ it follows that 
$V=V'\op V'{}^\pe$. This proves (b).

Let $g'=g|_{V'}$, $g''=g_{V''}$. We show:

(c) {\it $g'$ restricted to the generalized $1$-eigenspace of $g'$ is unipotent with a single Jordan block 
of size $a_1$; $-g'$ restricted to the generalized $(-1)$-eigenspace of $g'$ is unipotent with a single 
Jordan block of size $b_1$; $g''$ restricted to the generalized $1$-eigenspace of $g''$ is unipotent with 
Jordan blocks of sizes given by the nonzero numbers in $a_2,a_3,\do$; $-g''$ restricted to the generalized 
$(-1)$-eigenspace of $g''$ is unipotent with Jordan blocks of sizes given by the nonzero numbers in 
$b_2,b_3,\do$.}
\nl
As we have seen earlier we have $(g-1)^{a_1}(g+1)^{b_1}=0$ on $V'$ (even on $V$). Also $g'\in GL(V')$ is 
regular in the sense of Steinberg and $\dim V'=a_1+b_1$. This implies (c).

Let $a'_*$ be the sequence $a_1,0,0,\do$; let $b'_*$ be the sequence $b_1,0,0,\do$; let $a''_*$ be the 
sequence $a_2,a_3,\do$; let $b''_*$ be the sequence $b_2,b_3,\do$.

Now the proposition holds when $(g,L^1,L^2,\do)$ is replaced by 
$(g'',L^2,L^3,\do)\in\tcc^{V''}_{a''_*,b''_*}$ (by the induction hypothesis) or by 
$(g',L^1)\in\tcc^{V'}_{a'_*,b'_*}$ (we choose any $z^1\in L^1-\{0\}$ such that $(z^1_i,z^1_j)=1$ for 
$|i-j|=p_1$ and we apply 1.4). Hence the proposition holds for $(g,L^1,L^2,\do)$ (since $(V',V'')=0$).

{\it Case 2.} Next we assume that $k=0$, $\e=1$, $a_1>0,a_2>0$. Then $b_1=b_2=\do=0$. We have 
$a_1=2p_1+1,a_2=2p_2-1$. Let $V'=\op_{t\in[1,2],i\in[0,2p_t-1]}L^t_i\sub V$. We show that 
$$gV'=V'.\tag d$$
Let $N=g-1$. Then $V=\op_{t\in[1,\s+\k],i\in[0,2p'_t-1]}N^iL^t_0$ is a direct sum decomposition into lines 
and $p_i=p'_i$ if $i\in[1,2]$. Now $N^{2p_2-1}(V)$ contains the lines

($*$) $N^{2p_2-1+i}L^1_0 (i=0,1,\do,2p_1-2p_2)$ and $N^{2p_2-1}L^2_0$
\nl
(whose number is $2p_1-2p_2+2$); moreover, since $N$ has Jordan blocks of sizes $a_1=2p_1+1,a_2=2p_2-1$
and others of size $<a_2$, we see that $\dim N^{2p_2-1}(V)=2p_1-2p_2+2$ so that $N^{2p_2-1}(V)$ is equal to 
the subspace spanned by ($*$) and $N^{2p_2-1}(V)\sub V'$. Now $V'$ is the subspace of $V$ spanned by the 
lines $N^iL^t_0$ with $t\in[1,2],i\in[0,2p_t-1]$. It is enough to show that $NV'\sub V'$ or that 
$N^{2p_t}L^t_0\sub V'$ for $t=1,2$. But for $t=1,2$ we have 
$N^{2p_t}L^t_0\sub N^{2p_2-1}V\sub V'$ since $2p_t-2p_2+1\ge0$. This proves (d).

 \mpb

Let $V''=\op_{t\in[3,\s+\k],i\in[0,2p'_t-1]}L^t_i\sub V$. We show that

(e) $V''=V'{}^\pe$ (the perpendicular to $V'$) and $V=V'\op V'{}^\pe$.
\nl
For $t\in[1,2],r\in[3,\s]$, $i\in[0,2p_t-1]$ we have $(L_i^t,L_{p_r}^r)=0$. Thus 
$L_{p_r}^r\sub V'{}^\pe$ for $r\in[3,\s]$. Since $V'{}^\pe$ is $g$-stable it follows that 
$L^r_i\sub V'{}^\pe$ for $r\in[3,\s],i\in\ZZ$. If $\k=1$ we have $(L_i^t,L_0^{\s+1})=0$ for 
$i\in[0,2p_t-1],t\in[1,2]$. Thus $L_0^{\s+1}\in V'{}^\pe$. Hence $V''\sub V'{}^\pe$. But these two vector 
spaces have the same dimension so that $V''=V'{}^\pe$. Since $V=V'\op V''$ it follows that 
$V=V'\op V'{}^\pe$. This proves (e).

\mpb

Let $g'=g|_{V'}$, $g''=g_{V''}$. We show:

(f) {\it $g'$ is unipotent with exactly two Jordan blocks of size $a_1,a_2$. Moreover, $g''$ is unipotent 
with Jordan blocks of sizes given by the nonzero numbers in $a_3,a_4,\do$.}
\nl
Since $V'$ is the direct sum of the lines $N^iL^t_0$, $t\in[1,2],i\in[0,2p_t-1]$ and $V'$ is $N$-stable, we
see that the kernel of $N:V'@>>>V'$ has dimension $\le2$. Hence $N:V'@>>>V'$ has either a single Jordan 
block of size $2p_1+2p_2=a_1+a_2$ or two Jordan blocks of sizes $a'_1\ge a'_2$ where $a'_1+a'_2=a_1+a_2$. 
The first alternative does not occur since the Jordan blocks of $N:V'@>>>V'$ have sizes $\le a_1$ (by (e)). 
Thus the second alternative holds. Since $a'_1,a'_2$ must form a subsequence of $a_1>a_2>a_3>\do$ and 
$a'_1+a'_2=a_1+a_2$ it follows that $a'_1=a_1$, $a'_2=a_2$. This implies (f).

Let $a'_*$ be the sequence $a_1,a_2,0,\do$; let $a''_*$ be the sequence $a_3,a_4,\do$; let $b'_*=b''_*$ be 
the sequence $b_2,b_3,\do$. Now the proposition holds when $(g,L^1,L^2,\do)$ is replaced by 
$(g'',L^3,L^4,\do)\in\tcc^{V''}_{a''_*,b''_*}$ (by the induction hypothesis) or by 
$(g',L^1,L^2)\in\tcc^{V'}_{a'_*,b'_*}$ (we choose any $z^1\in L^1-\{0\}$ such that 
$(z^1_i,z^1_j)=(-1)^{p_1}$ for $|i-j|=p_1$ and any $z^2\in L^2-\{0\}$ such that $(z^2_i,z^2_j)=(-1)^{p_2}$ 
for $|i-j|=p_2$ and we apply 1.11 by possibly changing $z^2$ to $-z^2$). Hence the proposition holds for 
$(g,L^1,L^2,\do)$ (since $(V',V'')=0$).

{\it Case 3.} Next we assume that $k=0$, $\e=1$, $a_1>0,a_2=0$. Then $b_1=b_2=\do=0$ and $\s=1,\k=1$. We 
have $a_1=2p_1+1,p_2=0,p'_2=1/2$. We choose any $z^1\in L^1-\{0\}$ such that $(z^1_i,z^1_j)=(-1)^{p_1}$ for 
$|i-j|=p_1$ and any $z^2\in L^2-\{0\}$ such that $(z^2_i,z^2_j)=2$ for $|i-j|=p_2$ and we apply 1.9 by 
possibly changing $z^2$ to $-z^2$. We see that the proposition holds for $(g,L^1,L^2,\do)$.

{\it Case 4.} Finally assume that $k=0$, $\e=1$, $b_1>0$. Then $(-g,L^1,L^2,\do)\in\tcc^V_{b_*,a_*}$ is as 
in Case 2 or 3. Let $(z^t)$ be the corresponding sequence of vectors in $V$. This sequence is the desired 
sequence for $(g,L^1,L^2,\do)$. This completes the proof.

\subhead 1.19\endsubhead
In the setup of 1.1, we show that 1.1(b) holds. We must show that 

(a) any two elements $(g,L^1,L^2,\do,L^{\s+\k})$, $(g',L'{}^1,L'{}^2,\do,L'{}^{\s+\k})$ of $\tcc^V_{a_*,b_*}$
are in the same $Is(V)$-orbit.
\nl
Since $Is(V)$ acts transitively on $\cc^V_{a_*,b_*}$ we can assume that $g=g'$. Let $z^t\in L^t$ 
$(t\in[1,\s+\k])$ be as in 1.18. Let $z'{}^t\in L'{}^t$ $(t\in[1,\s+\k])$ be the analogous vectors for 
$(g,L'{}^1,L'{}^2,\do)$  instead of $(g,L^1,L^2,\do)$. By 1.18 we have 
$$(z^t_i,z^{t'}_j)=(z'{}^t_i,z'{}^{t'}_j)\tag b$$ 
for any $i,j\in\ZZ$ and any $t,t'\in[1,\s+\k]$. Since $\{z^t_i;t\in[1,\s+\k],i\in[0,2p'_t-1]\}$ and 
$\{z'{}^t_i;t\in[1,\s+\k],i\in[0,2p'_t-1]\}$ are bases of $V$ (see 1.0(b)) we see that there is a unique 
$\g\in GL(V)$ such that $\g(z^t_i)=z'{}^t_i$ for any $t\in[1,\s+\k],i\in[0,2p'_t-1]$. From (b) we see that 
$\g\in Is(V)$. We show that 
$$\g(z^t_{i+1})=z'{}^t_{i+1}\text{ for any }t\in[1,\s+\k],i\in[0,2p'_t-1].\tag c$$
When $i+1\in[0,2p'_t-1]$ this follows from the definition of $\g$. Thus we can assume that $i=2p'_t-1$ and 
we must show that $\g(z^t_{2p'_t})=z'{}^t_{2p'_t}$ for any $t\in[1,\s+\k]$. It is enough to show that 
$(\g(z^t_{2p'_t}),z'{}^{t'}_j)=(z'{}^t_{2p'_t},z'{}^{t'}_j)$ for any $t'\in[1,\s+\k],j\in[0,2p'_t-1]$ (we 
use again that $\{z'{}^t_i;t\in[1,\s+\k],i\in[0,2p'_t-1]\}$ is a basis of $V$). We have 
$(\g(z^t_{2p'_t}),z'{}^{t'}_j)=(\g(z^t_{2p'_t}),\g(z^{t'}_j))=(z^t_{2p'_t},z^{t'}_j)$ and this is equal to 
$(z'{}^t_{2p'_t},z'{}^{t'}_j)$ by (b). Thus (c) holds. From (c) we see that $\g(g(z^t_i))=g(\g(z^t_i))$ for 
any $t\in[1,\s+\k],i\in[0,2p'_t-1].$ It follows that $\g g=g\g$. From the definition it is clear that 
$\g(L^t)=L'{}^t$ for $t\in[1,\s+\k]$. Thus (a) holds (with $g'=g$). This proves 1.1(b).

\subhead 1.20\endsubhead
In the setup of 1.1, we show that 1.1(c) holds. Let $(g,L^1,L^2,\do,L^{\s+\k})\in\tcc^V_{a_*,b_*}$ and let 
$I$ be the set of all $\g\in Is(V)$ such that $\g g\g\i=g$, $\g(L^t)=L^t$ for $t\in[1,\s+\k]$. Let 
$z^t\in L^t (t\in[1,\s+\k])$ be as in 1.18. Let $\g\in I$. If $t\in[1,\s+\k]$, we have $\g(z^t)=\o^\g_tz^t$ 
where $\o^\g_t=\pm1$. If $\{t,t+1\}\sub[k+1,\s+\k],t=k+1\mod2$ and $\e=1$, we have $\o^\g_t=\o^\g_{t+1}$.
Indeed, for some $s\in\{1,-1\}$ we have 
$$\align&s2^{p_t-p_{t+1}-1}=(z^t_{-1},z^{t+1}_{p_{t+1}})=(\g(z^t_{-1}),\g(z^{t+1}_{p_{t+1}}))\\&=
\o^\g_t\o^\g_{t+1}(z^t_{-1},z^{t+1}_{p_{t+1}})=\o^\g_t\o^\g_{t+1}s2^{p_t-p_{t+1}-1}\endalign$$
hence $\o^\g_t\o^\g_{t+1}=1$ and our claim follows. Thus, $\g\m(\o^\g_t)$ is a homomorphism $\ps:I@>>>\ci$ 
(notation of 1.0). Assume that $\g$ is in the kernel of $\ps$. Then $\g$ restricts to the identity map 
$L^t@>>>L^t$ for $t\in[1,\s+\k]$. Since $\g$ commutes with $g$ it follows that $\g$ restricts to the 
identity map on each of the lines $g^iL^t$ ($t\in[1,\s+\k]$, $i\in\ZZ$). Since these lines generate $V$ (see
1.0(b)) we see that $\g=1$. Thus $\ps$ is injective. Now let $(\o_t)\in\ci$. We define $\g\in GL(V)$ by 
$\g(z^t_i)=\o_tz^t_i$ for $t\in[1,\s+\k],i\in[0,2p'_t-1]$. From the definitions we see that 
$$(\o_tz^t_i,\o_{t'}z^{t'}_j)=(z^t_i,z^{t'}_j)\tag a$$ 
for any $i,j\in\ZZ$ and any $t,t'\in[1,\s+\k]$.

From (a) we see that $\g\in Is(V)$. We show that 
$$\g(z^t_{i+1})=\o_tz^t_{i+1}\text{ for any }t\in[1,\s+\k],i\in[0,2p'_t-1].\tag b$$
When $i+1\in[0,2p'_t-1]$ this follows from the definition of $\g$. Thus we can assume that $i=2p'_t-1$ and 
we must show that $\g(z^t_{2p'_t})=\o_tz^t_{2p'_t}$ for any $t\in[1,\s+\k]$. It is enough to show that 
$(\g(z^t_{2p'_t}),\o_{t'}z^{t'}_j)=(\o_tz^t_{2p'_t},\o_{t'}z^{t'}_j)$ for any 
$t'\in[1,\s+\k],j\in[0,2p'_t-1]$ 
(we use again that $\{z^t_i;t\in[1,\s+\k],i\in[0,2p'_t-1]\}$ is a basis of $V$). We have 
$$(\g(z^t_{2p'_t}),\o_{t'}z^{t'}_j)=(\g(z^t_{2p'_t}),\g(z^{t'}_j))=(z^t_{2p'_t},z^{t'}_j)$$ 
and this is equal to $(z^t_{2p'_t},z^{t'}_j)$ by (a). Thus (b) holds.

From (b) we see that $\g(g(z^t_i))=g(\g(z^t_i)$ for any $t\in[1,\s+\k],i\in[0,2p'_t-1]$. It follows that
$\g g=g\g$. From the definition it is clear that $\g(L^t)=L^t$ for $t\in[1,\s+\k]$. Thus $\g\in I$. We see 
that $\ps$ is surjective hence an isomorphism. This proves 1.1(c).

\subhead 1.21\endsubhead
In the setup of 1.1, assume that $\nn$ is even $\ge2$ and $\e=1$. Let $\Om$ be the set of $Is(V)^0$-orbits 
on the set of $(\nn/2)$-dimensional subspaces of $V$ which are isotropic for $(,)$; note that $|\Om|=2$. If 
$(g,L^1,L^2,\do,L^\s)\in\tcc^V_{a_*,b_*}$ then the $(\nn/2)$-dimensional subspace 
$\op_{t\in[1,\s],i\in[p_t,2p_t-1]}L^t_i$ of $V$ is isotropic for $(,)$. Hence we have a partition
$$\tcc^V_{a_*,b_*}=\sqc_{\co\in\Om}\tcc^V_{a_*,b_*;\co}$$
where for $\co\in\Om$, $\tcc^V_{a_*,b_*;\co}$ is the set of all
$(g,L^1,L^2,\do,L^\s)\in\tcc^V_{a_*,b_*}$ such that $\op_{t\in[1,\s],i\in[p_t,2p_t-1]}L^t_i\in\co$. Now

(a) {\it the action 1.0(c) of $Is(V)$ restricts for any $\co\in\Om$ to an action of $Is(V)^0$ on 
$\tcc^V_{a_*,b_*;\co}$;}

(b) {\it if $\g\in Is(V)-Is(V)^0$ then the action of $\g$ on $\tcc^V_{a_*,b_*}$ maps $\tcc^V_{a_*,b_*;\co}$ 
onto $\tcc^V_{a_*,b_*;\Om-\co}$.}
\nl
For any $\co\in\Om$ we have the following variant of Theorem 1.1:

(c) $\tcc^V_{a_*,b_*;\co}\ne\emp$;

(d) {\it the action (a) of $Is(V)^0$ on $\tcc^V_{a_*,b_*;\co}$ is transitive;}

(e) {\it the isotropy group in $Is(V)^0$ at any point of $\tcc^V_{a_*,b_*;\co}$ is canonically isomorphic 
to $\ci$.}
\nl
Now (c) follows immediately from (b) and 1.1(a). We prove (d). Let 
$$(g,L^1,L^2,\do)\in\tcc^V_{a_*,b_*;\co},\qua (g',L'{}^1,L'{}^2,\do)\in\tcc^V_{a_*,b_*;\co}.$$
By 1.1(b) we can find $\g\in Is(V)$ which carries $(g,L^1,L^2,\do)$ to $(g',L'{}^1,L'{}^2,\do)$. By (b) we 
have automatically $\g\in Is(V)^0$. Hence (d) holds.

To prove (e) it is enough to show that if $\g$ is in the isotropy group in $Is(V)$ at $(g,L^1,L^2,\do)$, then
$\det(\g)=1$. Let $(\o_t)=\ps(\g)$ be as in 1.20. From the proof in 1.20 we see that
$\det(\g)=\prod_{t\in[1,\s]}\o_t^{2p_t}$. Since $\o_t=\pm1$ we see that $\det(\g)=1$, as required.

We now show:

(f) {\it If $a_1>0$, $b_1>0$ and $(g,L^1,L^2,\do)\in\tcc^V_{a_*,b_*;\co}$, then there exists 
$\g\in I'$ (the isotropy group in $Is(V)^0$ at $(g,L^1,L^2,\do)$) such that for $\d\in\{1,-1\}$, the
restriction of $\g$ to the generalized $\d$-eigenspace of $g$ has determinant $-1$.}
\nl
Define $(\o_t)$ by $\o_1=-1,\o_t=1$ for $t\in[2,\s]$. In our case we have $k\ge1$ hence $(\o_t)\in\ci$. Let 
$V'=\sum_{i\in\ZZ}L^1_i$, $V''=\sum_{t\in[2,\s],i\in\ZZ}L^t_i$. By 1.18, $V=V'\op V''$ (orthogonal direct 
sum). Define $\g\in I'$ by $\ps(\g)=(\o_t)$ (notation of 1.20). Then $\g$ acts as identity on $V''$ and as 
$-1$ times the identity on $V'$. It is enough to prove that the restriction of $\g$ to the generalized 
$\d$-eigenspace of $g_{V'}$ has determinant $-1$ or that this generalized $\d$-eigenspace has odd dimension. 
But this dimension is $a_1$ (if $\d=1$) and $b_1$ (if $\d=-1$) and $a_1,b_1$ are odd.

\subhead 1.22\endsubhead
In the setup of 1.1, assume that $\nn$ is odd (hence $\e=1$) and that $\cc^V_{a_*,b_*}\sub Is(V)^0$.
We have the following variant of Theorem 1.1:

(a) {\it the restriction of the action 1.0(c) to $Is(V)^0$ is transitive on $\tcc^V_{a_*,b_*}$;}

(b) {it the isotropy group in $Is(V)^0$ at any point of $\tcc^V_{a_*,b_*}$ is canonically isomorphic to a 
subgroup of $\ci$ of index $2$.}
\nl
Note that if $\g\in Is(V)-Is(V)^0$ then $-\g\in Is(V)^0$. Moreover $-1\in Is(V)$ acts trivially on
$\tcc^V_{a_*,b_*}$; hence (a) follows from 1.1(b). Now let $\g$ be in the isotropy group in $Is(V)$ at 
$(g,L^1,L^2,\do)$ and let $(\o_t)=\ps(\g)$ be as in 1.20. We have 
$$\det(\g)=\o_{\s+1}\prod_{t\in[1,\s]}\o_t^{2p_t}=\o_{\s+1}.$$
Thus the condition that $\g\in Is(V)^0$ is equivalent to the condition that $\o_{\s+1}=1$. This proves (b).

We now show:

(c) {\it If $a_1>0$, $b_1>0$ and $(g,L^1,L^2,\do)\in\tcc^V_{a_*,b_*}$ with $g\in Is(V)^0$ then there exists 
$\g\in I'$ (the isotropy group in $Is(V)^0$ at $(g,L^1,L^2,\do)$) such that for $\d\in\{1,-1\}$, the
restriction of $\g$ to the generalized $\d$-eigenspace of $g$ has determinant $-1$.}
\nl
Define $(\o_t)$ by $\o_1=-1,\o_t=1$ for $t\in[2,\s+1]$. In our case we have $k\ge1$ hence $(\o_t)\in\ci$. 
Let $V'=\sum_{i\in\ZZ}L^1_i$, $V''=\sum_{t\in[2,\s+1],i\in\ZZ}L^t_i$. By 1.18, we have $V=V'\op V''$ 
(orthogonal direct sum). Define $\g\in I'$ by $\ps(\g)=(\o_t)$ (notation of 1.20). Then $\g$ acts as 
identity on $V''$ and as $-1$ times the identity on $V'$. It is enough to prove that the restriction of 
$\g$ to the generalized $\d$-eigenspace of $g_{V'}$ has determinant $-1$ or that this generalized 
$\d$-eigenspace has odd dimension. But this dimension is $a_1$ (if $\d=1$) and $b_1$ (if $\d=-1$) and 
$a_1,b_1$ are odd.

\subhead 1.23\endsubhead
In the setup of 1.1, assume that $\nn\ge3$ and $\e=1$. When $\nn$ is odd we assume that 
$\cc^V_{a_*,b_*}\sub Is(V)^0$ and let $\p:\G@>>>Is(V)^0$ be a surjective morphism of algebraic groups with 
kernel of order $2$ such that $\G$ is connected and simply connected. When $\nn$ is even let 
$\p:\G@>>>Is(V)$ be a surjective morphism of algebraic groups 
with kernel of order $2$ such that $\p\i(Is(V)^0)$ is connected and simply connected. 

Let $\boc$ be a $\G^0$-conjugacy class contained in $\p\i(\cc^V_{a_*,b_*})$. (If $a_1b_1>0$ we have
$\boc=\p\i(\cc^V_{a_*,b_*})$; if $a_1b_1=0$ there are two choices for $\boc$.)
For $\nn$ odd let $X$ be the set of all $(\tg,L^1,L^2,\do,L^{\s+1})$ where $\tg\in\boc$ and 
$(\p(\tg),L^1,L^2,\do,L^{\s+1})\in\tcc^V_{a_*,b_*}$. For $\nn$ even let $X$ be the set of all
$(\tg,L^1,L^2,\do,L^{\s})$ where $\tg\in\boc$ and $(\p(\tg),L^1,L^2,\do,L^\s)\in\tcc^V_{a_*,b_*;\co}$.
Note that $X\ne\emp$. Now $\G^0$ acts on $X$ by
$$\g:(\tg,L^1,L^2,\do,L^{\s+\k})\m(\g\tg\g\i,\p(\g)L^1,\p(\g)L^2,\do,\p(\g)L^{\s+\k}).$$
We show:

(a) {\it This action is transitive.}
\nl
If $a_1b_1=0$, then (a) follows trivially from 1.21(d), 1.22(a). Assume now that $a_1b_1>0$. Let 
$(\tg,L^1,L^2,\do,L^{\s+\k})\in X$ and let $c$ be the nontrivial element in $\ker\p$. Let $g=\p(\tg)$. We 
define $\g$ in terms of $(g,L^1,L^2,\do,L^{\s+\k})$ as in 1.21(f) or 1.22(c). Let $\ti\g\in\p\i(\g)$. Since 
$\g g\g\i=g$ we see that either $\ti\g\tg\ti\g\i=\tg$ or $\ti\g\tg\ti\g\i=c\tg$. In the first case $\ti\g$ 
is in the centralizer in $\G^0$ of $\tg_s$ (the semisimple part of $\tg$). This centralizer is a connected 
algebraic group (by a result of Steinberg). Thus its image under $\p$ is connected hence it is contained in 
the connected centralizer of $g_s$ (the semisimple part of $g$) in $Is(V)^0$. Thus $\g=\p(\ti\g)$ is 
contained in the connected centralizer of $g_s$ in $Is(V)^0$. But then the restriction of $\g$ to the 
$1$-eigenspace of $g_s$ would have determinant $1$, contradicting the choice of $\g$. We see that we must 
have 

(b) $\ti\g\tg\ti\g\i=c\tg$.
\nl
Using 1.21(d), 1.22(a), we see that any $\G^0$-orbit on $X$ contains either $(\tg,L^1,L^2,\do,L^{\s+\k})$
or $(c\tg,L^1,L^2,\do,L^{\s+\k})$. From (b) and the definition of $\ti\g$ we see that the action of $\tg$ 
takes $(\tg,L^1,L^2,\do,L^{\s+\k})$ to $(c\tg,L^1,L^2,\do,L^{\s+\k})$. This shows that (a) holds.

\subhead 1.24\endsubhead
As in \cite{\WEU, \S3}, \cite{\WEIII, \S3} we see that 1.23 (resp. 1.1) implies that Theorem 0.3 holds when 
$G$ is $\G$ in 1.23 (resp. $G=Is(V)$ with $\nn\ge2$, $\e=-1$).

\head 2. Bilinear forms\endhead
\subhead 2.0\endsubhead
For any subset $S$ of $\ZZ$ we write $S''=S\cap(2\ZZ)$, $S'=S\cap(2\ZZ+1)$.

Let $V$ be a $\kk$-vector space of finite dimension $n$. Let $(,):V\T V^*@>>>\kk$ be the obvious pairing. 
Let $G_V=GL(V)$ and let $G_V^1$ be the set of all vector space isomorphisms $V@>\si>>V^*$. Note that an
element of $G_V^1$ can be viewed as a bilinear form $V\T V@>>>\kk$. For $\g\in G_V$ we define 
$\che\g\in G_{V^*}$ by $(\g(x),\che\g(\x))=(x,\x)$ for all $x\in V,\x\in V^*$. For $g\in G^1_V$ we define 
$\chg\in G^1_{V^*}$ by $(\chg z',gz)=(z,z')$ for any $z\in V,z'\in V^*$. There is a well defined group 
structure on $G:=G_V\sqcup G_V^1$ denoted by $*$ such that for $\g,\g'$ in $G_V$ and $g,g'$ in $G^1_V$ we 
have
$$\g*\g'=\g\g'\in G_V;\qua \g*g'=\che\g g'\in G^1_V;\qua g*g'=\chg g'\in G_V;\qua g*\g'=g\g'\in G^1_V.$$
Now let $g\in G^1_V$. For $i\in\ZZ$ let $g^{*i}$ be the $i$-th power of 
$g$ for the multiplication $*$. In particular we have $g^{*2}=g*g=\chg g$. For $i\in\ZZ''$ we have 
$g^{*i}\in G_V$. For $i\in\ZZ'$ we have $g^{*i}\in G^1_V$. For any $z\in V$ and $i\in\ZZ$ we set 
$z_i=g^{*i}z$; we have $z_i\in V$ if $i\in\ZZ''$ and $z_i\in V^*$ if $i\in\ZZ'$. Similarly, for any line $L$
in $V$ and $i\in\ZZ$ we set $L_i=g^{*i}L$; this is a line in $V$ if $i\in\ZZ''$ and a line in $V^*$ if 
$i\in\ZZ'$.

For any $z,z'$ in $V$ and any $i\in\ZZ'',j\in\ZZ',k\in\ZZ''$, we show:
$$(z_{i+k},z'_{j+k})=(z_i,z'_j),\tag a$$
$$(z_i,z'_j)=(z'_{-i},z_{-j}).\tag b$$
Indeed, we have 
$$\align&(z_i,z'_j)=(z_i,gz'_{j-1})=(z'_{j-1},(\chg)\i z_i)=(z'_{j-1},g(\chg g)\i z_i)\\&=
(z'_{j-1},gz_{i-2})=(z'_{j-1},z_{i-1}).\tag c\endalign$$
Repeating this we get $(z'_{j-1},z_{i-1})=(z_{i-2},z'_{j-2})$. Combining with (c) we get 
$(z_i,z'_j)=(z_{i-2},z'_{j-2})$; hence $(z_i,z'_j)=\ph(i-j)$ where $\ph:\ZZ'@>>>\kk$; by (c) we have 
$(z'_i,z_j)=\ph(j-i)$ for $i\in\ZZ'',j\in\ZZ'$. In particular, (a),(b) hold.

\mpb

Let $a_1\ge a_2\ge\do$, $b_1\ge b_2\ge\do$ be two sequences of integers $\ge0$ in $\NN$ such that 

if $i\ge1$, $a_i=a_{i+1}$ then $a_{i+1}=0$;

if $i\ge1$, $b_i=b_{i+1}$ then $b_{i+1}=0$;

if $a_i>0$, then $a_i\in\ZZ'$;

if $b_i>0$, then $b_i\in\ZZ''$;

$(a_1+a_2+\do)+(b_1+b_2+\do)=n$.
\nl
It follows that $a_i=0$ for large $i$ and $b_i=0$ for large $i$. Define $k\ge0$ by 
$\{i\ge1;a_ib_i>0\}=[1,k]$. We define $p_i\in\NN$ for $i\ge1$ as follows. If $i\in[1,k]$, we have 
$p_i=(a_i+b_i+1)/2$. If $i>k$ we define $p_i$ by requiring that for $s=1,3,5,\do$ we have:
$$(p_{k+s},p_{k+s+1})=(b_{k+s}/2,(b_{k+s+1}+2)/2)\text{ if }b_{k+s}>0;$$
$$(p_{k+s},p_{k+s+1})=((a_{k+s}+1)/2,(a_{k+s+1}+1))/2)\text{ if }a_{k+s}>0,a_{k+s+1}>0;$$
$$(p_{k+s},p_{k+s+1})=((a_{k+s}+1)/2,0)\text{ if }a_{k+s}>0,a_{k+s+1}=0;$$
$$(p_{k+s},p_{k+s+1})=(0,0)\text{ if }a_{k+s}=a_{k+s+1}=0.$$
We define $\s$ as follows. If $n=0$ we set $\s=0$. If $n\ge1$ let $\s$ be the largest $i$ such that $p_i>0$.
We have $p_1\ge p_2\ge\do\ge p_\s$ and
$$(2p_1-1)+(2p_2-1)+\do+(2p_\s-1)=n.$$

\mpb

Let $\cc^V_{a_*,b_*}$ be the set of all $g\in G^1_V$ such that $g^{*4}\in G_V$ is unipotent and such that on 
the generalized $1$-eigenspace of $g^{*2}$, $g^{*2}$ has Jordan blocks of sizes given by the nonzero numbers
in $a_1,a_2,\do$ and on the generalized $(-1)$-eigenspace of $g^{*2}$, $-g^{*2}$ has Jordan blocks of sizes 
given by the nonzero numbers in $b_1,b_2,\do$. 

For $g\in\cc^V_{a_*,b_*}$ let $\tcc^V_{g;a_*,b_*}$ be the set consisting of all $L^1,L^2,\do,L^\s$ where 
$L^t (t\in[1,\s])$ are lines in $V$ (the upper scripts are not powers) such that for $i\in\ZZ'',j\in\ZZ'$ we
have:

$(L^t_i,L^t_j)=0$ if $i-j\in[-2p_t+3,2p_t-3]'$, 

$(L^t_i,L^t_j)\ne0$ if $|i-j|=2p_t-1$ ($t\in[1,\s]$);

$(L^r_i,L^t_j)=0$ if $j-i\in[1-2p_r,4p_t-2p_r-3]'$, $1\le t<r\le\s$.
\nl
Here $L^t_i=g^{*i}L^t$. We then have:

(d) {\it $V=\op_{t\in[1,\s],i\in[0,2p_t-2]}L^t_i$.}
\nl  
(See \cite{\WEIII, 4.8(a)}.) Let $\tcc^V_{a_*,b_*}$ be the set of all $(g,L^1,L^2,\do,L^\s)$ such that
$g\in\cc^V_{a_*,b_*}$ and $(L^1,L^2,\do,L^\s)\in\tcc^V_{g;a_*,b_*}$.

Note that $G_V$ acts on $G^1_V$ by "twisted conjugation" that is by $\g:g\m\che\g g\g\i$. Also $G_V$ acts on
$\tcc^V_{a_*,b_*}$ by 
$$\g:(g,L^1,L^2,\do,L^\s)\m(\che\g g\g\i,\g(L^1),\g(L^2),\do,\g(L^\s)).\tag e$$
Now let $\ci$ be the subgroup of $\prod_{t\in[1,\s]}\{1,-1\}$ consisting of all $(\o_t)_{t\in[1,\s]}$ such 
that $\o_t=\o_{t+1}$ for any $t$ such that $\{t,t+1\}\sub[k+1,\s],t=k+1\mod2$, $b_t>0$. Thus $\ci$ is a 
finite elementary abelian $2$-group.

The folowing is the main result of this section.

\proclaim{Theorem 2.1} (a) $\tcc^V_{a_*,b_*}$ is nonempty;

(b) the action 2.0(e) of $G_V$ on $\tcc^V_{a_*,b_*}$ is transitive;

(c) the isotropy group in $G_V$ at any point of $\tcc^V_{a_*,b_*}$ is  canonically isomorphic to $\ci$.
\endproclaim 

\subhead 2.2\endsubhead
Let $a\in\NN',b\in\NN'',p\in\NN_{>0}$ be such that $a+b=2p-1$. For $e\in\NN''$ we define $n_e\in\ZZ$ by
$$(1-T^2)^a(1+T^2)^b=\sum_{e\in\NN''}n_eT^e.$$
We have $n_0=1,n_{4p-2-e}=-n_e$, $n_e=0$ if $e>4p-2$. We define $x_e\in\ZZ$ for $e\in\NN''$ by $x_0=1$ and
$$n_0x_e+n_2x_{e-2}+\do+n_ex_0=0\text{ for }e\ge2.\tag a$$  
For $h\in\ZZ'$ we set $x'_h=0$ if $|h|<2p-1$, $x'_h=x_{|h|-2p+1}$ if $|h|\ge2p-1$. We show:
$$\sum_{e\in\NN''}n_ex'_{e-j-1}=0\text{ for any }j\in[0,4p-4]''.\tag b$$
Assume first that $j\in[2p,4p-4]''$. We have 
$$e-j-1\le e-2p-1\le4p-2-2p-1\le2p-3.$$
 Hence we can assume that 
$e-j-1\le-2p+1$ so that $x'_{e-j-1}=x_{j+1-e-2p+1}$ and we must show that 
$$\sum_{e;e\le j+1-2p+1}n_ex_{j+1-e-2p+1}=0.$$
This holds since $j+1-2p+1\ge2$. Assume next that 
$j\in[0,2p-4]''$. We have $e-j-1\ge e-2p+4-1\ge-2p+3$. Hence we can assume that $e-j-1\ge2p-1$ so that 
$x'_{e-j-1}=x_{e-j-1-2p+1}$ and we must show that 
$$\sum_{e;e\ge j+1+2p-1}n_ex_{e-j-1-2p+1}=0$$
that is,
$$\sum_{e;e\ge j+1+2p-1}n_{4p-2-e}x_{e-j-1-2p+1}=0,$$ 
that is,
$$\sum_{e';4p-2-e'\ge j+1+2p-1}n_{e'}x_{4p-2-e'-j-1-2p+1}=0,$$
that is,
$$\sum_{e';e'\le 2p-2-j}n_{e'}x_{2p-2-e'-j}=0,$$
and this holds since $2p-2-j\ge2$. Assume next that $j=2p-2$. In the sum over $e$ we can assume that 
$e-j-1\ge2p-1$ or $e-j-1\le-2p+1$, that is $e\ge4p-2$ or $e\le0$. Thus $e=0$ or $e=4p-2$. Thus the sum is
$n_0x'_{-2p+1}+n_{4p-2}x'_{2p-1}=n_0+n_{4p-2}=0$.

\subhead 2.3\endsubhead
In the setup of 2.2 let $V$ be a $\kk$-vector space of dimension $2p-1$. Assume that we are given a basis 
$\{w_i;i\in[0,4p-4]''\}$ of $V$. Let $\{w_i;i\in[1,4p-3]'\}$ be the basis of $V^*$ such that
$$(w_i,w_j)=x'_{i-j}=x'_{j-i}\text{ if }i\in[0,4p-4]'',j\in[1,4p-3]'.$$
Thus $(w_i,w_j)=0$ if $|i-j|<2p-1$. We define $g\in G^1_V$ by $gw_i=w_{i+1}$ for $i\in[0,4p-4]''$. Let 
$\chg\in G^1_{V^*}$ be as in 2.0. We have 
$$\chg w_i=w_{i+1}\text{ if }i\in[1,4p-5]';$$
we must check that $(w_{i+1},w_{j+1})=(w_j,w_i)$ for $i\in[1,4p-3]',j\in[0,4p-4]''$; we use that 
$|i+1-(j+1)|=|j-i|$. 

We show:
$$\chg w_{4p-3}=\sum_{i\in[0,4p-4]''}n_iw_i$$ 
that is,
$$\sum_{i\in[0,4p-4]''}n_i(w_i,w_{j+1})=(w_j,w_{4p-3})\text{ for any }j\in[0,4p-4]'',$$
that is,
$$\sum_{i\in[0,4p-4]''}n_ix'(i-j-1)=x'(4p-3-j)\text{ for any }j\in[0,4p-4]'',$$
that is,
$$\sum_{i\in[0,4p-2]''}n_ix'(i-j-1)=0\text{ for any }j\in[0,4p-4]''.$$
This has been seen in 2.2(b).

We have $g^{*2}(w_i)=w_{i+2}$ for $i\in[0,4p-6]''$, $g^{*2}(w_{4p-4})=\sum_{i\in[0,4p-4]''}n_iw_i$. Hence 
$(g^{*2}-1)^a(g^{*2}+1)^b=0$ on $V$. Indeed this holds on $w_0$ and then it holds automatically on 
$w_i,i\in[0,4p-4]''$. Now $g^{*2}\in G_V$ is regular in the sense of Steinberg and satisfies 
$(g^{*2}-1)^a(g^{*2}+1)^b=0$ on $V$. Hence $V=V^+\op V^-$ where $g^{*2}$ acts on $V^+$ as a single unipotent
Jordan block of size $a$ and $-g^{*2}$ acts on $V^-$ as a single unipotent Jordan block of size $b$.

It follows that, if $L$ is the line in $V$ spanned by $w_0$ and $a_*=(a,0,0,\do)$, $b_*=(b,0,0,\do)$, then 
$(g,L)\in\tcc^V_{a_*,b_*}$; in particular, $\tcc^V_{a_*,b_*}\ne\emp$.

We now consider a variant of the situation above. Let $V'$ be a $\kk$-vector space of dimension $2p-1$ with 
a given element $g\in G^1_{V'}$ such that $g^{*4}=1$, on the generalized $1$-eigenspace of $g^{*2}$, 
$g^{*2}$ is a single unipotent Jordan block of size $a$ and on the generalized $(-1)$-eigenspace of 
$g^{*2}$, $-g^{*2}$ is a single unipotent Jordan block of size $b$. Moreover we assume that we are given 
$w\in V'$ such (with notation of 2.0) we have 
$$\align&(w_i,w_j)=0\text{ if }i\in\ZZ'',j\in\ZZ',|i-j|<2p-1\text{ and }\\&
(w_i,w_j)=1\text{ if }i\in\ZZ'',j\in\ZZ',|i-j|=2p-1.\endalign$$
We show:

(a) {\it for any $i\in\ZZ'',j\in\ZZ'$ we have $(w_i,w_j)=x'_{i-j}$.}
\nl
We can assume that $i=0$ and $j\ge1$. The equality in (a) is already known if $j\le2p-1$. It is enough to 
show that $(w_0,w_{2p-1+2t})=x_{2t}$ for $t\in\NN$. We argue by induction on $t$; for $t=0$ the result is 
already known. Now assume that $t\ge1$. Applying $(g^{*2}-1)^a(g^{*2}+1)^b=0$ to $w_{2t-2p+2}$ we obtain 
$\sum_{e\in[0,4p-2]''}n_ew_{2t-2p+2+e}=0$. Taking $(,w_1)$ we obtain
$$\sum_{e\in[0,4p-2]''}n_e(w_{2t-2p+2+e},w_1)=0$$ 
that is, 
$$\sum_{e\in[0,4p-2]''}n_e(w_0,w_{2t-2p+1+e})=0.$$
For $e$ in the sum we have $2t-2p+1+e\ge-2p+3$; hence we can assume that we have $2t-2p+1+e\ge2p-1$. Thus 
$$\sum_{e\in[0,4p-2]'';2t-2p+1+e\ge2p-1}n_e(w_0,w_{2t-2p+1+e})=0.$$
By the induction hypothesis this implies
$$\sum_{e\in[0,4p-4]'';2t-2p+1+e\ge2p-1}n_e x_{2t-4p+2+e}-(w_0,w_{2t-2p+1+e})=0.$$
It is then enough to show that
$$\sum_{e\in[0,4p-4]'';2t-2p+1+e\ge2p-1}n_e x_{2t-4p+2+e}-x_{2t}=0,$$
or that
$$\sum_{e\in[0,4p-2]'';2t-2p+1+e\ge2p-1}n_{4p-2-e}x_{2t-4p+2+e}=0,$$
or that 
$$\sum_{h,h'\in\NN'';h+h'=2t}n_hx_{h'}=0.$$
But this holds by the definition of $x_e$ since $2t\ge2$.

\subhead 2.4\endsubhead   
Let $p\in\NN_{>0}$. We define $n_e$ for $e\in\NN''$ by $n_e=\bin{2p}{e/2}$. We define $x_e$ for $e\in\NN''$ 
by $x_0=1,x_2=-(2p+1)$, and $n_0x_e+n_2x_{e-2}+\do+n_ex_0=0$ for $e\ge4$. For $e=2$ we have 
$$n_0x_e+n_2x_{e-2}+\do+n_ex_0=n_0x_2+n_2x_0=-(2p+1)+2p=-1.$$
For $d\in\ZZ'$ we set $\ph_p(d)=0$ if $|d|<2p-1$, $\ph_p(d)=x_{|d|-2p+1}$ if $|d|\ge2p-1$. We show for any 
$h\in\ZZ'$:
$$\sum_{e\in[0,4p]''}n_e\ph_p(e+h)=0\tag a$$ 
Assume that $h\le-1$. We set $h=-j-1$ so that $j\in\NN''$. Assume first that $j\ge2p+2$. We have 
$e-j-1\le e-2p-2-1\le 4p-2p-2-1\le 2p-3$. Hence we can assume that $e-j-1\le-2p+1$ so that
$\ph_p(e-j-1)=x_{j+1-e-2p+1}$ and we must show
$$\sum_{e\in\NN'';e\le j+1-2p+1}n_ex_{j+1-e-2p+1}=0.$$
This holds since $j+1-2p+1\ge4$.

Assume next that $j\le2p-4$. We have $e-j-1\ge e-2p+4-1\ge-2p+3$. Hence we can assume that $e-j-1\ge2p-1$ so
that $\ph_p(e-j-1)=x_{e-j-1-2p+1}$ and we must show:
$$\sum_{e\in\NN'';e\ge j+1+2p-1}n_ex_{e-j-1-2p+1}=0,$$
that is
$$\sum_{e\in\NN'';e\ge j+1+2p-1}n_{4p-e}x_{e-j-1-2p+1}=0,$$
that is
$$\sum_{e'\in\NN'';4p-e'\ge j+1+2p-1}n_{e'}x_{4p-e'-j-1-2p+1}=0,$$
that is
$$\sum_{e'\in\NN;e'\le2p-j}n_{e'}x_{2p-e'-j}=0$$
and this holds since $2p-j\ge4$.  

Assume next that $j=2p-2$. In the sum we can assume that $e-j-1\ge2p-1$ or $e-j-1\le-2p+1$ that is 
$e\ge4p-2$ or $e\le 0$. Thus $e=0$ or $e=4p-2$ or $e=4p$. Thus the sum is 
$$n_0\ph_p(-2p+1)+n_{4p-2}\ph_p(2p-1)+n_{4p}\ph_p(2p+1)=x_0+2px_0+x_2=x_2+2p+1=0.$$
Assume next that $j=2p$. In the sum we can assume that $e-j-1\ge2p-1$ or $e-j-1\le-2p+1$ that is $e\ge4p$ or 
$e\le2$. Thus $e=0,2$ or $4p$. Thus the sum is 
$$n_0\ph_p(-2p-1)+n_2\ph_p(-2p+1)+n_{4p}\ph_p(2p-1)=x_2+2p+1=0.$$
Thus the desired formula holds when $h\le-1$. Now assume that $h\ge1$. We have 
$$\align&\sum_{e\in\NN''}n_e\ph_p(e+h)=\sum_{e\in[0,4p]''}n_{4p-e}\ph_p(e+h)\\&=
\sum_{e\in[0,4p]''}n_e\ph_p(4p-e+h)=\sum_{e\in[0,4p]''}n_e\ph_p(-4p+e-h)\endalign$$ 
and this is $0$ by the first part of the proof since $-4p-h\le-1$.

\mpb

We have 
$$\sum_{e\in\NN'',j\in\NN''}n_eT^ex_jT^j=1-T^2,$$ 
hence 
$$(1+T^2)^{2p}\sum_{j\in\NN''}x_jT^j=1-T^2$$
and
$$\sum_{j\in\NN''}x_jT^j=(1-T^2)(1+T^2)^{-2p}=(1-T^2)(\sum_{k\ge0}(-1)^k\bin{2p-1+k}{2p-1}T^{2k}).$$
Thus,
$$\align&x_{2k}=(-1)^k\bin{2p-1+k}{2p-1}-(-1)^{k-1}\bin{2p-1+k-1}{2p-1}\\&=
(-1)^k((2p-1+k)!/k!+(2p+k-2)!/(k-1)!)/(2p-1)!\\&=(-1)^k (2p-2+k)!/k!(2p-1+k+k)/(2p-1)!\\&=
(-1)^k(2p+2k-1)(2p-2+k)(2p-2+k-1)\do(k+1)(2p-1)!\i.\endalign$$
We show for any $h\in\ZZ'$:
$$\ph_p(h)=(-1)^{(h+2p+1)/2}2h(h+2p-3)(h+2p-5)\do(h-2p+3)(4p-2)!!\i\tag b$$ 
where 
$$(4p-2)!!:=2\T4\T\do\T(4p-2)=2^{2p-1}(2p-1)!.$$
Assume first that $h=2d+1\ge2p-1$. We have
$$\align&\ph_p(h)=x_{2d+1-2p+1}=x_{2d-2p+2}\\&=
(-1)^{d-p+1}(2p+2d-2p+2-1)(2p-2+d-p+1)\\& \T(2p-2+d-p+1-1)\do(d-p+2)(2p-1)!\i\\&
=(-1)^{d-p+1}(2d+1)(p+d-1)(p+d-2)\do(d-p+2)(2p-1)!\i\endalign$$
so that the result holds in this case. Now both sides of (a) are invariant under $h\m-h$. Hence (a) also
holds if $h\le-2p+1$. If $h\in[-2p+3,2p-3]'$, both sides of (a) are zero. Hence (a) holds for any $h\in\ZZ'$.

In particular we have $\ph_p(2p+1)=-(2p+1)$.

\subhead 2.5\endsubhead
In the setup of 2.4, let $E$ be a $\kk$-vector space of dimension $2p$. Assume that we are given a basis
$\{w_i;i\in[0,4p-2]''\}$ of $E$. We define a basis $\{w_i;i\in[1,4p-1]'\}$ of $E^*$ by
$$(w_i,w_j)=\ph_p(i-j)=\ph_p(j-i)\text{ for }i\in[0,4p-2]'',j\in[1,4p-1]'.$$
Thus $(w_i,w_j)=0$ if $i\in[0,4p-2]'',j\in[1,4p-1]'$, $|i-j|<2p-1$. We define $g\in G^1_E$ by 
$gw_i=w_{i+1}$ for $i\in[0,4p-2]''$. We have 
$$\chg w_i=w_{i+1}\text{ if }i\in[1,4p-3]';$$
we must check that $(w_{i+1},w_{j+1})=(w_j,w_i)$ for $i\in[1,4p-1]',j\in[0,4p-2]''$; we use that 
$|i+1-(j+1)|=|j-i|$. 

We show:
$$\chg w_{4p-1}=-\sum_{i\in[0,4p-4]''}n_iw_i.$$
We must show for any $j\in[0,4p-2]''$ that 
$$-\sum_{i\in[0,4p-2]''}n_i(w_i,w_{j+1})=(w_j,w_{4p-1})$$
that is,
$$-\sum_{i\in[0,4p-2]''}n_i\ph_p(i-j-1)=\ph_p(4p-1-j),$$
that is,
$$\sum_{i\in[0,4p]''}n_i\ph_p(i-j-1)=0;$$
note that $n_{4p}=-1$. This has been seen in 2.4(a).

We have
$$g^{*2}(w_i)=w_{i+2}\text{ for }i\in[0,4p-4],$$
$$g^{*2}(w_{4p-2})=-\sum_{i\in[0,4p-2]''}n_iw_i.$$
Hence 
$$(g^{*2}+1)^{2p}=0\text{ on }E.\tag a$$ 
Indeed this holds on $w_0$ and then it holds automatically on $w_i,i\in[0,4p-2]''$. Now $g^{*2}\in GL(E)$
is regular in the sense of Steinberg and satisfies (a). Hence $-g^{*2}$ acts on $E$ as a single unipotent 
Jordan block of size $2p$.

\subhead 2.6\endsubhead
For $i\in\ZZ$ we write $w_i$ instead of $(w_0)_i$. This agrees with our earlier notation for $w_i$ when 
$i\in[0,4p-1]$. We show:
$$(w_i,w_j)=\ph_p(i-j)=\ph_p(j-i)\text{ for any }i\in\ZZ'',j\in\ZZ'.\tag a$$
By 2.0(a) there exists a function $f:\ZZ'@>>>\kk$ such that $(w_i,w_j)=f(i-j)$ for any $i\in\ZZ'',j\in\ZZ'$.
We must show that $f(h)=\ph_p(h)$ for $h\in\ZZ'$. We set $f'(h)=f(h)-\ph_p(h)$. We must show that $f'(h)=0$
for all $h\in\ZZ'$. This is clearly true when $h\in[-2p+1,2p-1]'$. Applying $\sum_{e\in[0,4p]''}n_eg^{*e}=0$
to $w_i, i\in\ZZ''$, we deduce 
$$\sum_{e\in[0,4p]''}n_ew_{i+e}=0;$$
hence 
$$\sum_{e\in[0,4p]''}n_e(w_{i+e},w_j)=0\text{ for }i\in\ZZ'',j\in\ZZ'.$$
Thus, $\sum_{e\in[0,4p]''}n_ef(i-j+e)=0$ for $i\in\ZZ'',j\in\ZZ'$ and $\sum_{e\in[0,4p]''}n_ef(h+e)=0$ for 
$h\in\ZZ''$. Combining this with $\sum_{e\in[0,4p]''}n_e\ph_p(h+e)=0$ for $h\in\ZZ''$, see 2.4(a), we deduce
$\sum_{e\in[0,4p]''}n_ef'(h+e)=0$ for $h\in\ZZ''$. We show that $f'(h)=0$ for $h\ge2p-1$ by induction on $h$.
For $h=2p-1$ this is already known. Now assume that $h\ge2p+1$. We have 
$\sum_{e\in[0,4p]''}n_ef'(h+e-4p)=0$. If $e\in[0,4p-2]''$ we have $h+e-4p\in[-2p+1,h-2]$ hence $f'(h+e-4p)=0$
and the sum over $e$ becomes $n_{4p}f'(h)=0$ so that $f'(h)=0$. This completes the induction. We now show 
that $f'(h)=0$ for $h\le-2p+1$ by descending induction on $h$. For $h=-2p+1$ this is known. Now assume that 
$h\le-2p-1$. If $e\in[2,4p]''$ we have $h+e\in[h+2,2p-1]$ hence $f'(e+h)=0$ and the equation
$\sum_{e\in[0,4p]''}n_ef'(h+e)=0$ becomes $n_0f'(h)=0$ so that $f'(h)=0$. This completes the descending 
induction and completes the proof of (a).

\subhead 2.7\endsubhead
We preserve the setup of 2.5. Let $\tw$ be a nonzero vector in $E$ such that

(a) $(\tw,w_i)=0$ for $i\in[1,4p-3]'$. 
\nl
Note that $\tw$ is uniquely determined up to a nonzero scalar. Then $\tw_i$ is defined for any $i\in\ZZ$
as in 2.0; in particular, $\tw_0=\tw,\tw_1=g\tw$. We have

(b) $(w_i,\tw_1)=0$ for $i\in[2,4p-2]''$.
\nl
Indeed, using 2.0(a),(b) we have $(w_i,\tw_1)=(\tw_{-i},w_{-1})=(\tw_0,w_{i-1})$ and this is zero since 
$i-1\in[1,4p-3]'$. 

We show that $(\tw_0,\tw_1)\ne0$. Let $E_1$ be the span of $\{w_i;i\in[2,4p-2]''\}$ and let $E'_1$ be the 
span of $\{w_i;i\in[1,4p-3]'\}$. The canonical pairing $(,):E\T E^*@>>>\kk$ restricts to a nondegenerate 
pairing $E_1\T E'_1@>>>\kk$ (by the formulas for $(w_i,w_j)$ in 2.5). Since $\tw_0$ is in the annihilator of
$E'_1$ in $E$, it follows that $\tw_0\n E_1$. Since $\tw_1$ is in the annihilator of $E_1$ in $E^*$, it 
follows that $\tw_0$ is not in the annihilator of $\tw_1$ in $E$. The claim follows. 

If $\tw$ is replaced by $a\tw$ with $a\in\kk^*$, then $(\tw_0,\tw_1)$ is replaced by $a^2(\tw_0,\tw_1)$ 
which, for a suitable $a$, is equal to $1$. Thus we can assume that
$$(\tw_0,\tw_1)=1.\tag c$$
Then $\tw_0$ is uniquely determined up to multiplication by $\pm1$. We have
$$\tw_0=\sum_{i\in[0,4p-2]''}c_iw_i$$
where $c_i\in\kk$ are uniquely determined. Since $\tw_0\n E_1$ we see that $c_*:=c_{4p-2}\ne0$. We set 
$\bc_i=c_ic_*\i\in\kk$. Note that $\bc_{4p-2}=1$. We have the following result.
$$\bc_i=-(n_0+n_2+\do+n_i)\text{ if }i\in[0,2p-2]'',\tag d$$ 
$$\bc_i=(n_0+n_2+\do+n_{4p-2-i})\text{ if }i\in[2p,4p-2]'',\tag e$$
$$c_*=\pm2^{-p}.\tag f$$
We can rewrite (a) as follows.
$$\sum_{i\in[0,4p-2]''}\bc_i\ph_p(i-h)=0\text{ for }h\in\in[1,4p-3]'.\tag $(*)$ $$
If $h=2p-1$ then $(*)$ is $\bc_0+1=0$. If $h\in[2p+1,4p-3]'$ then $(*)$ is 
$$\sum_{i\in[0,h-2p+1]''}\bc_i\ph_p(i-h)=0.$$ If $h\in[1,2p-3]'$ then $(*)$ is 
$$\sum_{i\in[h+2p-1,4p-2]''}\bc_i\ph_p(i-h)=0.$$ To prove (d),(e) it is enough to show:
$$-\sum_{i\in[0,h-2p+1]''}(n_0+n_2+\do+n_i)\ph_p(i-h)=0\text{ if }h\in[2p+1,4p-3]',\tag d${}'$ $$ 
$$\sum_{i\in[h+2p-1,4p-2]''}(n_0+n_2+\do+n_{4p-2-i})\ph_p(i-h)=0\text{ if }h\in[1,2p-3]'.\tag e${}'$ $$ 
We rewrite equation (e${}'$) using $i\m4p-2-i$ and $h\m4p-2-h$ as
$$\sum_{i\in[0,h-2p+1]''}(n_0+n_2+\do+n_i)\ph_p(h-i)=0\text{ if }h\in[2p+1,4p-3]'$$ 
which is the same as (d${}'$). Thus it is enough to prove (d${}'$). We argue by induction on $h$. If 
$h=2p+1$, equation (d${}$') is
$$n_0\ph_p(-2p-1)+(n_0+n_2)\ph_p(-2p+1)=0$$
that is $-(2p+1)+(1+2p)=0$, which is correct. If $h\ge2p+3$ we have 
$$\sum_{i\in[0,h-2p+1]''}n_ix_{h-i-2p+1}=0$$
since $h-2p+1\ge4$. Hence in this case (d${}'$) is equivalent to 
$$\sum_{i\in[2,h-2p+1]''}(n_0+n_2+\do+n_{i-2})\ph_p(i-h)=0$$ 
which is the same as equation (d${}'$) with $h$ replaced by $h-2$ (this holds by the induction hypothesis). 
This proves (d),(e).

The equation $(\tw_0,\tw_1)=1$ can be written as
$$1=(\tw_0,\sum_{i\in[0,4p-2]''}c_iw_{i+1})=(\tw_0,c_{4p-2}w_{4p-1})$$
that is,
$$1=c_*(\tw_0,w_{4p-1}).\tag g$$
We deduce that
$$1=c_*\sum_{i\in[0,4p-2]''}c_i(w_i,w_{4p-1})$$ 
that is,
$$c_*^{-2}=\sum_{i\in[0,4p-2]''}\bc_i\ph_p(4p-1-i).$$
We have $4p-i-1\ge-2p+3$ hence we can assume $4p-i-1\ge2p-1$. Thus
$$\align&c_*^{-2}=\sum_{i\in[0,2p]''}\bc_i\ph_p(4p-1-i)\\&
=-\sum_{i\in[0,2p-2]''}(n_0+n_2+\do+n_i)\ph_p(4p-1-i)+n_0+n_2+\do+n_{2p-2}\\&
=-\sum_{i\in[0,2p-2]''}(n_0+n_2+\do+n_i)x_{2p-i}+n_0+n_2+\do+n_{2p-2}.\endalign$$
Thus,
$$\align&c_*^{-2}=-\sum_{i\in\NN'',j\in\NN'';i+j\le 2p,j\ge2}n_ix_j+n_0+n_2+\do+n_{2p-2}\\&
=-\sum_{i\in\NN'',j\in\NN'';i+j\le 2p}n_ix_j+\sum_{i\in\NN'';i\le 2p}n_i+(n_0+n_2+\do+n_{2p-2})\\&
=-\sum_{k\in[0,2p]''}\sum_{i\in\NN'',j\in\NN'';i+j=k}n_ix_j+
\sum_{i\in\NN'';i\le 2p}n_i+(n_0+n_2+\do+n_{2p-2})\\&
=-\sum_{k\in[0,2p]'';k=0,2}\sum_{i\in\NN'',j\in\NN'';i+j=k}n_ix_j+
\sum_{i\in\NN'';i\le 2p}n_i+(n_0+n_2+\do+n_{2p-2})\\&
=-1+n_0x_2+n_2x_0+\sum_{i\in\NN'';i\le 2p}n_i+(n_0+n_2+\do+n_{2p-2})\\&
=-1-(n_2+1)+n_2+\sum_{i\in\NN'';i\le 2p}n_i+(n_0+n_2+\do+n_{2p-2})\\&=\sum_{i\in\NN'';i\le 2p}n_i
+(n_0+n_2+\do+n_{2p-2})\\&=n_0+n_2+\do+n_{2p}+n_{2p+2}+\do+n_{4p}=2^{2p}.\endalign$$
Thus $c_*^{-2}=2^{2p}$ and (f) follows.

If $\tw$ is replaced by $-\tw$ then $c_*$ is changed into $-c_*$. Hence $\tw$ can be chosen uniquely so that
$$c_*=2^{-p}.\tag f${}'$ $$

\subhead 2.8\endsubhead
We preserve the setup of 2.5. For $h\in\ZZ'$ we show
$$(\tw_0,w_h)=(-1)^{(h+1)/2}2^p(h-1)(h-3)\do(h-4p+3)(4p-2)!!\i\in2\ZZ.\tag a$$
We have $(\tw_0,w_h)=\sum_{i\in[0,4p-2]''}c_i\ph_p(i-h)$. Since $c_i=2^{-p}\bc_i$ it is enough to prove
$$\sum_{i\in[0,4p-2]''}\bc_i(-1)^{(h+1)/2}\ph_p(i-h)=2^{2p}(h-1)(h-3)\do(h-4p+3)(4p-2)!!\i.\tag b$$
It is also enough to prove this equality in $\ZZ$. For fixed $i$, $(-1)^{(h+1)/2}\ph_p(i-h)$ is a polynomial
in $h$ with rational coefficients of degree $\le2p-1$. Hence the left hand side of (b) is a polynomial in $h$
with rational coefficients of degree $\le2p-1$. Since $(\tw_0,w_h)=0$ for $h\in[1,4p-3]'$, this polynomial is
zero for $h\in[1,4p-3]'$ (that is for $2p-1$ values of $h$). It follows that 
$$(-1)^{(h+1)/2}\sum_{i\in[0,4p-2]''}\bc_i\ph_p(i-h)=a(h-1)(h-3)\do(h-4p+3)$$
for some rational number $a$. (The left hand side is $(-1)^{(h+1)/2}2^p(\tw_0,w_h)$.) For $h=4p-1$ we have 
$(\tw_0,w_h)=c_*\i=2^p$, see 2.7(g), hence 
$$2^{2p}=a(4p-2)(4p-4)\do2=a(4-2)!!$$
 that is, $a=2^{2p}(4p-2)!!\i$. It remains to show that
$$(-1)^{(h-1)/2}2^p(h-1)(h-3)\do(h-4p+3)(4p-2)!!\i\in2\ZZ.$$
Setting $h=2s+1$ it is enough to show that 
$$2^p(2s+1-1)(2s+1-3)\do(2s+1-4p+3)(4p-2)!!\i\in2\ZZ$$ 
or that
$$2^ps(s+1)\do(s-2p+2)(2p-1)!\i\in2\ZZ.$$
This is obvious since $p\ge1$.

\subhead 2.9\endsubhead
We preserve the setup of 2.5. We will show:
$$(\tw_0,\tw_h)=\sum_{k\in[1,p]}2^{2k-2}\ph_k(h)\in\kk\text{ for }h\in\ZZ';\tag a$$
$$(\tw_0,\tw_h)=1\text{ if }h\in[-2p+1,2p-1]';\tag b$$
$$(\tw_0,\tw_{2p+1})=1-2^{2p}.\tag c$$
We prove (a). We have
$$\align&(\tw_0,\tw_h)=\sum_{i\in[0,4p-2]''}c_i(\tw_0,w_{i+h})=
\sum_{i\in[0,4p-2]''}(-1)^{(i+h+1)/2}c_i2^p(i+h-1)\\&\T(i+h-3)\do(i+h-4p+3)\T(4p-2)!!\i.\tag d\endalign$$
Thus, (a) would follow from the equality
$$\align&\sum_{i\in[0,4p-2]''}(-1)^{i/2}\bc_i(i+h-1)(i+h-3)\do(i+h-4p+3)(4p-2)!!\i\\&=\sum_{k\in[1,p]}
(-1)^{(h+1)/2}2^{2k-2}\ph_k(h)\tag e\endalign$$
in $\kk$. It is enough to prove that (e) holds in $\ZZ$. We will do that assuming that (b) holds. Let 
$F_p(h)$ be the left hand side of (e). It can be viewed as a polynomial with rational coefficients in $h$ of
degree $\le2p-1$ in which the coefficient of $h^{2p-1}$ is
$$\align&(4p-2)!!\i\sum_{i\in[0,4p-2]''}\bc_i(-1)^{-i/2}\\&=
-(4p-2)!!\i\sum_{i\in[0,2p-2]''}(n_0+n_2+\do+n_i)(-1)^{i/2}\\&
+(4p-2)!!\i\sum_{i\in[2p,4p-2]''}(n_0+n_2+\do+n_{4p-2-i})(-1)^{i/2}\\&
=-(4p-2)!!\i\sum_{i\in[0,2p-2]''}(n_0+n_2+\do+n_i)(-1)^{i/2}\\&
+(4p-2)!!\i\sum_{i\in[0,2p-2]''}(n_0+n_2+\do+n_i)(-1)^{(4p-2-i)/2}\\&
=-2(4p-2)!!\i\sum_{i\in[0,2p-2]''}(n_0+n_2+\do+n_i)(-1)^{i/2}\\&
=-2(4p-2)!!\i(-1)^{p-1}(n_{2p-2}+n_{2p-6}+n_{2p-10}+\do)\\&
=-2(4p-2)!!\i(-1)^{p-1}2^{2p-2}\\&
=(-1)^p2^{2p-1}(4p-2)!!\i=(-1)^p(2p-1)!\i.\endalign$$
Thus, 
$$F_p(h)=(-1)^p(2p-1)!\i h^{2p-1}+\text{ lower powers of }h.$$
Note that $F_p(-h)=-F_p(h)$ for $h\in\ZZ'$. An equivalent statement is that 
$$(-1)^{(h+1)/2}(\tw_0,\tw_h)=-(-1)^{(-h+1)/2}(\tw_0,\tw_{-h})$$
which follows from
$(\tw_0,\tw_h)=(\tw_0,\tw_{-h})$, see 2.0. It follows that $F_p(-h)=-F_p(h)$ as polynomials in $h$. 
Specializing this for $h=0$ we see that 
$$F_p(0)=0.\tag g$$
In the case where $p=1$, from (f),(g) we see that $F_1(h)=-h$ so that (e) holds in this case (we have
$(-1)^{(h+1)/2}\ph_1(h)=-h$). We now assume that $p\ge2$. Now $F_p-F_{p-1}$ is a polynomial of degree $2p-1$
in $h$ whose value at $h\in[-2p+3,2p-3]'$ is $(-1)^{(h+1)/2}-(-1)^{(h+1)/2}=0$ (we use (b) for $p$ and 
$p-1$) and whose value at $0$ is $0$ (see (e)); moreover the coefficient of $h^{2p-1}$ in $F_p-F_{p-1}$ is 
$(-1)^p(2p-1)!\i$ (see (f)). It follows that $F_p-F_{p-1}=(-1)^{(h+1)/2}2^{2p-2}\ph_p(h)$. From this we see 
by induction on $p$ that (e) holds.

It remains to prove (b),(c). To prove (b) we can assume that $h\ge1$ (we use that 
$(\tw_0,\tw_h)=(\tw_0,\tw_{-h})$, see 2.0). Thus it is enough to prove (b) for $h\in[1,2p-1]'$ and (c). If 
$h=1$, (b) holds by the definition of $\tw_0$. Assume now that $h\in[3,2p+1]'$. In the right hand side of
(e) the sum over $i$ can be restricted to those $i$ such that $i+h\n\{1,3,\do,4p-3\}$ hence such that 
$i+h\ge4p-1$; for such $i$ we have $i\ge4p-1-h\ge(4p-1)-(2p+1)$ hence $i\ge2p-2$. Moreover, if $i=2p-2$, 
then we must have $h=2p+1$. Thus we have
$$\align&(\tw_0,\tw_h)=(-1)^{(h+1)/2}\sum_{i\in[4p-1-h,4p-2]''}(-1)^{i/2}\bc_i(i+h-1)\\&\T
(i+h-3)\do(i+h-4p+3)(4p-2)!!\i\\&
=(-1)^{(h+1)/2}\sum_{i\in[4p-1-h,4p-2]'';i\ge2p}(-1)^{i/2}(n_0+n_2+\do+n_{4p-2-i})\\&
\T(i+h-1)(i+h-3)\do(i+h-4p+3)(4p-2)!!\i\\&
-(-1)^{(h+1)/2}\d_{h,2p+1}(-1)^{(2p-2)/2}(n_0+n_2+\do+n_{2p-2})\\&
=(-1)^{(h+1)/2}\sum_{i\in[4p-1-h,4p-2]''}(-1)^{i/2}(n_0+n_2+\do+n_{4p-2-i})\\&
\T(i+h-1)(i+h-3)\do(i+h-4p+3)(4p-2)!!\i\\&
-(-1)^{p+1}\d_{h,2p+1}(-1)^{p-1}(n_0+n_2+\do+n_{2p})\\&
-(-1)^{p+1}\d_{h,2p+1}(-1)^{p-1}(n_0+n_2+\do+n_{2p-2})\\&
=x-\d_{h,2p+1}(n_0+n_2+\do+n_{2p}+n_0+n_2+\do+n_{2p-2})\\&
=x-\d_{h,2p+1}(n_0+n_2+\do+n_{2p}+n_{2p+2}+\do+n_{4p})\\&
=x-\d_{h,2p+1}2^{2p}\endalign$$
where
$$\align&x=(-1)^{(h+1)/2}\sum_{i\in[4p-1-h,4p-2]''}(-1)^{i/2}(n_0+n_2+\do+n_{4p-2-i})\\&
\T(i+h-1)(i+h-3)\do(i+h-4p+3)(4p-2)!!\i.\endalign$$
It remains to show that $x=1$. Setting $h=2h'+1,i=4p-2-2i'$ we have
$$\align&x=\sum_{i'\in[0,h']}(-1)^{i'+h'}(n_0+n_2+\do+n_{2i'})\\&\T(2p-i+h'-1)(2p-i+h'-2)\do(h'-i+1)
(2p-1)!\i\\&=\sum_{i'\ge0,u\ge0;i'+u=h'}(-1)^u(n_0+n_2+...+n_{2i'})r_u\endalign$$
where $r_u=(u+1)(u+2)...(u+2p-1)/(2p-1)!$. Note that 
$$\sum_{i\ge0,u\ge0;i+u=e}(-1)^in_{2i}r_e=\d_{e,0}$$
for any $e\in\NN$. Hence
$$\align&x=\sum_{i'\ge0,j\ge0,r\ge0,u\ge0;i'=j+r,i'+u=h'}(-1)^{h'+j+r}n_{2j}r_u\\&=
\sum_{r\in[0,h']}(-1)^{h'+r}\sum_{j,u\ge0;j+u=h'-r}(-1)^jn_{2j}r_u\\&=
\sum_{r\in[0,h']}(-1)^{h'+r}\d_{h'-r}=(-1)^{h'+h'}=1.\endalign$$
This completes the proof of (a),(b),(c).

\subhead 2.10\endsubhead
We fix two integers $p_1,p_2$ such that $p_1\ge p_2\ge1$. Let $V',V''$ be two $\kk$-vector spaces of 
dimension $2p_1,2p_2-2$ respectively and let $V=V'\op V''$. We identify $V^*=V'{}^*\op V''{}^*$ in the 
obvious way. Let $(,):V\T V^*@>>>\kk$ be the obvious pairing. Assume that $V'$ has a given basis 
$\{z_i;i\in[0,4p_1-2]''\}$ and that $V''$ has a given basis $\{v_i;i\in[0,4p_2-6]''\}$. There is a unique 
basis $\{z_i;i\in[1,4p_1-1]'\}$ of $V'{}^*$ and a unique basis $\{v_i;i\in[1,4p_2-5]'\}$ of $V''{}^*$ such 
that
$$(z_i,z_j)=\ph_{p_1}(i-j)\text{ for }i\in[0,4p_1-2]'',j\in[1,4p_1-1]',$$
$$(v_i,v_j)=\ph_{p_2-1}(i-j)\text{ for }i\in[0,4p_2-6]'',j\in[1,4p_2-5]'.$$
(Notation of 2.4; the basis of $V''$ and $V''{}^*$ is empty when $p_2=1$.) We define  
$g\in G^1_V$ by $gz_i=z_{i+1}$ for $i\in[0,4p_1-2]''$, $gv_i=v_{i+1}$ for $i\in[0,4p_2-6]''$. We have 
$$g^{*2}(z_i)=z_{i+2}\text{ for }i\in[0,4p_1-4],\qua g^{*2}(v_i)=v_{i+2}\text{ for }i\in[0,4p_2-8],$$
$$(g^{*2}+1)^{2p_1}=0\text{ on }V',\qua (g^{*2}+1)^{2p_2-2}=0\text{ on }V''.$$
(See 2.5). Hence $-g^{*2}$ acts on $V'$ as a single unipotent Jordan block of size $2p_1$ and on $V''$ as a 
single unipotent Jordan block of size $2p_2-2$. (When $p_2=1$, $-g^{*2}=0$ on $V''=0$.)

For $i\in\ZZ$ we write $z_i$ instead of $(z_0)_i$ (as in 2.0); when $p_2\ge2$ we write $v_i$ instead of 
$(v_0)_i$. This agrees with our earlier notation for $z_i$ when $i\in[0,4p_1-1]$ and $v_i$ for 
$i\in[0,4p_2-5]$. We have
$$(z_i,z_j)=\ph_{p_1}(i-j)\text{ for }i\in\ZZ'',j\in\ZZ';$$
$$(v_i,v_j)=\ph_{p_2-1}(i-j)\text{ for }i\in\ZZ'',j\in\ZZ'\text{ (assuming $p_2\ge2$)}.$$
(See 2.6(a).) If $p_2\ge2$ we clearly we have
$$(z_i,v_j)=0, (v_i,z_j)=0\text{ for }i\in\ZZ'',j\in\ZZ'.$$
As in 2.7, 2.8, there is a unique vector $\tz\in V'$ such that for any $h\in\ZZ'$ we have
$$(\tz_0,z_h)=2^{p_1}(-1)^{(h+1)/2}(h-1)(h-3)\do(h-4p_1+3)(4p_1-2)!!\i,$$
Similarly if $p_2\ge2$, there is a unique vector $\tv\in V''$ such that for any $h\in\ZZ'$ we have
$$(\tv_0,v_h)=2^{p_2-1}(-1)^{(h+1)/2}(h-1)(h-3)...(h-4p_2+7)(4p_2-6)!!\i.$$
(Notation of 2.0.) If $p_2=1$ we set $\tv_i=0$ for all $i\in\ZZ$. As in 2.9, we have
$$(\tz_0,\tz_h)=\sum_{k\in[1,p_1]}2^{2k-2}\ph_k(h),\tag a$$
$$(\tv_0,\tv_h)=\sum_{k\in[1,p_2-1]}2^{2k-2}\ph_k(h),\text{ (if $p_2\ge2$)},\tag b$$
$$(\tz_0,\tz_h)=1\text{ if }h\in[-2p_1+1,2p_1-1]';\qua (\tz_0,\tz_{2p_1+1})=1-2^{2p_1},\tag c$$
$$(\tv_0,\tv_h)=1\text{ if }h\in[-2p_2+3,2p_2-3]'; \qua(\tv_0,\tv_{2p_2-1})=1-2^{2p_2-2}
\text{ (if $p_2\ge2$)}.\tag d$$
Let $\z\in\kk$ be such that $\z^2=-1$. We set
$$\x=2^{-p_2+1}\tz_{-2p_2}+2^{-p_2+1}\z\tv_0\in V.$$
Let $h\in\ZZ'$. We show:
$$\align&(\x_0,z_h)
=2^{p_1-p_2+1}(-1)^{(h+2p_2+1)/2}(h+2p_2-1)(h+2p_2-3)\do\\&\T(h+2p_2-4p_1+3)(4p_1-2)!!\i\in2\ZZ.\endalign$$
Indeed,
$$\align&(\x_0,z_h)=2^{-p_2+1}(\tz_{-2p_2},z_h)=2^{-p_2+1}(\tz_0,z_{2p_2+h})
=2^{-p_2+1}2^{p_1}(-1)^{(2p_2+h+1)/2}\\&\T(2p_2+h-1)(2p_2+h-3)\do(2p_2+h-4p_1+3)(4p_1-2)!!\i,\endalign$$
as desired. In particular we have
$$(\x_0,z_h)=0\text{ if }h\in[1-2p_2,4p_1-2p_2-3]'.$$
Let $h\in\ZZ'$. From the definitions we have $(\x_0,\x_h)=2^{-2p_2+2}((\tz_0,\tz_h)-(\tv_0,\tv_h))$. From 
this we deduce using (a)-(d) that
$$(\x_0,\x_h)=\sum_{k\in[p_2,p_1]}2^{2k-2p_2}\ph_k(h)\in\ZZ\text{ for }h\in\ZZ',$$
$$(\x_0,\x_h)=0\text{ if }h\in[-2p_2+3,2p_2-3]'; (\x_0,\x_{2p_2-1})=1.$$
It follows that, if $L$ is the line in $V$ spanned by $z_0$, $L'$ is the line in $V$ spanned by $\x_0$ and 
$a_*=(0,0,0,\do)$, $b_*=(2p_1,2p_2-2,0,\do)$, then $(g,L,L')\in\tcc^V_{a_*,b_*}$; in particular, 
$\tcc^V_{a_*,b_*}\ne\emp$.

\subhead 2.11\endsubhead
Let $p_1,p_2$ be integers such that $p_1\ge p_2\ge1$. We consider a $\kk$-vector space $V$ of dimension 
$2p_1+2p_2-2$ with a given bilinear form $g\in G^1_V$ such that (with notation of 2.0)
$-g^{*2}\in G_V$ is unipotent with a single Jordan block of size $2p_1$ (if $p_2=1$) or with two Jordan 
blocks, one of size $2p_1$ and one of size $2p_2-2$ (if $p_2\ge2$). We assume given two vectors $z,\x$ in 
$V$ such that (with notation of 2.0), setting for $h\in\ZZ'$:
$$\a_h=(z_i,z_j), \b_h=(\x_i,\x_j),\g_h=(\x_i,z_j)\text{ where }i\in\ZZ'',j\in\ZZ', h=j-i,$$
we have
$$\a_h=0\text{ if }h\in[-2p_1+3,2p_1-3]',\a_{2p_1-1}=1,$$
$$\b_h=0\text{ if }h\in[-2p_2+3,2p_2-3]',\b_{2p_2-1}=1,$$
$$\g_h=0\text{ if }h\in[1-2p_2,4p_1-2p_2-3]'.$$
We show:

(a) {\it After possibly replacing $\x$ by $-\x$, the following equalities hold for any $h\in\ZZ'$:

(a1) $\a_h=\ph_{p_1}(h)\in\ZZ$,

(a2) $\b_h=\sum_{k\in[p_2,p_1]}2^{2k-2p_2}\ph_k(h)\in\ZZ$ for $h\in\ZZ'$,

(a3) $\g_h=2^{p_1-p_2+1}(-1)^{(h+2p_2+1)/2}(h+2p_2-1)(h+2p_2-3)\do(h+2p_2-4p_1+3)(4p_1-2)!!\i\in2\ZZ$.}
\nl
($\ph_p$ as in 2.4.) We prove (a1). If $|h|\le2p_1-1$, then (a1) is clear. Thus we can assume that 
$|h|\ge2p_1+1$. Since $\a_h=\a_{-h}$ we can also assume that $h\ge1$ (hence $h\ge2p_1+1$). We must only 
prove that

(b) $\a_h=x_{h-2p_1+1}$ if $h\ge2p_1-1$ is odd,
\nl
where $x_e$ is as in 2.4 (with $p=p_1$). We have $(g^{*2}+1)^{2p_1}=0$ on $V$ hence applying to $z_0$, we 
have$\sum_{j\in[0,2p_1]}r_jz_{2j}=0$ where $r_j=\bin{2p_1}{j}$. Taking $(,z_{2p_1+2s-1})$ we get 
$\sum_{j\ge0}r_j\a_{2p_1+2s-1-2j}=0$. The coefficient of $T^s$ $(s\in\NN$) in 
$$(\sum_{j\in\NN}r_jT^j)(\sum_{u\in\NN}\a_{2p_1-1+2u}T^u)$$ 
is $$k_s=\sum_{j\in[0,s]}r_j\a_{2p_1-1+2s-2j}.$$
If $s\ge2,j>s,j\le2p_1$, we have $\a_{2p_1-1+2s-2j}=0$ since $2p_1-3\ge2p_1-1+2s-2j\ge-2p_1+3$; hence 
$k_s=\sum_{j\ge0}r_j\a_{2p_1-1+2s-2j}$ for $s\ge2$. We have 
$$r_0\a_{2p_1+1}+r_1\a_{2p_1-1}+r_{2p_1}\a_{-2p_1+1}=0$$ hence $\a_{2p_1+1}=-(2p_1+1)$ and 
$$k_1=r_0\a_{2p_1+1}+r_1\a_{2p_1-1}=-1.$$ Also $k_0=1$. Thus $\sum_{s\ge0}k_sT^s=1-T$. The left hand side is
$$(\sum_{j\ge0}r_jT^j)(\sum_{u\ge0}\a_{2p_1-1+2u}T^u).$$
Thus $\sum_{u\ge0}\a_{2p_1-1+2u}T^u=(1-T)(1+T)^{-2p_1}$. On the other hand from the definition of $x_{2u}$ 
we have $\sum_{u\ge0}x_{2u}T^u=(1-T)(1+T)^{-2p_1}$. This proves (b) hence (a1).

Note that

(c) $\{z_i;i\in[0,4p_1-4]''\})$ together with $\{\x_i;i\in[0,4p_2-4]''\}$ form a basis of $V$.

\subhead 2.12\endsubhead
We show:

(a) $\{z_{2i};i\in[0,2p_1-1]\}$ are linearly independent.
\nl
Assume that this is not true. Then $z_{4p_1-2}\in E$, the span of $\{z_i;i\in[0,4p_1-4]''\})$ hence $E$ is 
$g^{*2}$-stable and the annihilator $(gE)^\pe$ of $gE$ in $V$ is $g^{*2}$-stable. For $i\in[0,2p_1-2]$ we 
have $(\x_{2p_2},z_{2i+1})=0$ hence $\x_{2p_2}\in(gE)^\pe$. Since $(gE)^\pe$ is $g^{*2}$-stable we see that
$\x_i\in(gE)^\pe$ for all $i\in\ZZ''$. Thus $E'$, the span of $\{\x_i;i\in[0,4p_2-4]''\}$, is contained in 
$(gE)^\pe$. Now $E'$ has dimension $2p_2-1$ which is the same as $\dim(gE)^\pe$. Hence $E'=(gE)^\pe$. Since 
$V=E\op E'$ (see 2.11(c)) we see that $V=E\op(gE)^\pe$ with both summands $g^{*2}$-stable. Now $-g^{*2}$ 
acts on $E$ as a single Jordan block of size $2p_1-1$. Thus $-g^{*2}:V@>>>V$ has a Jordan block of size 
$2p_1-1$. This contradicts the assumption that the Jordan blocks of $-g^{*2}:V@>>>V$ have even sizes. This 
proves (a).

We set $N=g^{*2}+1,e=p_1-p_2$. Let $\cl$ be the span of $\{N^iz_0;i\in[2p_2-1,2p_1-1]\}$ or equivalently the
span of $\{N^{2p_2-1}z_i;i\in[0,4e]''\}$. We show that

(b) $\dim\cl=2e+1$.
\nl
Let $\cl'$ be the span of $\{N^iz_0;i\in[2p_2-1,2p_1-2]\}$. We have $\dim\cl'=2e$ since
$\{N^iz_0;i\in[0,2p_1-2]\}$ is a linearly independent set. If (b) is false we would have 
$N^{2p_1-1}z_0\in\cl'$. Then the span of $\{N^iz_0;i\in[0,2p_1-2]\}$ is $N$-stable. Hence the span of 
$\{g^{*(2i)}z_0;i\in[0,2p_1-2]\}$ is $g^{*2}$-stable. This contradicts the proof of (a).

We show:

(c) $N^{2p_2-1}\x_0\in\cl$.
\nl
From the structure of Jordan blocks of $N:V@>>>V$ we see that $\dim N^{2p_2-1}V=2e+1$. Clearly,
$\cl\sub N^{2p_2-1}V$. Hence using (b) it follows that $\cl=N^{2p_2-1}V$ so that (c) holds.

Using (c) we deduce 
$$N^{2p_2-1}\x_0=\sum_{i\in[0,2e]}c_{2i}N^{2p_2-1}z_{2i}\tag d$$
where $c_{2i}\in\kk$ ($i\in[0,2e]$) are uniquely determined. 

\subhead 2.13\endsubhead
For $j\in\NN$ we set $m_j=\bin{2p_2-1}{j}$ so that $N^{2p_2-1}=\sum_{j\in[0,2p_2-1]}m_jg^{*(2j)}$. From 
2.13(d) we deduce
$$\sum_{j\in[0,2p_2-1]}m_j\x_{2j}=\sum_{i\in[0,2e],j\in[0,2p_2-1]}c_{2i}m_jz_{2i+2j}.\tag a$$
Taking $(,z_u)$ with $u\in\ZZ'$ we deduce
$$\sum_{j\in[0,2p_2-1]}m_j\g_{u-2j}=\sum_{i\in[0,2e],j\in[0,2p_2-1]}c_{2i}m_j\a_{u-2i-2j}.\tag b$$
We show:

(c1) If $u\in[2p_2-1,4p_1-2p_2-3]'$ then the left hand side of (b) is $0$.

(c2) If $u=4p_1-2p_2-1$ then the left hand side of (b) is $\g_{4p_1-2p_2-1}$.
\nl
For (c1) it is enough to show: if $u$ is as in (c1) and $j\in[0,2p_2-1]$ then $u-2j+2p_2\in[1,4p_1-3]$. 
Indeed we have 
$$u-2j+2p_2\le4p_1-2p_2-3+2p_2=4p_1-3$$ and 
$$u-2j+2p_2\ge2p_2-1-4p_2+2+2p_2=1.$$
 For (c2) it is 
enough to show: if $j\in[1,2p_2-1]$ then $4p_1-2p_2-1-2j+2p_2\in[1,4p_1-3]$ or that 
$4p_1-1-2j\in[1,4p_1-3]$. This is clear.

If $u\in[2p_2-1,2p_1-3]'$ then in the right hand side of (b) we have $u-2i-2j<2p_1-1$; we can assume then 
that $u-2i-2j\le-2p_1+1$ hence 
$$2i\ge u-2j+2p_1-1\ge2p_2-1-(4p_2-2)+2p_1-1=2e$$
and $i\ge e$. Thus in this case (b) becomes (using (c1) and setting $u=2p_1-1-2t$):
$$\sum_{i\in[e,2e],j\in[0,2p_2-1]}c_{2i}m_j\a_{2p_1-1-2t-2i-2j}$$
for $t\in[1,e]$.
Setting $c'_h=c_{4e-h}$ for $h\in[0,2e]''$ and with the change of variable $j\m2p_2-1-j$, $i\m2e-i$ we obtain
$$\sum_{i\in[0,e],j\in[0,2p_2-1]}c'_{2i}m_j\a_{-2p_1+1-2t+2i+2j}=0\text{ for }t\in[1,e].\tag d$$
In the last sum we have $-2p_1+1-t+2i+2j<2p_1-1$. Indeed, we have 
$$-2p_1+1-2t+2i+2j\le-2p_1-1+2p_1-2p_2+4p_2-2=2p_2-3<2p_1-1.$$
 Hence we can restrict the sum to indices such 
that $-2p_1+1-2t+2i+2j\le-2p_1+1$ that is $-t+i+j=-2s$ where $s\ge0$. Thus we have 
$$\sum_{i\in[0,e],j\ge0,s\ge0,i+j+s=t}c'_{2i}m_j\a_{-2p_1+1-2s}=0\text{ for }t\in[1,e].$$
Hence 
$$(\sum_{i\in[0,e]}c'_{2i}T^i)(\sum_{j\ge0}m_jT^j)(\sum_{s\ge0}\a_{-2p_1+1-s}T^s)
=c'_0+\text{ terms of degree $>e$ in }T.$$
Using results in 2.11 this can be written as
$$(\sum_{i\in[0,e]}c'_{2i}T^i)(1+T)^{2p_2-1}(1-T)(1+T)^{-2p_1}=c'_0+\text{ terms of degree $>e$ in }T$$
that is,
$$(\sum_{i\in[0,e]}c'_{2i}T^i)(1+T)^{-2e-1}(1-T)=c'_0+\text{ terms of degree $>e$ in }T,$$
hence
$$\sum_{i\in[0,e]}c'_{2i}T^i=(1-T)\i(1+T)^{2e+1}(c'_0+\text{ terms of degree $>e$ in }T).$$
We have $(1+T)^{2e+1}=\sum_{j\in[0,2e+1]}l_jT^j$ where $l_j=\bin{2e+1}{j}$. Hence
$$(1-T)\i(1+T)^{2e+1}=\sum_{j\in[0,e]}(l_0+l_1+\do+l_j)T^j+\text{ terms of degree $>e$ in }T.$$
We see that 

(e) $c'_{2i}=c'_0(l_0+l_1+\do+l_i)$ for $i\in[0,e]$.
\nl
In the remainder of this subsection we assume that $e>0$. If $u=2p_1-1$ then in the right hand side of (b) 
we have $u-2i-2j\in[-2p_1+1,2p_1-1]$; we can then assume that $u-2i-2j$ is $-2p_1+1$ or $2p_1-1$. Hence 
$i+j$ is $2p_1-1$ or $0$ and $(i,j)$ is $(2e,2p_2-1)$ or $(0,0)$. Thus in this case (b) becomes (using (c1))
$c_0+c_{4e}=0$ that is $c_0=-c'_0$. (The left hand side of (b) is $0$ by (c1); here we use that $e>0$.)

If $u\in[2p_1+1,4p_1-2p_2-3]'$ then in the right hand side of (b) we have $u-2i-2j>-2p_1+1$; we can then
assume that $u-2i-2j\ge2p_1-1$ hence 
$$2i\le u-2j-2p_1+1\le 4p_1-2p_2-3-2p_1+1=2e-2$$ 
and $i\le e-1$. Using this and (c1) we see that (b) becomes (setting $u=2p_1-1+2t$):
$$\sum_{i\in[0,e-1],j\in[0,2p_2-1]}c_{2i}m_j\a_{2p_1-1+2t-2i-2j}=0\text{ for }t\in[1,e-1].$$
Note that in the sum we have $2p_1-1+2t-2i-2j>-2p_1+1$. (Indeed we have 
$2p_1-1+2t-2i-2j\ge2p_1+1-2p_1+2p_2+2-4p_2+2=-2p_2+5>-2p_1+1$.) Hence we can restrict the sum to indices 
such that $2p_1-1+2t-2i-2j\ge2p_1-1$ that is $2p_1-1+2t-2i-2j=2p_1-1+2s$ where $s\ge0$. Thus we have 
$$\sum_{i\in[0,e-1],j\ge0,s\ge0;i+s+j=t}c_{2i}m_j\a_{2p_1-1+2s}=0\text{ for }t\in[1,e-1].$$
For such $t$ we have also
$$\sum_{i\in[0,e-1],j\ge0,s\ge0;i+s+j=t}c'_{2i}m_j\a_{-2p_1+1-2s}=0$$
as we have seen earlier; the index $i$ cannot take the value $e$ since $i\le t$. Adding the last two 
equations and using $\a_{2p_1-1+2s}=\a_{-2p_1+1-2s}$ we obtain
$$\sum_{i\in[0,e-1],j\ge0,s\ge0;i+s+j=t}(c_{2i}+c'_{2i})m_j\a_{-2p_1+1-2s}=0\text{ for }t\in[1,e-1].
\tag $*$ $$
We show that $c_{2i}+c'_{2i}=0$ for $i\in[0,e-1]$. For $i=0$ this is already known; the general case 
follows from ($*$) by induction on $i$. Using also (e), we see that 
$$c_{2i}=-c'_0(l_0+l_1+\do+l_i)\text{ for }i\in[0,e-1].\tag f$$
(In the case where $i=0$, this is just $c_0=-c'_0$ which is already known.)

\subhead 2.14\endsubhead
If $u=4p_1-2p_2-1$, then using 2.13(b) and 2.13(c2) we have
$$\g_{4p_1-2p_2-1}=\sum_{i\in[0,2e],j\in[0,2p_2-1]}c_{2i}m_j\a_{4p_1-2p_2-1-2i-2j}.\tag a$$
Taking $(,\x_{2p_2-1})$ with 2.13(a) we obtain
$$\sum_{j\in[0,2p_2-1]}m_j\b_{2p_2-1-2j}=\sum_{i\in[0,2e],j\in[0,2p_2-1]}c_{2i}m_j\g_{2i+2j-2p_2+1}.$$
In the left hand side only the contribution of $j=0$ and $j=2p_2-1$ is $\ne0$; it is $1$; in the right hand 
side we have $2i+2j-2p_2+1\ge-2p_2+1$ hence we can assume that $2i+2j-2p_2+1>4p_1-2p_2-3$, that is 
$2i+2j\ge4p_1-2$; hence we have $i=2e,j=2p_2-1$ and the right hand side is $c_{4e}\g_{4p_1-2p_2-1}$. Thus
$$2=c'_0\g_{4p_1-2p_2-1}.\tag b$$
We see that $c'_0\ne0$ and using (a),(b) we have
$$2c'_0{}\i=\sum_{i\in[0,2e],j\in[0,2p_2-1]}c_{2i}m_j\a_{4p_1-2p_2-1-2i-2j}.$$
In the right hand side we have $4p_1-2p_2-1-2i-2j\ge-2p_1+1$; we can assume then that either 
$4p_1-2p_2-1-2i-2j=-2p_1+1$ (hence $i=2e,j=2p_2-1$) or $4p_1-2p_2-1-2i-2j\ge2p_1-1$ (hence $i\le e$). The 
first case can arise only if $e=0$ hence it is included in the second case. Thus
$$2c'_0{}\i=\sum_{i\in[0,e],j\in[0,2p_2-1]}c_{2i}m_j\a_{4p_1-2p_2-1-2i-2j}.\tag c$$
Assume now that $e>0$. From 2.13(d) with $t=e$ we have
$$0=\sum_{i\in[0,e],j\in[0,2p_2-1]}c'_{2i}m_j\a_{-4p_1+2p_2+1+2i+2j}.\tag d$$
We now add (c) and (d) and use that $c_{2i}+c'_{2i}=0$ if $i\in[0,e-1]$ and $c_e=c'_e$. We get
$$2c'_0{}\i=2c'_{2e}\sum_{j\in[0,2p_2-1]}m_j\a_{2p_1-1-2j}.$$
If $j\in[1,2p_2-1]$ we have $2p_1-1-2j\in[-2p_1+3,2p_1-3]$ hence $\a_{2p_1-1-j}=0$. Thus
$2c'_0{}\i=2c'_{2e}=2c'_02^{2e}$ and $c'_0{}^2=2^{-2e}$. Changing if necessary $\x$ by $-\x$ we can 
therefore assume that 
$$c'_0=2^{-e}.\tag e$$
Assume now that $e=0$. We have $c'_0=c_0$ and (c) becomes
$$2c_0\i=\sum_{j\in[0,2p_2-1]}c_0m_j\a_{2p_1-1-2j}$$
that is, $2c_0\i=2c_0$ hence $c_0^2=1$. Changing if necessary $\x$ by $-\x$ we can therefore assume that 
$c_0=1$. Thus (e) holds without the assumption $e>0$.     

Using (e) we rewrite 2.13(e), 2.13(f) as follows:
$$c_{2e-i}=2^{-e}(l_0+l_1+\do+l_i)\text{ for }i\in[0,e],\tag f$$ 
$$c_i=-2^{-e}(l_0+l_1+\do+l_i)\text{ for }i\in[0,e-1].\tag g$$
When $z_i,\x_i$ are replaced by the vectors with the same name in 2.10, the quantities $c_{2i}$ become the 
quantities $c_{2i}^0$. (Here $i\in[0,2e]$.) We show that
$$c_{2i}=c_{2i}^0\text{ for }i\in[0,2e]. \tag h$$
By the analogue of (b) we have $2=c_{4e}^0\g^0_{4p_1-2p_2-1}$. By results in 2.10 we have 
$\g^0_{4p_1-2p_2-1}=2^{e+1}$. Hence $c_{4e}^0=2^{-e}$. Using this and the analogues of 2.13(e), 2.13(f) we 
see that $c^0_{2i}$ are given by the same formulas as $c_{2i}$ in (e),(f). This proves (h).

\subhead 2.15\endsubhead
Let $C=\sum_{t\ge0}\g_{4p_1-2p_2-1+2t}T^t$, $C^0=\sum_{t\ge0}\g^0_{4p_1-2p_2-1+2t}T^t$. If 
$u=4p_1-2p_2-1+2t,t\ge0$ then for any $j$ that contributes to the left hand side of 2.13(b) we have 
$u-2j\ge-2p_2+1$. Indeed,
$$u-2j\ge4p_1-2p_2-1-2j\ge4p_1-2p_2-1-4p_2+2\ge-2p_2+1$$
hence we can assume that 
in the left hand side of 2.13(b) we have $u-2j\ge4p_1-2p_2-1$. Muliplying both sides of 2.13(b) with $T^t$
and summing over all $t\ge0$ we thus obtain
$$\align&\sum_{t\ge0}\sum_{j\in[0,2p_2-1];t-j\ge0}m_j\g_{4p_1-2p_2-1+2t-2j}T^t\\&=
\sum_{t\ge0}\sum_{i\in[0,2e],j\in[0,2p_2-1]}c_{2i}m_j\a_{4p_1-2p_2-1+2t-2i-2j}T^t.\endalign$$
The left hand side equals
$$(\sum_{j\in[0,2p_2-1]}m_jT^j)(\sum_{t'\ge0}\g_{4p_1-2p_2-1+2t'}T^{t'})=(1+T)^{2p_2-1}C.$$
Thus, 
$$C=(1+T)^{-2p_2+1}\sum_{t\ge0}\sum_{i\in[0,2e],j\in[0,2p_2-1]}c_{2i}m_j\a_{4p_1-2p_2-1+2t-2i-2j}T^t.$$
Similarly we have
$$C^0=(1+T)^{-2p_2+1}\sum_{t\ge0}\sum_{i\in[0,2e],j\in[0,2p_2-1]}c^0_{2i}m_j\a^0_{4p_1-2p_2-1+2t-2i-2j}T^t.$$
By 2.14(h) we have $c_{2i}=c^0_{2i}$. By 2.11(a1) we have
$$\a_{4p_1-2p_2-1+2t-2i-2j}=\a^0_{4p_1-2p_2-1+2t-2i-2j}$$
 for all $i,j,t$. It follows that $C=C^0$ hence
$$\g_{4p_1-2p_2-1+2t}=\g^0_{4p_1-2p_2-1+2t}\text{ for any }t\ge0.\tag a$$
We set $C'=\sum_{t\ge0}\g_{2p_2-3-2t}T^t$, $C'{}^0=\sum_{t\ge0}\g^0_{2p_2-3-2t}T^t$. If $u=2p_2-3-2t,t\ge0$,
then for any $j$ that contributes to the left hand side of 2.13(b) we have $u-2j\le4p_1-2p_2-3$ (indeed,
$u-2j\le2p_2-3-2j\le2p_2-3\le4p_1-2p_2-3$) hence we can assume that in the left hand side of 2.13(b) we have
$u-2j\le-2p_2-1$. With the substitution $j\m2p_2-1-j$ the previous inequality becomes $j-t\le0$ and the left
hand side of 2.13(b) becomes
$$\sum_{j\in[0,2p_2-1]}m_j\g_{u-4p_2+2+2j}=\sum_{j\in[0,2p_2-1]}m_j\g_{-2p_2-1+2(j-t)}.$$
Muliplying both sides of 2.13(b) with $T^t$ and summing over all $t\ge0$ we thus obtain
$$\sum_{t\ge0,j\ge0;t-j\ge0}m_j\g_{-2p_2-1+2(j-t)}T^t
=\sum_{t\ge0}\sum_{i\in[0,2e],j\in[0,2p_2-1]}c_{2i}m_j\a_{2p_2-3-2t-2i-2j}T^t.$$
The left hand side equals
$$(\sum_{j\in[0,2p_2-1]}m_jT^j)(\sum_{t'\ge0}\g_{-2p_2-1-2t'}T^{t'})=(1+T)^{2p_2-1}C'.$$
Thus,
$$C'=(1+T)^{-2p_2+1}\sum_{t\ge0}\sum_{i\in[0,2e],j\in[0,2p_2-1]}c_{2i}m_j\a_{2p_2-3-2t-2i-2j}T^t.$$
Similarly we have
$$C'{}^0=(1+T)^{-2p_2+1}\sum_{t\ge0}\sum_{i\in[0,2e],j\in[0,2p_2-1]}c^0_{2i}m_j\a^0_{2p_2-3-2t-2i-2j}T^t.$$
By 2.14(h) we have $c_{2i}=c^0_{2i}$. By 2.11(a1) we have $$\a_{2p_2-3-2t-2i-2j}=\a^0_{2p_2-3-2t-2i-2j}$$
for all $i,j,t$. It follows that $C'=C'{}^0$ hence
$$\g_{2p_2-3-2t}=\g^0_{2p_2-3-2t}\text{ for any }t\ge0.\tag b$$
Clearly, (a),(b) imply 2.11(a3).

\subhead 2.16\endsubhead
We set $B=\sum_{s\ge0}\b_{2p_2-1+2s}T^s$, $B^0=\sum_{s\ge0}\b^0_{2p_2-1+2s}T^s$. Let $t\ge1$. Taking 
$(,\x_{2p_2-1+2t})$ with 2.13(a) we obtain
$$\sum_{j\in[0,2p_2-1]}m_j\b_{2p_2-1+2t-2j}=\sum_{i\in[0,2e],j\in[0,2p_2-1]}c_{2i}m_j
\g_{2i+2j-2p_2+1-2t}.\tag a$$
For any $j$ that contributes to the left hand side of (a) we have $2p_2-1+2t-2j\ge-2p_2+3$ (indeed, 
$2p_2-1+2t-2j\ge2p_2+1-4p_2+2=-2p_2+3$) hence we can assume that in the left hand side of (a) we have 
$2p_2-1+2t-2j\ge2p_2-1$ that is $t\ge j$. Multiplying both sides of (a) by $T^t$ and summing over all 
$t\ge1$, we thus obtain
$$\align&\sum_{t\ge1}\sum_{j\in[0,2p_2-1];t\ge j}m_j\b_{2p_2-1+2t-2j}T^t\\&
=\sum_{t\ge1}\sum_{i\in[0,2e],j\in[0,2p_2-1]}c_{2i}m_j\g_{2i+2j-2p_2+1-2t}T^t.\endalign$$
The left hand side equals
$$-1+(\sum_{j\in[0,2p_2-1]}m_jT^j)(\sum_{t'\ge0}\b_{2p_2-1+t'}T^{t'})=1+(T+1)^{2p_2-1}B.$$
Thus,
$$B=(T+1)^{2p_2-1}(1+\sum_{t\ge1}\sum_{i\in[0,2e],j\in[0,2p_2-1]}c_{2i}m_j\g_{2i+2j-2p_2+1-2t}T^t).$$
Similarly we have
$$B^0=(T+1)^{2p_2-1}(1+\sum_{t\ge1}\sum_{i\in[0,2e],j\in[0,2p_2-1]}c^0_{2i}m_j\g^0_{2i+2j-2p_2+1-2t}T^t).$$
By 2.14(h) we have $c_{2i}=c_{2i}^0$. By 2.11(a3) we have 
$$\g_{2i+2j-2p_2+1-2t}=\g^0_{2i+2j-2p_2+1-2t}$$
for any $i,j,t$. It follows that $B=B^0$. Hence 
$$\b_{2p_2-1+2s}=\b^0_{2p_2-1+2s}$$
for any $s\ge0$. This clearly implies 2.11(a2).

\subhead 2.17\endsubhead
We preserve the setup of 2.1. We prove 2.1(a) by induction on $n$. If $n=0$ we have $V=0$ and 
$a_i=b_i=p_i=0$ for all $i$. We take $g=0$ and $(L^t)$ to be the empty set of lines. We obtain an element of
$\tcc^V_{a_*,b_*}$. Now assume that $n>0$.

Assume first that $a_1\ge1$. We can find a direct sum decomposition $V=V'\op V''$ such that 
$\dim V'=a_1+b_1=2p_1-1$. We identify $V^*=V'{}^*\op V''{}^*$ in the obvious way. Let $a'_*$ be the sequence
$a_1,0,0,\do$; let $b'_*$ be the sequence $b_1,0,0,\do$; let $a''_*$ be the sequence $a_2,a_3,\do$; let 
$b''_*$ be the sequence $b_2,b_3,\do$. By the induction hypothesis we have 
$\tcc^{V''}_{a''_*,b''_*}\ne\emp$.
By 2.3 we have $\tcc^{V'}_{a'_*,b'_*}\ne\emp$. Let $(g',L^1)\in\tcc^{V'}_{a'_*,b'_*}$ and let 
$(g'',L^2,L^3,\do)\in\tcc^{V''}_{a''_*,b''_*}$. Here $g'\in G^1_{V'}$, $g''\in G^1_{V''}$. Let 
$g=g'\op g''\in G^1_V$. Clearly, $(g,L^1,L^2,\do)\in\tcc^V_{a_*,b_*}$ hence 2.1(a) holds in this case. Thus 
we may assume that $a_1=a_2=\do=0$ and $b_1>0$. We see that $-g^{*2}$ is unipotent. We can find a direct sum
decomposition $V=V'\op V''$ such that $\dim V'=b_1+b_2$. We identify $V^*=V'{}^*\op V''{}^*$ in the obvious 
way. Let $b'_*$ be the sequence $b_1,b_2,0,\do$; let $b''_*$ be the sequence $b_3,b_4,\do$; let $a'_*=a''_*$
be the sequence $0,0,\do$. By the induction hypothesis we have $\tcc^{V''}_{a''_*,b''_*}\ne\emp$. By 2.11 
we have $\tcc^{V'}_{a'_*,b'_*}\ne\emp$. Let $(g',L^1,L^2)\in\tcc^{V'}_{a'_*,b'_*}$ and let 
$(g'',L^3,L^4,\do)\in\tcc^{V''}_{a''_*,b''_*}$. Here $g'\in G^1_{V'}$, $g''\in G^1_{V''}$. 
Clearly, $(g'\op g'',L^1,L^2,\do)\in\tcc^V_{a_*,b_*}$ hence 2.1(a) holds in this 
case. This completes the proof of 2.1(a).

In the following result we preserve the setup of 2.1. 

\proclaim{Proposition 2.18} Let $(g,L^1,L^2,\do,L^\s)\in\tcc^V_{a_*,b_*}$. Let $\ph_r$ be as in 2.4. There 
exist vectors $z^t\in L^t-\{0\}$ for $t\in[1,\s]$ such that (i),(ii) below hold for $i\in\ZZ'',j\in\ZZ'$.

(i) Assume that $t\in[1,\s],a_t>0$. Then $(z^t_i,z^t_j)=x'_{i-j}$ ($x'_h$ as in 2.2 with $p=p_t$); 
$(z^t_i,z^{t'}_j)=0$ if $t'\in[1,\s],t'\ne t$.

(ii) Assume that $\{t,t+1\}\sub[k+1,\s],t=k+1\mod2$ and $a_t=0$. Then 

$(z^t_i,z^t_j)=\ph_{p_t}(i-j)$,

$(z^{t+1}_i,z^{t+1}_j)=\sum_{r\in[p_{t+1},p_t]}2^{2r-2p_{t+1}}\ph_r(i-j)$,
$$\align&(z^t_i,z^{t+1}_j)
=2^{p_t-p_{t+1}+1}(-1)^{(i-j+2p_2+1)/2}(i-j+2p_{t+1}-1)(i-j+2p_{t+1}-3)\do\\&
\T(i-j+2p_{t+1}-4p_t+3)(4p_t-2)!!\i,\endalign$$
$(z^t_i,z^{t'}_j)=0$ if $t'\in[1,\s],t'\n\{t,t+1\}$.
\endproclaim
We argue by induction on $n$. When $n=0$ the result is obvious. Now assume that $n\ge1$.

{\it Case 1.} Assume first that $a_1\ge1$. We have $a_1+b_1=2p_1-1$. Let 
$V'=\op_{i\in[0,4p_1-4]''}L^1_i\sub V$. We show that 
$$g^{*2}V'=V'.\tag a$$
It is enough to show that $g^{*2}L^1_{4p_1-4}\sub V'$. Since $g^{*i}L_0^1\in V'$ for $i\in[0,4p_1-4]''$ and 
$a_1+b_1=2p_1-1$ it is enough to show that $(g^{*2}-1)^{a_1}(g^{*2}+1)^{b_1}L^1_0=0$. It is also enough to 
show that $(g^{*2}-1)^{a_1}(g^{*2}+1)^{b_1}=0$ on $V$. But this follows from the fact that 
$g\in\cc^V_{a_*,b_*}$.

Now let $V''=\op_{t\in[2,\s],i\in[0,2p_t-2]}L^t_i\sub V$. We show that

(b) {\it $V''=(gV')^\pe$, the annihilator of $gV'$ in $V$. Hence $V''$ is $g^{*2}$-stable and
$V=V'\op(gV')^\pe$.} 
\nl
We have $(L_{2p_r}^r,L^1_{i+1})=0$ for $r\in[2,\s]$, $i\in[0,4p_1-4]''$. Thus 
$L_{2p_r}^r\sub(gV')^\pe$. Since $(gV')^\pe$ is $g^{*2}$-stable (we use (a) and 2.0(a)) it 
follows that $L_i^r\sub(gV')^\pe$ for any $i\in\ZZ'',r\in[2,\s]$. Thus $V''\sub(gV')^\pe$. But 
these two vector spaces have the same dimension so that $V''=(gV')^\pe$ and (b) follows.

We identify $V^*=V'{}^*\op V''{}^*$ in the obvious way. From (a),(b) we see that $g\in G^1_V$ restricts to
an isomorphism $g':V'@>>>V'{}^*$ and to an isomorphism $g'':V''@>>>V''{}^*$. We show:

(c) {\it $g'{}^{*2}$ restricted to the generalized $1$-eigenspace of $g'{}^{*2}$ is unipotent with a single 
Jordan block of size $a_1$; $-g'{}^{*2}$ restricted to the generalized $(-1)$-eigenspace of $g'{}^{*2}$ is 
unipotent with a single Jordan block of size $b_1$ (if that eigenspace is $\ne0$). Moreover, $g''{}^{*2}$ 
restricted to the generalized $1$-eigenspace of $g''{}^{*2}$ is unipotent with Jordan blocks of sizes given 
by the nonzero numbers in $a_2,a_3,\do$; $-g''{}^{*2}$ restricted to the generalized $(-1)$-eigenspace of 
$g''{}^{*2}$ is unipotent with Jordan blocks of sizes given by the nonzero numbers in $b_2,b_3,\do$.}
\nl
As we have seen earlier we have $(g^{*2}-1)^{a_1}(g^{*2}+1)^{b_1}=0$ on $V'$ (even on $V$). Also 
$g'{}^{*2}\in GL(V')$ is regular in the sense of Steinberg and $\dim V'=a_1+b_1$. This implies (c).

Let $a'_*$ be the sequence $a_1,0,0,\do$; let $b'_*$ be the sequence $b_1,0,0,\do$; let $a''_*$ be the 
sequence $a_2,a_3,\do$; let $b''_*$ be the sequence $b_2,b_3,\do$. Now the proposition holds when 
$(g,L^1,L^2,\do)$ is replaced by $(g'',L^2,L^3,\do)\in\tcc^{V''}_{a''_*,b''_*}$ (by the induction
hypothesis) or by $(g',L^1)\in\tcc^{V'}_{a'_*,b'_*}$ (we choose any $z^1\in L^1-\{0\}$ such that 
$(z^1_i,z^1_j)=1$ for $i\in\ZZ'',j\in\ZZ'$, $|i-j|=2p_1-1$ and we apply 2.3). Hence the proposition holds 
for $(g,L^1,L^2,\do)$ (we use (b)).

{\it Case 2.} Next we assume that $k=0$, $b_1>0$. Then $a_1=a_2=\do=0$. We have $b_1=2p_1,b_2=2p_2-2$. Let 
$V'=\op_{t\in[1,2],i\in[0,4p_t-4]''}L^t_i\sub V$. We show that 
$$g^{*2}V'=V'.\tag d$$
Let $N=g^{*2}+1$. Then $V=\op_{t\in[1,\s],i\in[0,4p_t-4]''}N^{i/2}L^t_0$ is a direct sum decomposition into
lines. Now  $N^{2p_2-2}(V)$ contains the lines

($*$) $N^{2p_2-2+(i/2)}L^1_0 (i\in[0,4p_1-4p_2]'')$ and $N^{2p_2-2}L^2_0$
\nl
(whose number is $2p_1-2p_2+2$); moreover, since $N$ has Jordan blocks of sizes $b_1=2p_1,b_2=2p_2-2$ 
and others of size $<b_2$ we see that $\dim N^{2p_2-2}(V)=2p_1-2p_2+2$ so that $N^{2p_2-2}(V)$ is equal to 
the subspace spanned by ($*$) and $N^{2p_2-2}(V)\sub V'$. Now $V'$ is the subspace of $V$ spanned by the 
lines $N^iL^t_0$ with $t\in[1,2],i\in[0,2p_t-2]$. It is enough to show that $NV'\sub V'$ or that 
$N^{2p_t-1}L^t_0\sub V'$ for $t=1,2$. But for $t=1,2$ we have $N^{2p_t-1}L^t_0\sub N^{2p_2-2}V\sub V'$ since
$2p_t-2p_2+1\ge0$. This proves (d).

\mpb

Let $V''=\op_{t\in[3,\s],i\in[0,4p_t-4]''}L^t_i\sub V$. We show that 

(e) {\it $V''=(gV')^\pe$, the annihilator of $gV'$ in $V$. Hence $V''$ is $g^{*2}$-stable and
$V=V'\op(gV')^\pe$.}
\nl
We have $(L_{2p_r}^r,L^t_{i+1})=0$ for $t\in[1,2],r\in[3,\s]$, $i\in[0,4p_t-4]''$. Thus
$L_{2p_r}^r\sub(gV')^\pe$. Since $(gV')^\pe$ is $g^{*2}$-stable (we use (d) and 2.0(a)) it 
follows that $L_i^r\sub(gV')^\pe$ for any $i\in\ZZ'',r\in[3,\s]$. Thus $V''\sub(gV')^\pe$. But 
these two vector spaces have the same dimension so that $V''=(gV')^\pe$ and (e) follows.

\mpb

We identify $V^*=V'{}^*\op V''{}^*$ in the obvious way. From (d),(e) we see that $g:V@>>>V^*$ restricts to
an isomorphism $g':V'@>>>V'{}^*$ and to an isomorphism $g'':V''@>>>V''{}^*$. We show:

(f) {\it $-g'{}^{*2}$ is unipotent with a single Jordan block of size $b_1$ (if $b_2=0$) or with
two Jordan blocks of size $b_1,b_2$ (if $b_2>0$). Moreover, $-g''{}^{*2}$ is unipotent with Jordan blocks of 
sizes given by the nonzero numbers in $b_3,b_4,\do$.}
\nl
Since $V'$ is the direct sum of the lines $N^iL^t_0$, $t\in[1,2],i\in[0,2p_t-2]$, and $V'$ is $N$-stable, we
see that the kernel of $N:V'@>>>V'$ has dimension $\le2$. Hence $N:V'@>>>V'$ has either a single Jordan 
block of size $2p_1+2p_2-2=b_1+b_2$ or two Jordan blocks of sizes $b'_1\ge b'_2$ where $b'_1+b'_2=b_1+b_2$. 
In the first case we must have $b_2=0$ (since the Jordan blocks of $N:V'@>>>V'$ have sizes $\le b_1$ (by 
(e)). In the second case, since $b'_1,b'_2$ must form a subsequence of $b_1>b_2>b_3>\do$ and 
$b'_1+b'_2=b_1+b_2$ it follows that $b'_1=b_1$, $b'_2=b_2$. This implies (f). This completes the proof.

\subhead 2.19\endsubhead
In the setup of 2.1, we show that 2.1(b) holds. We must show that 

(a) any two elements $(g,L^1,L^2,\do,L^\s)$, $(g',L'{}^1,L'{}^2,\do,L'{}^{\s})$ of $\tcc^V_{a_*,b_*}$ are in
the same $G_V$-orbit.
\nl
Since $G_V$ acts transitively on $\cc^V_{a_*,b_*}$ we can assume that $g=g'$. Let $z^t\in L^t$ 
$(t\in[1,\s])$ be as in 2.18. Let $z'{}^t\in L'{}^t$ $(t\in[1,\s])$ be the analogous vectors for 
$(g,L'{}^1,L'{}^2,\do)$  instead of $(g,L^1,L^2,\do)$. By 2.18 we have 
$$(z^t_i,z^{t'}_j)=(z'{}^t_i,z'{}^{t'}_j)\tag b$$ 
for any $i\in\ZZ'',j\in\ZZ'$ and any $t,t'\in[1,\s]$. Since $\{z^t_i;t\in[1,\s],i\in[0,4p_t-4]''\}$ and 
$\{z'{}^t_i;t\in[1,\s],i\in[0,4p_t-4]''\}$ are bases of $V$ (see 2.0(d)) we see that there is a unique 
$\g\in GL(V)$ such that $\g(z^t_i)=z'{}^t_i$ for any $t\in[1,\s],i\in[0,4p_t-4]$. We show that 
$$\che\g(z^t_{j+1})=z'{}^t_{j+1}\text{ for any }t\in[1,\s],j\in[0,4p_t-4]''.\tag c$$
It is enough to show that $(z'{}^{t'}_i,z'{}^t_{j+1})=(z'{}^{t'}_i,\che\g(z^t_{j+1}))$ that is,
 $(z'{}^{t'}_i,z'{}^t_{j+1})=(z^{t'}_i,z^t_{j+1})$ for any $t,t'\in[1,\s]$ and any $i,j\in[0,4p_t-4]''$. 
This follows from (b). From (c) we see that $\che\g(g(z^t_j))=g(\g(z^t_j))$ for any 
$t\in[1,\s],j\in[0,4p_t-4]''$. It follows that $\che\g g=g\g$. From the definition it is clear that 
$\g(L^t)=L'{}^t$ for $t\in[1,\s]$. Thus (a) holds (with $g'=g$). This proves 2.1(b).

\subhead 2.20\endsubhead
In the setup of 2.1, we show that 2.1(c) holds. Let $(g,L^1,L^2,\do,L^\s)\in\tcc^V_{a_*,b_*}$ and let $I$ be
the set of all $\g\in G_V$ such that $\che\g g\g\i=g$, $\g(L^t)=L^t$ for $t\in[1,\s]$. Let 
$z^t\in L^t (t\in[1,\s])$ be as in 2.18. Let $\g\in I$. If $t\in[1,\s]$ we have $\g(z^t)=\o^\g_tz^t$ where 
$\o^\g_t\in\kk-\{0\}$. Since $\g$ commutes with $g^{*2}$, it follows that $\g(z^t_i)=\o^\g_tz^t_i$ for 
$i\in\ZZ''$. For $t\in[1,\s],j\in\ZZ'$ we have 
$$\che\g(z^t_j)=\che\g(g(z^t_{j-1}))=g(\g(z^t_{j-1}))=g(\o^\g_tz^t_{j-1})=\o^\g_tz^t_j;$$
thus, $\che\g(z^t_j)=\o^\g_tz^t_j$. For any $t,t'\in[1,\s]$, $i\in\ZZ'',j\in\ZZ'$ we have 
$$(z^{t'}_i,\o^\g_tz^t_j)=(z^{t'}_i,\che\g(z^t_j))=(\g\i(z^{t'}_i),z^t_j)=(\o^\g_{t'})\i(z^{t'}_i,z^t_j).$$
Thus, $(\o^\g_t-(\o^\g_{t'})\i)(z^{t'}_i,z^t_j)=0$. Taking $t'=t,i-j=2p_t-1$ we deduce that 
$\o^\g_t-(\o^\g_t)\i=0$ hence $\o^\g_t=\pm1$. Taking $t'=t+1$ (where 
$\{t,t+1\}\sub[k+1,\s],t=k+1\mod2,a_t=0$) and using that
$$(z^{t+1}_i,z^t_j)=(z^t_{-i},z^{t+1}_{-j})=\pm2^{p_t-p_{t+1}+1}\text{ if }j-i+2p_{t+1}=-1$$
we see that $(\o^\g_t-(\o^\g_{t+1})\i)2^{p_t-p_{t+1}+1}=0$ hence $\o^\g_t-(\o^\g_{t+1})\i=0$ and 
$\o^\g_t=\o^\g_{t+1}$.
We see that $\g\m(\o^\g_t)$ is a homomorphism $\ps:I@>>>\ci$ (notation of 2.0). Assume that $\g$ is in the 
kernel of $\ps$. Then $\g$ restricts to the identity map $L^t@>>>L^t$ for $t\in[1,\s]$. Since $\g$ commutes 
with $g^{*2}$ it follows that $\g$ restricts to the identity map on each of the lines
$g^{*i}L^t$ ($t\in[1,\s]$, $i\in\ZZ''$). Since these lines generate $V$ (see 2.0) we see that $\g=1$. 
Thus, $\ps$ is injective. Now let $(\o_t)\in\ci$. We define $\g\in GL(V)$ by $\g(z^t_i)=\o_tz^t_i$ for
$t\in[1,\s],i\in[0,4p_t-4]''$. From the definitions we see that 
$$(\o_tz^t_i,\o_{t'}z^{t'}_j)=(z^t_i,z^{t'}_j)\tag a$$ 
for any $i\in\ZZ'',j\in\ZZ'$ and any $t,t'\in[1,\s]$. We show that 
$$\che\g(z^t_{i+1})=\o_tz^t_{i+1}\text{ for any }t\in[1,\s],i\in[0,4p_t-4]''.\tag b$$
It is enough to show that $(\g(z^{t'}_j),\o_tz^t_{i+1})=(z^{t'}_j,z^t_{i+1})$ for any 
$t'\in[1,\s],j\in[0,4p_{t'}-4]''$ or that $(\o_{t'}z^{t'}_j,\o_tz^t_{i+1})=(z^{t'}_j,z^t_{i+1})$ or that 
$$(\o_{t'}\o_t-1)(z^{t'}_j,z^t_{i+1})=0.$$ 
The second factor is zero unless either $t=t'$ or $t'=t+1$ (where $\{t,t+1\}\sub[k+1,\s],t=k+1\mod2,a_t=0$)
in which case the first factor is zero. This proves (b).

From (b) we see that $\che\g(g(z^t_i))=g(\g(z^t_i))$ for any $t\in[1,\s],i\in[0,4p_t-4]''$. It follows that 
$\che\g g=g\g$. From the definition it is clear that $\g(L^t)=L^t$ for $t\in[1,\s]$. Thus $\g\in I$. We see
that $\ps$ is surjective hence an isomorphism. This proves 2.1(c).

\subhead 2.21\endsubhead
We now assume that $n\ge1$.
We denote by $\ovs n\to V$ (resp. $\ovs n\to V^*$) the $n$-th exterior power of $V$ (resp. $V^*$); we have
naturally $\ovs n\to V^*=(\ovs n\to V)^*$. Any $\g\in G_V$ induces an element 
$\ovs n\to\g:\ovs n\to V@>\si>>\ovs n\to V$; any $g\in G^1_V$ induces an element 
$\ovs n\to g:\ovs n\to V@>\si>>\ovs n\to V^*$. For any $\th\in\ovs n\to V-\{0\}$ we denote by $\th^*$ the 
unique element in $\ovs n\to V^*-\{0\}$ such that $(\th,\th^*)=1$.

We show:

(a) For any $g\in G^1_V$ we have $\chg g\in SL(V)$.
\nl
Let $(e_i)$ be a basis of $V$; let $(e_i^*)$ be the dual basis of $V^*$. We have $ge_i=\sum_jx_{ij}e^*_j$,
$\chg e^*_k=\sum_hy_{kh}e_h$ where $X=(x_{ij}),Y=(y_{ij})$ are square matrices. Now 
$$\d_{ki}=(\chg e^*_k,ge_i)=(\sum_hy_{kh}e_h,\sum_jx_{ij}e^*_j)=\sum_hy_{kh}x_{ih}.$$
Thus $YX^t=I$ where $X^t$ is the transpose of $X$. We have $\chg ge_i=\sum_{j,h}x_{ij}y_{jh}e_h$. Thus the 
matrix of $\chg g$ is $XY$. We have 
$$\det(XY)=\det(X)\det(Y)=\det(X^t)\det(Y)=\det(YX^t)=1,$$
as required.

We now fix $\th\in \ovs n\to V-\{0\}$ and we set
$$\G^1=\{g\in G^1_V;\ovs n\to g\text{ takes $\th$ to $\th^*$}\}.$$ 
\nl
If $g\in\G^1$ then, using (a), we see that $\ovs n\to \chg$ takes $\th^*$ to $\th$. We see that 
$\G:=SL(V)\sqc\G^1$ is a subgroup of $G_V\sqc G^1_V$. Let $SL(V)'=\{\g\in G_V;\det(\g)=\pm1\}$.

We show:

(b) {\it Let $g,g'\in\ G^1_V$, $\g\in G_V$ be such that $\che\g g\g\i=g'$. If $g,g'\in\G^1$ then 
$\g\in SL(V)'$. Conversely, if $g\in\G^1$ and $\g\in SL(V)'$ then $g'\in\G^1$. }
\nl
Replacing $V$ by $\ovs n\to V$ we can assume that $n=1$. We have $g\th=\th^*$, $g'\th=\th^*$, $\g\th=a\th$ 
where $a\in\kk-\{0\}$. We have $\th^*=\che\g g\g\i(\th)=\che\g g a\i\th=\che\g a\i\th^*=a^{-2}\th^*$ hence 
$a^2=1$ and $a=\pm1$ proving the first assertion of (b). The second assertion is proved similarly.

\subhead 2.22\endsubhead
Assuming that $a_1>0$ we show:

(a) $\cc^V_{a_*,b_*}\cap\G^1$ is a single $SL(V)$-conjugacy class in $\G$. 
\nl
Let $g,g'\in\cc^V_{a_*,b_*}\cap\G^1$. From Theorem 2.1(b) we see that $\che\g g\g\i=g'$ for some $\g\in G_V$.
Using 2.21(b) we see that $\det(\g)=\pm1$. If $\det(\g)=1$ then $g,g'$ are in the same $SL(V)$-conjugacy 
class, as required. Assume now that $\det(\g)=-1$. We complete $g$ to an element
$(g,L^1,L^2,\do)\in\tcc^V_{a_*,b_*}$ and we write $V=V'\op V''$, $V^*=V'{}^*\op V''{}^*$ as in the
proof of 2.18 (Case 1). Let $\g_0\in GL(V)$ be such that $\g_0|_{V'}=-1$, $\g_0|_{V''}=1$. Since $\dim V'$ 
is odd we have $\det(\g_0)=-1$. We have $\che\g_0 g\g_0\i=g$ hence $\che\g\che\g_0 g\g_0\i\g\i=g'$. We have 
$\g\g_0\in SL(V)$ so that $g,g'$ are in the same $SL(V)$-conjugacy class, as required. 

\subhead 2.23\endsubhead
Assuming that $a_1=0$ (hence $b_1>0$) we show:

(a) $\cc^V_{a_*,b_*}\cap\G^1$ is a union of two $SL(V)$-conjugacy classes in $\G$.
\nl
Let $g\in\cc^V_{a_*,b_*}\cap\G^1$. Let $C(g)$ (resp. $C'(g)$) be the set of elements of the form 
$\che\g g\g\i=g'$ for some $\g\in G_V$ such that $\det(\g)=1$ (resp. $\det(\g)=-1$). It is clear that 
$C(g)$ and $C'(g)$ are $SL(V)$-conjugacy classes. As in the proof of 2.22 we see, using 2.1(b) and 2.21(b), 
that $\cc^V_{a_*,b_*}\cap\G^1=C(g)\cup C'(g)$. It remains to prove that $C(g)\cap C'(g)=\emp$. Assume that 
$C(g)\cap C'(g)\ne\emp$. It follows that there exists $\g_0\in G_V$ such that $\che\g_0 g\g_0\i=g$ and 
satisfies $\det(\g_0)=-1$. Let $g_s$ be the semisimple part of $g$. Then $\g_0$ is in the centralizer of 
$g_s$ in $G_V$ which is a symplectic group all of whose elements have necessarily determinant $1$. This 
contradicts $\det(\g_0)=-1$. 

\subhead 2.24\endsubhead
Let $\boc$ be an $SL(V)$-conjugacy class contained in $\cc^V_{a_*,b_*}\cap\G^1$. (See 2.22(a), 2.23(a).) Let
$X$ be the set of all $(g,L^1,L^2,\do,L^\s)\in\tcc^V_{a_*,b_*}$ where $g\in\boc$. Note that $X\ne\emp$. Now 
$SL(V)'$ acts on $X$ by the restriction of the $G_V$-action on $\tcc^V_{a_*,b_*}$ (see 2.21(b)). Using 
2.1(b) and 2.21(b), we see that this $SL(V)'$-action is transitive. We now restrict this action to $SL(V)$.

We show:

(a) {\it This $SL(V)$-action is transitive.}
\nl
Let $(g,L^1,L^2,\do,L^\s)\in X$, $(g',L'{}^1,L'{}^2,\do,L'{}^\s)\in X$. We must show that these two
sequences are in the same $SL(V)$-orbit. As we have seen, we can find $\g\in SL(V)'$ which conjugates 
$(g,L^1,L^2,\do,L^\s)$ to $(g',L'{}^1,L'{}^2,\do,L'{}^\s)$. If $a_1=0$ this implies by the argument in 2.3 
that $\det(\g)=1$ so that in this case (a) holds. We can thus assume that $a_1>0$. If $\det(\g)=1$ then the 
proof is finished. We now assume that $\det(\g)=-1$. Let $\g_0\in G_V$ be as in 2.22. We have 
$\det(\g_0)=-1$ and $\g_0$ conjugates $(g,L^1,L^2,\do,L^\s)$ to itself. Hence $\g\g_0$ conjugates 
$(g,L^1,L^2,\do,L^\s)$ to $(g',L'{}^1,L'{}^2,\do,L'{}^\s)$. We have $\g\g_0\in SL(V)$. This proves (a).

\subhead 2.25\endsubhead
Assume that $n\ge3$. As in \cite{\WEIII, \S4} we see that 2.24(a) implies that Theorem 0.3 holds for $\G$ 
instead of $G$. 

\head 3. Exceptional groups\endhead
\subhead 3.1\endsubhead
In this section we assume that $G=G^0$ (as in 0.2) is simple of exceptional type. 
In the case where $\boc$ is a distinguished unipotent class this follows from \cite{\WEH} where it was 
proved by a reduction to a computer calculation. In the non-unipotent case the same method works
but it uses instead of \cite{\WEU, 1.2(c)}, the more general formula \cite{\WEC, 2.3(a)}.
The needed computer calculation was actually done at the time of preparing \cite{\WEC}. (I thank Frank 
L\"ubeck for providing to me tables of Green functions for groups of rank $\le8$ in GAP 
format. I also thank Gongqin Li for her help with programming in GAP to perform the computer calculation.)

We will describe below the result in the form of a list of 
rows in each case; each row corresponds to an $\e_D$-elliptic $\e_D$-conjugacy class in $W$. For 
example, the row
$$12;\Ph_{20}; (E_8(a_2))_{E_8},(E_7(a_2)A_1)_{E_7A_1},(J_{11}J_5)_{D_8}$$
in type $E_8$ corresponds to the elliptic conjugacy class $C$ in $W$ such that the characteristic polynomial
in the reflection representation of any $w\in C$ is the cyclotomic polynomial $\Ph_{20}$ and the length of
any element in $C_{min}$ is $d_C=12$. The row also includes the names of the three distinguished conjugacy 
classes $\boc$ such that $C\clu\boc$; 
for example, $(E_7(a_2)J_2)_{E_7A_1}$ is the conjugacy class of $su=us$
where $s$ is a semisimple element with $Z_G(s)^0$ of type $E_7A_1$ (in the subscript) and $u$ is a unipotent
element of $Z_G(s)^0$ whose $E_7$ component is of type $E_7(a_2)$ (notation as in \cite{\WEU, 4.3}) and whose
$A_1$-component has a single Jordan block of size $2$ in the standard representation of $A_1$. On the other 
hand, $(J_{11}J_5)_{D_8}$ is the conjugacy class of $su=us$ where $s$ is a semisimple element with 
$Z_G(s)^0$ of type $D_8$ and $u$ is a unipotent element of $Z_G(s)^0$ with Jordan blocks of sizes $11,5$ in 
the standard representation of $D_8$. 

{\it Type $E_8$.}
$$\align&8;\Ph_{30}; (E_8)_{E_8}, (E_7J_2)_{E_7A_1}, (E_6J_3)_{E_6A_2},(J_9J_1J_4)_{D_5A_3}, 
(J_5J_5)_{A_4A_4},\\&
(J_6J_3J_2)_{A_5A_2A_1}, (J_9)_{A_8},(J_8J_2)_{A_7A_1},(J_{15}J_1)_{D_8},\endalign$$
$$10;\Ph_{24};(E_8(a_1))_{E_8}, (E_7(a_1)J_2)_{E_7A_1}, (E_6(a_1)J_3)_{E_6A_2},(J_7J_3J_4)_{D_5A_3}, 
(J_{13}J_3)_{D_8}.$$
$$12;\Ph_{20}; (E_8(a_2))_{E_8},(E_7(a_2)J_2)_{E_7A_1},(J_{11}J_5)_{D_8},$$
$$14;\Ph_6\Ph_{18}; (E_7A_1)_{E_8}, (E_7(a_3)J_2)_{E_7A_1},(J_9J_7)_{D_8},$$
$$16;\Ph_{15}; (D_8)_{E_8},(E_7(a_4)J_2)_{E_7A_1},$$
$$18;\Ph_2^2\Ph_{14}; (E_7(a_1)A_1)_{E_8},$$   
$$20;\Ph_{12}^2;(D_8(a_1))_{E_8},(J_7J_5J_3J_1)_{D_8},$$
$$22;\Ph_6^2\Ph_{12}; (E_7(a_2)A_1)_{E_8},(E_7(a_5)J_2)_{E_7A_1},$$
$$24;\Ph_{10}^2; (A_8)_{E_8},$$
$$28;\Ph_3\Ph_9;(D_8(a_3))_{E_8},$$
$$40;\Ph_6^4;(2A_4)_{E_8}.$$

{\it Type $E_7$.}
$$7;\Ph_2\Ph_{18};(E_7)_{E_7},(J_{11}J_1J_2)_{D_6A_1},(J_6J_3)_{A_5A_2},(J_4J_4J_2)_{A_3A_3A_1},(J_8)_{A_7},
$$
$$9;\Ph_2\Ph_{14};(E_7(a_1))_{E_7}, ((J_9J_3)A_1)_{D_6A_1},$$
$$11;\Ph_2\Ph_6\Ph_{12};(E_7(a_2))_{E_7}, (J_7J_5J_2)_{D_6A_1},$$
$$13;\Ph_2\Ph_6\Ph_{10};(D_6A_1)_{E_7},$$
$$17;\Ph_2\Ph_4\Ph_8; (D_6(a_1)A_1)_{E_7},$$
$$21;\Ph_2\Ph_6^3;(D_6(a_2)A_1)_{E_7}.$$

{\it Type $E_6$.}
$$6;\Ph_3\Ph_{12};(E_6)_{E_6},(J_6J_2)_{A_5A_1},(J_3J_3J_3)_{A_2A_2A_2},$$
$$8;\Ph_9;(E_6(a_1))_{E_6},$$
$$12;\Ph_3\Ph_6^2; (A_5A_1)_{E_6}.$$

{\it Type $F_4$.}
$$4;\Ph_{12}; (F_4)_{F_4},(J_6J_2)_{C_3A_1},(J_3J_3)_{A_2A_2},(J_4J_2)_{A_3A_1}, (J_9)_{B_4},$$
$$6;\Ph_8; (F_4(a_1))_{F_4},(J_4J_2J_2)_{C_3A_1},$$
$$8;\Ph_6^2;(F_4(a_2))_{F_4},(J_5J_3J_1)_{B_4},$$
$$12;\Ph_4^2;(F_4(a_3))_{F_4}.$$

{\it Type $G_2$.}
$$2;\Ph_6;(G_2)_{G_2},(J_3)_{A_2},(J_2J_2)_{A_1A_1},$$
$$4;\Ph_3;(G_2(a_1))_{G_2}.$$

\enddocument